\documentstyle[amssymb,amsmath,eucal]{article}

\begin{document}

\def\R {{\Bbb R }}
\def\C {{\Bbb C }}
\def\KK{{\Bbb K}}
\def\HH{{\Bbb H}}

\def\K{{\bold{K}}}
\def\G{{\bold{G}}}

\def\W{{\rm W}}

\def\SS{{\cal S}}

\def\U{{\rm U}}

\def\H{{\frak B}}
 \def\h{{\cal B}}
 \def\T{{\cal T}}
\def\const{{\rm const}}
\def\B{{\rm B}}

\def\W{{\rm W}}
\def\SW{{\rm SW}}
\def\Pos{{\rm Pos}}
\def\Pol{{\rm Pol}}
\def\Fl{{{\rm Fl}_{p,q}}}
\def\Gr{{\rm Gr}}
\def\IGr{{\rm IGr}}
\def\P{{\bold P}}

\def\OO{{\rm O}}
\def\SO{{\rm SO}}
\def\GL{{\rm GL}}
\def\SL{{\rm SL}}
\def\SU{{\rm SU}}
\def\Sp{{\rm Sp}}
\def\SOS{\SO^*}

\def\BB{{\cal B}}
\def\V{{\cal V}}

\def\Mat{{\rm Mat}}

\def\Subsection{}
\def\Subsections{}

\newcommand{\rk}{\mathop{\rm rk}\nolimits}

\def\kvadrat{\hfill~$\boxtimes$}

\def\ov{\overline}
`\def\phi{\varphi}
\def\epsilon{\varepsilon}
\def\kappa{\varkappa}
\def\le{\leqslant}
\def\ge{\geqslant}

\newcommand{\Res}{\mathop{\rm Res}\nolimits}

\renewcommand{\Re}{\mathop{\rm Re}\nolimits}
\renewcommand{\Im}{\mathop{\rm Im}\nolimits}
\newcommand{\tr}{\mathop{\rm tr}\nolimits}

\def\dual{{\widehat \G_{\rm sph}}}

\def\L{{\cal L}}
\def\M{{\bold M}}

\vspace{22pt}

\begin{center}

{\Large\bf Matrix balls, radial analysis of Berezin kernels,
and hypergeometric determinants}

\vspace{22pt}

{\large\sc
Yurii A. Neretin}\footnote%
{supported by the
grants RFBR 98-01-00303 and
NWO 047-008-009}

\end{center}

\vspace{22pt}

 The subject of this paper is a natural
one-parametric family of Hilbert spaces  $V_\alpha$
and representations interpolating $L^2$ on a Riemannian
noncompact symmetric space $G/K$ and $L^2$ on the
dual Riemannian compact symmetric space
$G_{{\rm comp}}/K$. The spaces $V_\alpha$ consist of
holomorphic functions on some symmetric space
$\widetilde G/\widetilde K\supset G/K$, the limit of
$V_\alpha$ as $\alpha\to+\infty$ is $L^2(G/K)$,
and the limit of $V_\alpha$ as $\alpha$  tends
to $-\infty$ taking integer values is   $L^2(G_{{\rm comp}}/K)$.
Even the case in which   $G/K$
is the Lobachevskii plane $\U(1,1)/\U(1)\times\U(1)$
and
 $G_{{\rm comp}}/K$
 is the two-dimensional sphere $\U(2)/\U(1)\times\U(1)$,
is nontrivial and is still not satisfactory clear today  (see \cite{Ner8}).

 More formally, we consider the problem
of the restriction of a unitary highest weight representation
of
a semisimple group $\widetilde G$
 to a symmetric subgroup (see  4.7 below).
  This formulation
 is short and simple
 but it does not explain why this restriction problem
 differs from a nondenumerable collection
  of other restriction problems.

The objects of this paper first appeared in the
 short
note of Berezin \cite{Ber2} in 1978. Berezin perished two
years later and a complete text of his work never
was published.
For the second time,
 these  objects occurred in a joint work
of G.I.Olshanskii and author, but it was published only
partially in two short notes \cite{Ner} and \cite{Ols1}.

As a result,
 this problem was "lost"  and  it became
 visible
only in 1994-95 in \cite{UU} and \cite{NO}.
Some other references
in the last 5 years are \cite{vD}, \cite{Hil},
\cite{Ner3}--\cite{Ner8},  \cite{OO}, \cite{OZ},
\cite{Zha} (these papers can also  be the source
for many other references).

 The present paper has three purposes.

1) We intend to  survey
 phenomena  arising in
the analysis of the Berezin kernels.
We consider a special case
$$G=\U(p,q)$$
In this case, there exist some specific
tools that allow to avoid main difficulties
existing in a general case. For instance, the complete
Plancherel formula%
\footnote{In the case of Hermitian
symmetric spaces $G/K$ the Plancherel formula
for large values of the parameter $\alpha$
(see below) was announced
in the Berezin work \cite{Ber2}, proof was published
by Upmeier and Unterberger in 1994 \cite{UU};
for rank 1 case it was obtained by
van Dijk and Hille \cite{vD}, for general
situation it was obtained in \cite{Ner6},\cite{Ner8}}
 can be proved in a very simple
way\footnote{Partially this was done by Hille \cite{Hil}}.

2) For the case $G=\U(p,q)$, we can also obtain some
results
that today are not achieved in the general
situation. The  main new result
 is an explicit construction  of
 a unitary intertwining operator
 from $L^2\bigl(\U(p,q)/\U(q)\times\U(p)\bigr)$
 to the Berezin deformation of $L^2$.

3) Many multivariate special functions appear
 in a natural way in the harmonic analysis of the Berezin
 kernels (in particular, multivariate continuous dual Hahn,
 dual Hahn, Meixner--Pollachek, Krawtchouk, Laguerre,
 Jacobi orthogonal polynomials, Jack polynomials,
Jack zonal functions,
  matrix $\Gamma$-function and $\B$-function,
 matrix Bessel functions, Heckman--Opdam multivariate
 hypergeometric functions). Also, it seems to me
 that the analysis of the Berezin kernels leads
 to some "new" special functions (to the $\Lambda$-function,
 see Section 11, and to generalizations of Gross--Richards
 kernels, see section 10).
 The relatively simple picture for $\U(p,q)$
allows to touch some of these objects in a simple way.

It seems to me that   the subject of this paper
 is elementary, this is
some topic in analysis of the matrix variable.
For this reason, I try to follow a
simple approach to noncommutative
 harmonic analysis
in the spirit of \cite{GN}, \cite{Hua}, \cite{Zhe}
whenever I can.

I am very grateful to Grigory Olshanski,
Bent \O rsted and Vladimir Molchanov for meaningful
discussion. I thanks the administration of
the Erwin Schr\"odinger Institute for Mathematical
Physics, where this text was prepared, for their hospitality.

\begin{center}

{\large\bf Contents}

\end{center}

{\sc

\noindent
  1. Preliminaries: positive definite kernels

\noindent
  2.   Groups $\U(p,q)$ and matrix balls

\noindent
  3.   Spaces of holomorphic functions
on matrix balls. Berezin scale

\noindent
   4.   Kernel representations  and spaces $\V_\alpha$

\noindent
    5.  Index hypergeometric transform: preliminaries

\noindent
    5.  Helgason transform and spherical transform

\noindent
    7.  Plancherel formula for kernel representations

\noindent
    8.  Boundary behavior of
holomorphic functions and separation of spectra

\noindent
    9.  Interpolation between
  $L^2\bigl(\U(p,q)/\U(p)\times\U(q)\bigr)$
  and
  $L^2\bigl(\U(p+q)/\U(p)\times\U(q)\bigr)$.

   Pickrell formula

\noindent
    10.  Radial part of the spaces
$\V_\alpha$ and Gross--Richards kernels

 \noindent
   11.  \O rsted problem. Identification
of kernel representations and $L^2(\G/\K)$

\smallskip

\noindent
\quad   Addendum. Pseudoriemannian symmetric spaces, Berezin forms,
   and some problems of non $L^2$-harmonic analysis
}

 \bigskip

\newcounter{sec}
 \renewcommand{\theequation}{\arabic{sec}.\arabic{equation}}

\newcounter{fact}
\def\fact{\addtocounter{fact}{1}{\sc  \arabic{sec}.\arabic{fact}}}

\newcounter{punkt}
\def\punct{\addtocounter{punkt}{1}{\bf  \arabic{sec}.\arabic{punkt}}}

 \def\sect{\addtocounter{sec}{1}
       \setcounter{equation}{0}
       \setcounter{fact}{0}
        \setcounter{punkt}{0}
        {\large\bf \arabic{sec}. }}

{\large\bf \sect Preliminaries: positive definite kernels}

\nopagebreak

\medskip

Positive definite kernels machinery is a usual tool
for work upon
Hilbert spaces (\cite{Sch}, \cite{Kre}).
 This section contains simple
general facts concerning this subject.

\smallskip

{\bf \punct. Positive definite kernels.}
Let $X$ be a set, let $H$ be a Hilbert space, and
let $\langle\cdot,\cdot\rangle$ denote the scalar product in $H$.
 Consider
a map $x\mapsto v_x$ from $X$ to $H$. We define the function
$L(x,y)$ on $X\times X$  by
$$L(x,y)=  \langle v_x,v_y\rangle $$

Obviously,
$
L(x,y)= \overline{L(y,x)}
$ and
for any $x_1,\dots,x_n\in X$ the matrix
\begin{equation}
\begin{pmatrix}
L(x_1,x_1) &\dots & L(x_1,x_n)\\
\vdots &\ddots &\vdots \\
L(x_n,x_1) &\dots & L(x_n,x_n)\label{1.1}
\end{pmatrix}
\end{equation}
is positive semidefinite\footnote{An $n\times n$
Hermitian matrix
$Q=\{q_{ij}\}$ is called {\it positive semidefinite} if
$\sum q_{ij}w_i\ov w_j\ge0$ for any $w_1,\dots,w_n\in\C$.}.

 A function $L(x,y)$ on the set $X\times X$ is called
{\it a positive definite kernel} on $X$ if
$L(x,y)= \overline{L(y,x)}$ and for  all
 $x_1,\dots,x_n\in X$ the matrix (\ref{1.1})  is positive semidefinite.

\smallskip

  Nontrivial examples
appear in  Sections 3--4,

\smallskip

{\sc Theorem \fact.} \cite{Sch} {\it Let $L(x,y)$
be a positive definite kernel on $X$. Then}

a) {\it There exists a Hilbert space $H$ and a system of vectors
$v_x\in H$ enumerated by  $x\in X$ such that
$L(x,y)=  \langle v_x,v_y\rangle $ and the linear span
of the vectors $v_x$ is dense in $H$.}

b) {\it  A space $H$ is unique in the following sense.
Let $H'$ be another Hilbert space and let $v'_x\in H'$
be another system of vectors satisfying the same conditions.
Then there
exists a unique unitary operator $U:H\to H'$ such that
$Uv_x=v'_x$   for all $x$.}

\smallskip

{\sc Proof.} Let $v_x$, where $x\in X$,
 be formal symbols. Consider
the linear space $\widetilde H$  consisting of all the formal
finite linear
combinations
$\sum_{k=1}^n c_k v_{x_k}$, where $c_k\in \C$.
 We define a scalar product
in    $\widetilde H$ by
$$ \langle\sum_{k=1}^n c_k v_{x_k},
     \sum_{l=1}^{m} c'_l v_{y_l}\rangle
     =\sum_{k=1,l=1}^{n,m}c_k\overline c_l' L(x_k, y_l)$$
The space $\widetilde H$ is a pre-Hilbert space and $H$ is
 the Hilbert space
associated with $\widetilde H$.   \kvadrat

\smallskip

We denote\footnote{
in honor of S.Bergman, V.Bargmann and F.A.Berezin}
by $\H[L]=\H[L;X]$ the space $H$
equipped with the distinguished system of  vectors
$v_x$.
The set of  vectors $v_x\in \H[L]$ is called
  {\it a supercomplete system},
or {\it an overfilled basis}, or {\it a system of coherent states}.

\smallskip

{\sc Remark.} A supercomplete system is nothing but
a system of vectors with explicitly known pairwise
scalar products. Nevertheless the knowledge of these data
can be an effective tool
 for  workin  in  a given Hilbert space.
 \hfill $\square$

\smallskip

If $X$ is a separable metric space and $L(x,y)$
 is continuous on $X\times X$, then
the Hilbert space $\H[L]$ is separable and the map
$x\mapsto v_x$ is continuous.
In all natural cases, these conditions are satisfied.
For certain technical reasons,
the general definition given above is more convenient.

\smallskip

{\bf \punct. Functional realization $\H^\circ[L]$.}
Consider the space $\H[L]$ defined by a positive definite kernel
$L$. To each $h\in \H[L]$ we assign a function $f_h(x)$ on $X$ by
$$f_h(x)=\langle h, v_x\rangle_{\H[L]}$$
The linear span of the vectors $v_x$
is dense in $\H[L]$ and hence the map
$h\mapsto f_h$ is an embedding of $\H[L]$ to the space
of functions on $X$. We denote  by $\H^\circ[L]$
the image of this embedding.
 A scalar product in   $\H^\circ[L]$
is defined by
\begin{equation}
\langle f_h, f_q\rangle_{\H^\circ[L]}
: =\langle h,q\rangle_{\H[L]}\label{1.2}
\end{equation}
The function $\theta_a$
 corresponding to an element
$v_a$ of the supercomplete system is given by
$$\theta_a(x)=L(a,x)$$

{\sc Proposition \fact.} (Reproducing property)
{\it  For each $f\in\H^\circ[L]$ }
      \begin{equation}
\langle f, \theta_a\rangle_{\H^\circ[L]}=f(a)
\label{1.3}
\end{equation}

{\sc Proof.}
$$\langle f_h,\theta_a\rangle_{\H^\circ[L]} =
\langle h,v_a\rangle_{\H[L]}=f_h(a)$$

{\bf \punct. Reconstruction of the kernel $L$
from $\H^\circ[L]$.}
Let $\H^\circ$ be some Hilbert space  consisting
 of functions\footnote{The space $L^2$ does not
 consist of functions!}
  on $X$, and let
 the linear functional
$u_x(f):=f(x)$ be continuous on $\H^\circ$ for each $x\in X$ .
 Then (see \cite{RS}, Theorem 2.4) for any $x\in X$
 there exists a unique function $\theta_x\in \H^\circ$ such that
 $$\langle f, \theta_x\rangle=f(x)$$
 We define a positive definite kernel $L(x,y)$ on $X$ by
 $$L(x,y):=\langle \theta_x,\theta_y \rangle_{\H^\circ}=
 \theta_x(y)=\ov{\theta_y(x)}$$
 for all $f\in \B^\circ$.
 Then the space $\H^\circ$ coincides with the space $\H^\circ[L]$
 defined by the positive definite kernel $L$.

 The kernel $L$ is called the {\it reproducing kernel}
 of the space $\H^\circ$.

 \smallskip

 {\bf \punct. Convergence in $\H^\circ[L]$.}

 \nopagebreak

 {\sc Lemma \fact.} (see, for instance, \cite{Ner8})
  {\it Assume $X$ is a complete metric space,
 and the kernel $L$ is continuous. If a sequence
 $f_j\in \H^\circ[L]$ converges to $f$ in $\H^\circ[L]$, then
 it converges uniformly on compact sets in $X$.}

 \smallskip

 {\sc Corollary \fact.} {\it Any element of $\H^\circ[L]$
 can be approximated by finite sums $\sum c_k L(a_k,x)$
 in the topology of uniform convergence on compact sets.}

 \smallskip

 {\bf \punct. Hilbert spaces of holomorphic functions.}
 Let $\Omega$ be a bounded open domain in $\C^n$.
Assume a positive definite kernel $L(z,u)$
on $\Omega$ is antiholomorphic
 in $z$ and holomorphic in $u$. Then, by Corollary 1.4,
 elements of the Hilbert space $\H^\circ[L]$
  are holomorphic functions\footnote{Recall
  the {\it Weierstrass
  theorem}.  Let $f_n$ be holomorphic and
  $f_n$ converge to $f$ uniformly on compact sets.
  Then partial derivatives
   $\frac{\partial}{\partial z_\alpha}f_n$ converge to
   $\frac{\partial}{\partial z_\alpha}f$
  uniformly on compact sets. In particular, $f$ is holomorphic.
   }.

 Let $\zeta(z)$
 be a  continuous function on
 $\Omega$ and let $\zeta(z)>0$ for any $z\in \Omega$.
   Denote by ${\cal B}(\Omega, \zeta)$
 the space of holomorphic functions on $\Omega$ satisfying the
 condition
 $$\int_\Omega |f(z)|^2 \zeta(z)\, \{dz\}<\infty$$
 where
 $$\{dz\}:=\prod_{j=1}^n d(\Re z_j)\,\prod_{j=1}^n d(\Im z_j)$$
  denotes the Lebesgue measure on $\Omega$.
 Consider the $L^2$ scalar product
 \begin{equation}
 \langle f,g\rangle=
 \int_\Omega f(z)\overline{g(z)} \zeta(z)\, \{dz\}
 \label{1.4}
 \end{equation}
 in the space   ${\cal B}(\Omega, \zeta)$.

 \smallskip

  {\sc Theorem \fact} a) {\it
The space    ${\cal B}(\Omega, \zeta)$
 is closed in $L^2(\Omega, \zeta)$.}

 b){\it For any $u\in\Omega$
  the linear functional $f\mapsto f(u)$
 is continuous on    ${\cal B}(\Omega, \zeta)$. }

 \smallskip

 {\sc  Proof.}  b) Let
  $\Omega$ be a polydisk
 $|z_j|<r_j$, $\zeta(z)=1$ and  $u=0$. Then
 $f(0)$ coincides (up to a factor) with $\langle f,1\rangle$.
 Thus in this case the statement is obvious. Consider
 an arbitrary point $w\in \Omega$ and a small polydisk
 $D\subset\Omega$
 with center in $w$. Then
 $$\int_\Omega |f(z)|^2\zeta(z) \,\{dz\}
 \ge \min_{z\in D} \zeta(z)
 \int_D |f(z)|^2 \{dz\}$$
 Hence the convergence $f_n\to f$
  in    ${\cal B}(\Omega, \zeta)$ implies
 the convergence  $f_n\to f$ in $L^2(D)$.
 Thus $f_n(w)$ converges to $f(w)$. This implies b)

 The assertion of  a) is a consequence of b) and Lemma 1.3.
  \kvadrat

  \smallskip

  The spaces ${\cal B}(\Omega,\tau)$ are called
   {\it weight spaces of holomorphic functions.}

  \smallskip

  {\bf\punct. Another functional realization $\H^\star[L]$.}
  Let $X$ be a smooth manifold. Let a
positive definite kernel $L$ be
  $C^\infty$-smooth. Let ${\cal E}'(X)$ be
  the space of compactly supported distributions on $X$.
  We define a scalar product in  ${\cal E}'(X)$
  by
  $$
  \langle \chi_1,\chi_2\rangle=
  \{ L , \chi_1\otimes\ov\chi_2\}
  $$
  where $\{\cdot,\cdot\}$ denotes the pairing of
   smooth functions and distributions.
  We define the space $\H^\star[L]$ as the
 Hilbert
  space associated with the pre-Hilbert space
    ${\cal E}'(X)$.

    The canonical unitary operator $J:\H[L]\to\H^\star[L]$
    is defined by the condition:
    $Jv_x$ is the $\delta$-function supported by $x$.

    {\sc Remark.} In general,
    distributions do not represent all
   the elements of $\H^\star[L]$.

{\bf \punct. Some operations with positive definite kernels.}

{\sc Proposition \fact.} {\it Let $L_1$, $L_2$
be positive definite kernels on $X$. Then
 $L_1 L_2$ and $L_1+L_2$ are positive definite kernels.}

 {\sc Proof.} Let $v_x$ be the supercomplete system in $\H[L_1]$
 and $w_x$ be the supercomplete system in   $\H[L_2]$.
 Consider the space $\H[L_1]\otimes\H[L_2]$
 with the supercomplete system $v_x\otimes w_x$
 and the space
               $\H[L_1]\oplus\H[L_2]$
 with the supercomplete system $v_x\oplus w_x$.
 Both statements are now  obvious.\kvadrat

 \smallskip

 {\sc Proposition \fact.} {\it Let $L(x,y)$ be a positive
 definite kernel on $X$.
  Then for an arbitrary function $\gamma(x)$ on $X$

  {\rm a)} the
  kernel $M(x,y):=L(x,y)\gamma(x)\ov{\gamma(y)}$
   is positive definite.

   {\rm b)} the operator $\H^\circ[L]\to   \H^\circ[M]$
   given by $f(x)\mapsto f(x)\gamma(x)$ is unitary.}

  \smallskip

  {\sc Proof.} Consider the vectors $\gamma(x)v_x\in\H[L]$.\kvadrat

  \bigskip

{\large\bf \sect Groups $\U(p,q)$ and matrix balls}

\medskip

       Assume that $p\le q$.

{\bf \punct. Groups $\U(p,q)$.}
 Consider the space $\C^p\oplus\C^q$
 equipped with the indefinite Hermitian form $Q$ defined by
 the $(p+q)\times (p+q)$ matrix
 $\bigl(\begin{smallmatrix} 1&0\\0&-1\end{smallmatrix}\bigr)$.
 We denote by $e_1^+,\dots,e_p^+,e_1^-,\dots e_q^-$
 the standard basis in $\C^p\oplus\C^q$.
   The {\it pseudounitary group}\footnote{We also denote by $\U(n)$
   the group of unitary $n\times n$ matrices.}
    $\G=\U(p,q)$
is the group of all linear operators  in $\C^p\oplus\C^q$
preserving the form $Q$. In other words,  the group
$\U(p,q)$ consists of all $(p+q)\times(p+q)$ matrices
 satisfying the condition
\begin{equation}
\begin{pmatrix}a&b\\ c&d\end{pmatrix}
\begin{pmatrix}1&0\\0&-1\end{pmatrix}
\begin{pmatrix}a&b\\ c&d\end{pmatrix}^*=
\begin{pmatrix}1&0\\0&-1\end{pmatrix}
\label{2.1}
\end{equation}

{\bf\punct. Another realization of $\U(q,q)$.}  Assume $p= q$.
Consider the new basis
\begin{align}
w_1^+:=\tfrac1{\sqrt 2} (e_1^++e_1^-),\dots,
 w_q^+:=\tfrac1{\sqrt 2}  (e_q^++e_q^-),
 \label{2.2}
 \\
   w_1^-:=\tfrac1{\sqrt 2} (e_1^+-e_1^-),\dots,
   w_q^-=\tfrac1{\sqrt 2}  (e_q^+-e_q^-)
   \label{2.3}
   \end{align}
 in $\C^q\oplus \C^q$. The matrix of the
 Hermitian form $Q$ in this basis is
  $\bigl(\begin{smallmatrix} 0&1\\1&0\end{smallmatrix}\bigr)$.
  Thus the group $\U(q,q)$ can be represented as the group of
all  $(q + q)\times (q+q)$ matrices
  $\bigl(\begin{smallmatrix} A&B\\C&D\end{smallmatrix}\bigr)$
  satisfying
  \begin{equation}
  \begin{pmatrix}A&B\\ C&D\end{pmatrix}
\begin{pmatrix}0&1\\1&0\end{pmatrix}
\begin{pmatrix}A&B\\ C&D\end{pmatrix}^*=
\begin{pmatrix}0&1\\1&0\end{pmatrix}
\label{2.4}
\end{equation}

 {\bf\punct. Another realization of $\U(p,q)$.}
Let us realize the groups $\U(p,q)$ and $\U(q,q)$ as in 2.1.
 Consider the natural embedding
of $\U(p,q)$ to $\U(q,q)$ given by
$$h\mapsto
 \begin{pmatrix} 1_{q-p}&0\\0&h\end{pmatrix}$$
 where $1_{q-p}$ denotes the unit matrix
 of the size $q-p$.
 Consider the realization of  $\U(q,q)$  described in
  \Subsection 2.2. Then $\U(p,q)$ becomes the group
 of
 $((q-p)+p+(q-p)+p)\times ((q-p)+p+(q-p)+p)$
 block matrices $g$ satisfying the conditions
 \begin{equation}
 g\begin{pmatrix} 0&0&1&0\\
                   0&0&0&1\\
                   1&0&0&0\\
                   0&1&0&0\end{pmatrix} g^*
                   =
                \begin{pmatrix} 0&0&1&0\\
                   0&0&0&1\\
                   1&0&0&0\\
                   0&1&0&0\end{pmatrix};\qquad
 g\begin{pmatrix} v\\0\\v\\0\end{pmatrix}=
         \begin{pmatrix} v\\0\\v\\0\end{pmatrix}
\label{2.5}
\end{equation}
for all $v\in \C^{q-p}$.
Sometimes this  model is useful.

\smallskip

{\bf \punct. The symmetric space
 $\G/\K=\U(p,q)/\U(p)\times\U(q)$.}
We say that a $p$-dimensional subspace
 $R\subset\C^p\oplus\C^q$
is {\it positive}, if the Hermitian form $Q$ is
 positive definite on
$R$. Denote by ${\Gr}_{p,q}^+$ the space of
all positive $p$-dimensional subspaces in  $\C^p\oplus\C^q$.
The group $\U(p,q)$ acts on $\C^p\oplus\C^q$ and hence
it acts on ${\Gr}^+_{p,q}$. Obviously
(by the Witt theorem),
this action is transitive. The stabilizer of the subspace
$\C^p\oplus 0\in {\Gr}^+_{p,q}$ consists
of matrices having the form
\begin{equation}
\begin{pmatrix} a&0\\ 0&d \end{pmatrix};\qquad\qquad
a\in\U(p),\,\, b\in \U(q)
\label{2.6}
\end{equation}
Hence the space ${\Gr}^+_{p,q}$
is the homogeneous space
$$\G/\K=\U(p,q)/\U(p)\times\U(q)$$
{\it In this paper we fix the notation}
$$\G=\U(p,q);\qquad \K=\U(p)\times\U(q)$$

\smallskip

 {\bf\punct. Cartan matrix balls.} Assume $p\le q$. A matrix
 ball
 $\B_{p,q}$ is the space of all complex $p\times q$ matrices
 with norm\footnote{The term norm in this paper
 means the norm
 of an operator
 from Euclidean space $\C^p$ to Euclidean space $\C^q$.
 The notation $A>0$ means that the operator $A$ is positive
 i.e., $\langle Ax,x\rangle> 0$ for all $x\ne 0$. The notation
 $A>B$ means $A-B>0$.}
  less than 1.
 Also, $z\in\B_{p,q}$ iff $zz^*<1$.

 For $z\in\B_{p,q}$ we define the subspace
 $Graph(z)\subset\C^p\oplus\C^q$
 consisting of all the vectors of the form $(v,vz)$,
 where $v\in\C^p$. It can easily be checked that the map
 $z\mapsto Graph(z)$ is a bijection
 $\B_{p,q}\to\Gr^+_{p,q}$. Thus
  the space $\Gr^+_{p,q}\simeq\G/\K$
 is parametrized by points of the
 matrix ball $\B_{p,q}$.
  The action of the group
 $\G$ in these coordinates is given
by the {\it linear-fractional transformations}
\begin{equation}
z\mapsto z^{[g]}:=(a+zc)^{-1}(b+zd)
\label{2.7}
\end{equation}
where $g=\bigl(\begin{smallmatrix}a&b\\c&d\end{smallmatrix}\bigr)
\in\G$.

 The complex Jacobian of the transformation (\ref{2.7}) is
  \begin{equation}
  \det(a+zc)^{-p-q}
  \det
  \begin{pmatrix}a&b\\ c&d\end{pmatrix}
\label{2.8}
  \end{equation}

{\sc Lemma \fact.}
 {\it The $\G$-invariant measure on $\B_{p,q}$ is given by
  \begin{equation}
  \det(1-zz^*)^{-p-q} \{dz\}
\label{2.9}
  \end{equation}
  where $\{dz\}$  is the Lebesgue measure on $\B_{p,q}$.}

{\sc Proof.}
  This follows from
  (\ref{2.8}) and the simple identity
  \begin{equation}
  1-z^{[g]}\bigl(u^{[g]}\bigr)^*=
     (a+zc)^{-1}(1-zu^*)(a^*+c^*u^*)^{-1}
  \label{2.10}
  \end{equation}

  \smallskip

  {\bf \punct. Structure of the boundary
  of the matrix ball.} Denote by $\overline \B_{p,q}$
the closure of $\B_{p,q}$ in $\C^{pq}$,
 i.e the set of matrices with norm
$\le 1$.  Denote by $M_k$ the set of all matrices
$z\in\overline \B_{p,q}$ satisfying the condition
    $$\rk (1-zz^*)=k;\qquad k=0,\dots, p-1$$
  The following statement is trivial 

  {\sc Lemma \fact.}
  {\it   The sets $M_k$ are exactly the orbits
  of $\U(p,q)$ on the boundary of $\B_{p,q}$.}

  {\sc Remark.}  The orbit $M_0$
  is the Shilov boundary of $\B_{p,q}$, see \cite{FK}, X.5.
  \hfill $\square$

  \smallskip

  {\bf\punct. Siegel realizations of $\G/\K$: matrix wedges.}
  The Cartan realization of $\G/\K$ will be basic for us.
  Nevertheless in some places
  we shall need of the  Siegel realizations.%
  \footnote{The action of the compact subgroup $\K\subset\G$
   has the simplest form
   in the Cartan realization; if we want to write
    formulas related to the parabolic subgroup
    (see  Subsection 6.2 below),
   then the Siegel realizations are more convenient.}

  Assume $p=q$.
    Denote by $\W_q$ the space of  all
  $q\times q$ complex matrices $T$ satisfying the condition
  $$ T+T^*>0$$

  Consider the basis (\ref{2.2})--(\ref{2.3}),
   denote by $Z_+$ the subspace in $\C^q\oplus\C^q$
   spanned
  by $w_j^+$ and denote by $Z_-$ the subspace spanned by
  $w_j^-$. For $T\in \W_q$ consider the operator
   $T:Z_+\to Z_-$.
  It can easily be checked that
   the map $T\mapsto Graph (T)$ is a bijection
  $\W_q\to {\Gr}^+_{q,q}$

  The action of the group $\G=\U(q,q)$ on the wedge
    $\W_q$  is given by
    \begin{equation}
    T\mapsto T^{[g]}:=(A+ TC )^{-1}(B+TD)
    \label{2.11}
    \end{equation}
    where the matrix
    $g=
    \bigl(\begin{smallmatrix}A&B\\C&D\end{smallmatrix}\bigr)$
    satisfies (\ref{2.4}).

 The spaces $\B_{q,q}$  and $\W_q$ are identified by
 the {\it Cayley transform}
 $$ T= -1+  (1+ z)^{-1}$$

 {\bf \punct. Siegel realization in the case $q\ne p$:
 sections of wedges.}
 We define the space $\SW_{p,q}$ as the subset
 in $\W_q$ consisting of all block
 $((q-p)+p)\times ((q-p)+p)$-matrices
of the form
 \begin{equation}
 S= \begin{pmatrix} 1&0\\2K & L \end{pmatrix}\in \W_q
 \label{2.12}
 \end{equation}
 The group $\U(p,q)$ acts on the space $\SW_{p,q}$ by the
 same formula (\ref{2.11}) for a matrix $g$ satisfying
 equations (\ref{2.5}), for details see \cite{Ner5}.

 \smallskip

 {\bf\punct. Radial part of the Lebesgue measure.}
   The subgroup
 $\K=\U(p)\times\U(q)\subset\G$ acts on
 $\B_{p,q}$ by the transformations $z\mapsto a^{-1}zd$
 (see (\ref{2.6}), (\ref{2.7})). Obviously,
 any element  $z\in\B_{p,q}$
 can be reduced by such transformations
 to the form
 $$\begin{pmatrix}
 \lambda_1& 0 & \dots &0& 0&\dots &0 \\
 0& \lambda_2 & \dots &0& 0&\dots &0 \\
\vdots & \vdots &\ddots &\vdots  &\vdots&\vdots&\vdots\\
 0& 0 &\dots & \lambda_p & 0&\dots & 0
 \end{pmatrix}  $$
 where
 \begin{equation}
 1>\lambda_1\ge\lambda_2\ge \dots \ge \lambda_p\ge 0
 \label{2.13}
 \end{equation}
 are the eigenvalues of $\sqrt{zz^*}$.
 We denote the simplex (\ref{2.13}) by $\Lambda_p$.

 Consider the
 map $\pi:\B_{p,q}\to \Lambda_p$ taking each
 $z\in\B_{p,q}$ to the collection of
 the eigenvalues of $\sqrt{zz^*}$.
The image of the Lebesgue measure on $\B_{p,q}$ with respect
 to the map $\pi$ is the measure on $\Lambda_p$ given by
 $$
 \const\cdot
 \prod_{i\le j} (\lambda_i^2-\lambda_j^2)^2 \prod_j \lambda_j^{2(q-p)+1}
 \prod_j d\lambda_j
 $$
 (see \cite{Hel1}, X.1),
 see also numerous calculations of this type
 in \cite{Hua}, Chapter 3).

 \smallskip

 {\it Here and below the symbol '$\const$' denotes a constant
 depending only on $p,q$.}

 \smallskip

  Thus the  image of the $\G$-invariant measure (\ref{2.9})
 is
 \begin{equation}
 \const\cdot
 \prod_{1\le i< j\le p}
  (\lambda_i^2-\lambda_j^2)^2 \prod_{j=1}^p \lambda_j^{2(q-p)+1}
 (1-\lambda_j^2)^{-p-q}
 \prod_{j=1}^p d\lambda_j
 \label{2.14}
 \end{equation}

 It will be convenient for us to define
new coordinates
\begin{equation}
x_j= {\lambda_j^2}/{(1-\lambda_j^2)}
\label{2.15}
\end{equation}
Then the simplex $\Lambda_p$ transforms
to the simplicial cone ${\cal X}_p$
\begin{equation}
{\cal X}_p:x_1\ge x_2\ge\dots\ge x_p\ge 0
\label{2.16}
\end{equation}
The measure (\ref{2.14})
 in the coordinates $x_j$ is of the form
\begin{equation}
\prod_{j=1}^p x_j^{q-p}
\prod_{1\le k<l\le p} (x_k-x_l)^2\prod_{j=1}^p dx_k
\label{2.17}
\end{equation}

{\it Below we consider $\K$-invariant functions on $\G/\K$
as symmetric  functions in
the variables $x_j\ge 0$
or as functions on the simplicial cone $\cal X$.}

\smallskip

{\sc Theorem \fact.} (Hua)
\begin{equation}
\int_{\B_{p,q}}  \det(1-zz^*)^{\alpha-p-q}\{dz\}=
\omega\prod_{k=1}^p
\frac{\Gamma(\alpha-q-k+1)}{\Gamma(\alpha-k+1)}
\label{2.18}
\end{equation}
 {\it where $\omega$ is the volume of $\B_{p,q}$.}

The problem is reduced
 to integrating of the function $\prod(1+x_k)^{-\alpha}$
over the measure (\ref{2.17}).
This is a special case of the Selberg integral,
see  \cite{AAR}, chapter 8.
 Hua's calculations (\cite{Hua}, chapter 2)  are interesting
and important by themself. For other ways
of calculation see  \cite{FK},
\cite{Ner5}.

\smallskip

{\small
{\bf \punct. Comments.} A {\it matrix ball}
(see \cite{Ner1}) $\B$ is

---   the set of all
$p\times q$ matrices over $\R$, $\C$ or
the quaternion algebra $\Bbb H$
such that $zz^*<1$

---  or the set of all $n\times n$ matrices
 over $\R$, $\C$, $\Bbb H$ satisfying
 $zz^*< 1$ and  the natural symmetry condition.

 The natural conditions of symmetry are

 $z=z^t$ and $z=-z^t$ over $\R$;

  $z=z^t$, $z=-z^t$,  and $z=z^*$ over $\C$;

  $z=z^*$ and $z=-z^*$ over $\HH$.

We consider the group $G$ of linear-fractional   transformations
(\ref{2.7})
 preserving $\B$. Denote by $K$ the stabilizer of the point
$0\in \B$. All {\it classical Riemannian
 noncompact symmetric spaces}
$G/K$
can be obtained in this way. The table is contained in
\cite{Ner2}, Appendix A, see also \cite{Ner5}.

Tools that are used below cannot be applied to arbitrary matrix
ball.  }

\bigskip

{\large\bf \sect Spaces of holomorphic functions
in matrix balls. Berezin scale}

\nopagebreak

\medskip

A subject of this  section is standard. Usually we
give sketches of proofs.

{\bf \punct. Berezin scale: large values of parameter.}
Let $\alpha>p+q-1$. Consider the weight space
$$H_\alpha:=
\h\bigl(\B_{p,q},\det(1-zz^*)^{\alpha-p-q}\{dz\}\bigr)$$
consisting of  holomorphic functions on $\B_{p,q}$
(see  \Subsection 1.5).
 The scalar product  in $H_\alpha$ is defined by
  \begin{equation}
  \langle f, g \rangle_\alpha=
 C(\alpha)^{-1}
 \int_{\B_{p,q}} f(z)\overline{g(z)}
  \det(1-zz^*)^{\alpha-p-q}\{dz\}
  \label{3.0}
  \end{equation}
We define the normalization constant $C(\alpha)$ by
(\ref{2.18}), then $\langle 1,1\rangle_\alpha=1$.

  Let us define the operators
  $\tau_\alpha(g)$, where
  $g=\bigl(\begin{smallmatrix}a&b\\c&d\end{smallmatrix}\bigr)
  \in\G=\U(p,q)$, in
  $H_\alpha$ by
  \begin{equation}
  \tau_\alpha(g)f(z)=f\bigl((a+zc)^{-1}(b+zd) \bigr)\det(a+zc)^{-\alpha}
  \label{3.1}
  \end{equation}
  A simple calculation
  (with applying (\ref{2.8}), (\ref{2.10}))
   shows that the operators $\tau_\alpha(g)$
  are unitary in  $H_\alpha$.
   Hence we obtain a unitary representation
  of $\G=\U(p,q)$ in $H_\alpha$.

  \smallskip

  {\sc Remark.} If $\alpha$ is  integer, then
   the expression $\det(a+zc)^{-\alpha}$ is well defined.
   Otherwise
   the function
  $$\det(a+zc)^{-\alpha}=\det(a)^{\alpha}
  \det\bigl[(1+zca^{-1})^{-\alpha}\bigr]=
   e^{\alpha(\ln\det a+2\pi ki)}
                          \det\bigl[(1+zca^{-1})^{-\alpha}\bigr]
  $$
  on $\B_{p,q}$ has an infinite
  number of holomorphic branches,
  which differ by  a constant factor. Indeed,
  (\ref{2.1}) implies $\|a\|>\|c\|$, hence $\|ca^{-1}\|<1$,
  hence the expression $[\dots]$ is well defined.
  Thus,  for a noninteger $\alpha$,
  the operators $\tau_\alpha(g)$ are defined up to a scalar
  factor $e^{2\pi i \alpha k}$, and the representation
  $\tau_\alpha$ is projective (see \cite{Kir}, 14).
We can also consider the representation $\tau_\alpha$
 as a
  representation of the universal covering
  group of the group $\U(p,q)$, see \cite{Kir}, Corollary
  to Theorem 14.3.1.       \hfill$\square$

  \smallskip

 {\sc Theorem \fact.}(Berezin, \cite{Ber1})
  {\it The reproducing kernel
 of the space $H_\alpha$ is  }
 $$K_\alpha(z,u)=\det(1- uz^*)^{-\alpha}$$

{\sc Proof.} We must find functions $\theta_a\in H_\alpha$
 such that $\langle f,\theta_a\rangle=f(a)$.
Expanding $f(z)$ into its Taylor series,
we obtain
\begin{multline*}
\frac {1}{C(\alpha)}\int_{\B_{p,q}}f(z)\det(1-zz^*)^{\alpha-p-q}\{dz\}=
\frac {1}{C(\alpha)}\int_{\B_{p,q}}\bigl(f(0)+\dots\bigr)
\det(1-zz^*)^{\alpha-p-q}\{dz\}=  \\=
\frac {1}{C(\alpha)}\int_{\B_{p,q}}f(0)
\det(1-zz^*)^{\alpha-p-q}\{dz\}=
f(0)
\end{multline*}

 Thus  $\theta_0(z)=1$.
  Let us evaluate
  $\langle \tau_{\alpha}(g)f, \theta_0\rangle_\alpha$
  in two ways. First, it equals
  $$\langle
  f(z^{[g]})\det(a+zc)^{-\alpha},\theta_0\rangle_\alpha
  =f(a^{-1}b)\det(a)^{-\alpha}$$
  The operators $\tau_\alpha(g)$ are unitary,
  hence it is equal to
$$\langle f, \tau_\alpha(g^{-1})\theta_0\rangle_\alpha=
\langle f,(a-b^*z)^{-\alpha}\rangle_\alpha$$
and we obtain an explicit expression for
 $\theta_{a^{-1}b}(z)$.\kvadrat

  \smallskip

Formula (\ref{3.0})
 defining the scalar product in $H_\alpha$
 has sense for $\alpha>p+q-1$. Nevertheless we shall see
(Theorem 3.8) that the reproducing kernel  $K_\alpha$
is positive definite for $\alpha>p-1$ and also for
$\alpha=0,1,\dots,p-1$. Subsections 3.2--3.3 contains
the preparation to the statement and proof of Theorem 3.8.

\smallskip

{\bf\punct.  Preliminaries: $\K$-harmonics in the space
of polynomials.}
Let
$$\mu_1\ge\mu_2\ge\dots \ge \mu_n\ge 0$$
be integers. The {\it Schur function} $s_\mu$ is defined by
$$s_\mu(y_1,\dots,y_n)=   s^n_\mu(y_1,\dots,y_n)
:=\frac{\det [y_i^{\mu_j+n-j}]} {\det [y_i^{n-j}]}=
         \frac{\det [y_i^{\mu_j+n-j}]}{\prod_{k<l}(y_k-y_l)}
$$
where $y_i\in\C$ (see \cite{Mac}, I.3; \cite{Zhe}, \S73).
The numerator is zero on the hyperplanes $y_j=y_k$
and hence $s_\mu^n$ is a symmetric polynomial in  $y_1,\dots, y_n$.
It is easy to prove that
\begin{equation}
s_{\mu_1,\dots,\mu_n,0,\dots,0}^{n+m}(y_1,\dots,y_n,0,\dots,0)=
s^n_{\mu_1,\dots,\mu_n}(y_1,\dots,y_n)
\label{3.2}
\end{equation}

Let $A$ be a $n\times n$ matrix. Let $y_1,\dots, y_n$ be
the eigenvalues of $A$.
We define the {\it Schur function} $S_\mu(A)$ by
$$S_\mu(A)=s_\mu(y_1,\dots,y_n)$$

By $\rho_\mu=\rho^n_{\mu_1,\dots,\mu_n}$
 we denote the representation
of $\GL(n,\C)$ with the signature $\mu$
(see \cite{Zhe}, \S49; \cite{Mac};
these objects are used only in this Subsection).
 The Schur
function $S_\mu(A)$ is the character of $\rho_\mu$, i.e.,
$S_\mu(A)=\tr \rho_\mu(A)$.
Recall (see \cite{Zhe}, \S73; see also
\cite{Hua}, Theorem 1.2.4 on eliminating of undeterminacy)
 that the dimension of the representation
 $\rho_\mu$ is given by
the Weyl formula
\begin{equation}
\dim \rho^n_\mu= s_\mu(1,\dots,1)=
   \prod_{n\ge i> j\ge 1}\frac{ \mu_i-\mu_j+i-j}{i-j}
   \label{3.3}
   \end{equation}

 Denote by $\Mat_{p,q}$ the space of all
 complex
 $p\times q$ matrices. Denote by $\Pol_{p,q}$ the
 space of all holomorphic polynomials on $\Mat_{p,q}$.
 The linear-fractional transformations
 (\ref{2.7}) for
 $g=\bigl(\begin{smallmatrix} a&0\\0&d \end{smallmatrix}
  \bigr)\in \K$ reduce to the form
 $$\tau_\alpha(g)f(z)=f(a^{-1}zd)\det(a)^{-\alpha}$$
 Let us omit the unessential scalar factor
 $\det (a)^{-\alpha}$ and consider the action
 of $\U(p)\times\U(q)$ given by
     $$ \begin{pmatrix} a&0\\0&d \end{pmatrix}
    :\,\,   f(z)\mapsto f(a^{-1}zd)=f(a^*zd)$$
It is natural to extend this action to a
$\GL(p,\C)\times\GL(q,\C)$-action given by the  formula
$f(z)\mapsto f(a^{t}zd)$.%
\footnote{The groups  $\U(p)\times\U(q)$   and
           $\GL(p,\C)\times\GL(q,\C)$
have the same invariant subspaces in $\Pol_{p,q}$.
Indeed a subspace of homogeneous polynomials
 of a given degree is invariant with respect to
 both groups. It remains to apply the Weyl's unitary  trick,
 see \cite{Zhe}, \S42.}

{\sc
 Theorem \fact.} (Hua, \cite{Hua}, see also \cite{Zhe}, \S56)
\begin{equation}
\Pol_{p,q}
\simeq \bigoplus\limits_{\mu_1\ge\mu_2\ge\dots\ge\mu_p\ge0}
  \rho^p_{\mu_1,\dots,\mu_p}
  \otimes\rho^q_{\mu_1,\dots,\mu_p,0,\dots,0}
\label{3.4}
\end{equation}

This theorem is a consequence of the following Lemma.

\smallskip

{\sc Lemma \fact.}
  {\it Denote by $\Delta_j$ the determinant of
the left  upper $j\times j$ block of  $z$.
All $\GL(p,\C)\times\GL(q,\C)$-highest vectors  in $\Pol_{p,q}$
have the form $\prod_{j=1}^p \Delta_j^{\mu_{j}-\mu_{j+1}}$.}

\smallskip

{\sc Proof of the Lemma.}
Denote by $N_q$ the group of all $d\in \GL(q,\C)$ with units
on the diagonal and zeros under the diagonal.
Denote by $N'_p$  the group of
all $a\in\GL(p,\C)$  with units on the diagonal and
zeros over the diagonal.
 Denote by $A_p$, $A_q$  the diagonal subgroups
 in $\GL(p,\C)$, $\GL(q,\C)$.
 We want to find all eigenfunctions of the subgroup
 $A_pN_p'\times A_q N_q\subset \GL(p,\C)\times\GL(q,\C)$.
 Denote by $\cal R$ the set of all
 $w\in\Mat_{p,q}$ such that $w_{ij}=0$  for $i\ne j$

 A highest vector is
a $N'_p\times N_q$-invariant function
on $\Mat_{p,q}$.
For a  matrix $z\in \Mat_{p,q}$ in a general position
there exists
$H\in {\cal R}$ that can be represented
in the form $H=T z S$, where $T\in N_p'$,
 $S\in N_q$
(it is sufficient to apply the Gauss
 elimination algorithm).
  Moreover, the diagonal matrix elements
 $h_{jj}$ of $H$  are given by
 $h_{jj}=\Delta_j/\Delta_{j-1}$.
 Hence a highest vector is a rational function depending on
  $\Delta_j/\Delta_{j-1}$, $j=1,2,\dots,p$.

  Now we call to mind that a highest vector is an
  $A_p\times A_q$-eigenvector.
  \kvadrat

  \smallskip

  Denote by $\Pol^{\mu}_{p,q}$ the subspaces in $\Pol_{p,q}$
  corresponding  to the summands in (\ref{3.4}).

  \smallskip

{\sc Lemma \fact.}
a) {\it  The function $S_\mu^p(uz^*)$ is a positive definite kernel
on $\Pol_{p,q}$.}

b) $S^p_{\mu_1,\dots,\mu_p}(uz^*)=
S^q_{\mu_1,\dots,\mu_p,0,\dots,0}(z^*u)$

c) {\it The subspace $\Pol_{p,q}^{\mu}$
coincides with the space $\H^\circ[S^p_\mu; \Mat_{p,q}]$
defined by the positive definite kernel $S^p_\mu(uz^*)$.}

\smallskip

{\sc Proof.} a)  Obviously, the representation
$\rho^p_\mu$ extends canonically
 to the representation   of
the semigroup $\Gamma$ of all operators $\C^p\to\C^p$,
and $\rho^p_\mu(A^*)= \rho^p_\mu(A)^*$
(see, for instance \cite{Mac}, Chapter I, Appendix A,
 \cite{Ner2}, 3.3).

Let us show that
$\tr \rho_\mu(AB^*)$   is a positive definite kernel on $\Gamma$.
Let $v_j$ be an orthonormal basis in $\rho^p_\mu$.
Then
$$ \tr \rho_\mu(AB^*)=
\sum_j \langle \rho_\mu(AB^*)v_j,v_j\rangle=
 \langle \rho_\mu(B^*)v_j,\rho_\mu(A^*)v_j\rangle$$
The summands are positive definite
 (since $v_j\mapsto \rho_\mu(A^*)v_j$
 is a function from $\Gamma$
 to Euclidean space, see 1.1) By Proposition 1.6,
  the sum also is positive definite.

 Now we embed $\Mat_{p,q}$ to $\Gamma$ adding $(q-p)$ zero
 rows at the bottom of the matrix.

 b) The nonzero eigenvalues of $uz^*$  and $u^*z$
 coincide (let $h$ be a vector;
 then $zu^*h=\lambda h$ implies $(u^*z)u^*h=\lambda u^*h$).
 Thus the statement follows from (\ref{3.2}).

 c) First, the kernel $S_\mu(uz^*)$ is $\K$-invariant
 and hence it defines a $\K$-invariant subspace.

The elements of $\Pol_{p,q}$
 are holomorphic functions and hence they
are uniquely determined by their restrictions to
the submanifold $M_0$ defined in \Subsection 2.6.
 The elements
of $M_0$ are isometric embeddings $\C^p\to\C^q$. Consider the natural
map $\U(q)\to M_0$: we take a unitary matrix and delete its
last $q-p$ rows. Hence a functions on $M_0$ can be regarded
as a function on $\U(q)$.

 The pullback of the kernel  $S_\mu(z^*u)$ from $M_0$
 to $\U(q)$ is
$$K_\mu(g_1,g_2):=S_\mu\Bigl(g_1 \,\,\,
 \underbrace{\phantom{\begin{matrix}.\\ .\end{matrix}}}_p
 \!\!\!\!\!   \!\!\!\!\!   \!\!\!\!\!
 \begin{pmatrix}\,\,\,\, 1&0\\ \,\,\,\, 0& 0\end{pmatrix}
 g_2^{-1}\Bigr)
 $$
 The kernel $L_\mu(g_1,g_2):=S_\mu(g_1g_2^{-1})$
 on $\U(q)$ defines
 the space $\rho^q_{\mu_1,\dots,\mu_p,0,\dots,0} \otimes
            \rho^q_{\mu_1,\dots,\mu_p,0,\dots,0}$.
By Corollary  1.4, the space $H^\circ[L_\mu]$ is the linear span
of the functions $S_\mu(gA)$, where $A$ ranges over invertible matrices.
This space is finite-dimensional, and thus it is closed with
respect to pointwise convergence. Hence it contains functions
$S_\mu(gA)$   for all noninvertible matrices $A$. Thus
$\H^\circ[K_\mu;\U(q)]\subset H^\circ[L_\mu,\U(q)]$.
 Hence the representation of $\U(q)$ in
 $\H^\circ[S_\mu;\Mat_{p,q}]$ is  the direct sum
 of several copies of $\rho_{\mu_1,\dots,\mu_p,0,\dots,0}$.
 It remains to apply Theorem 3.2.\kvadrat

 \smallskip

 {\bf\punct.  Expansion of the Berezin kernel $K_\alpha$.}

 \nopagebreak

{\sc Theorem \fact} (Hua \cite{Hua})
$$
\prod_{j=1}^p (1-y_j)^{-\alpha}=\sum_{\mu_1\ge\dots\ge\mu_p\ge 0}
c(\mu_1,\dots,\mu_p; \alpha ) \,
s_{\mu_1,\dots,\mu_p}(y_1,\dots,y_p)
 $$
{\it where }
\begin{equation}
   c(\mu_1,\dots,\mu_p; \alpha )=
   \dim \rho^p_{\mu_1,\dots,\mu_p}
 \prod\limits_{j=1}^p   \frac
 {\Gamma(\alpha+\mu_j-j+1)(p-j)!}{\Gamma(\alpha-j+1)(\mu_j+p-j)!}
\label{3.5}
\end{equation}

Proof is contained in \cite{Hua}, Theorem 1.2.5, \cite{FK}.

\smallskip

{\sc Corollary \fact.} (Berezin \cite{Ber1})
\begin{equation}
\det (1-uz^*)^{-\alpha}=
             \sum_{\mu_1\ge\dots\ge\mu_p\ge 0}
c(\mu_1,\dots,\mu_p; \alpha ) \,
s_{\mu_1,\dots,\mu_p}(uz^*)
\label{3.6}
\end{equation}
{\it where  the
constants $c(\mu;\alpha)$
are the same as above {\rm(\ref{3.5})}}.

\smallskip

{\sc Corollary \fact.}
{\it The kernels $c(\mu;\alpha)S_\mu(uz^*)$
and $\det(1-uz^*)^{-\alpha}$
define the same scalar product
in the subspace $\Pol^\mu_{p,q}$.}

\smallskip

{\sc Proof.}
The $H_\alpha$-scalar product in
$\Pol^\mu_{p,q}$ is $\U(p)\times\U(q)$-invariant.
The scalar product in $\Pol^\mu_{p,q}$ defined
by  $S_\mu(uz^*)$ also is $\U(p)\times\U(q)$-invariant.
The $\U(p)\times\U(q)$-module  $\Pol^\mu_{p,q}$
is irreducible. Hence these scalar products differs by a
scalar factor.
Denote by $\sigma_\mu S_\mu(uz^*)$
the kernel defining the $H_\alpha$-scalar product in
$\Pol^\mu_{p,q}$.
The subspaces $\Pol^\mu_{p,q}$ are pairwise
orthogonal,  and hence (see the proof of Proposition 1.6)
  we have $\det(1-uz^*)^{-\alpha}=
 \sum\sigma_\mu S_\mu(uz^*)$.\kvadrat

\smallskip

{\bf \punct. Berezin scale: general values of parameter $\alpha$.}

{\sc Theorem \fact.} (Berezin \cite{Ber1}, and also
\cite{Gin2}, \cite{RV}, \cite{Wal})
a) {\it  The function
$K_\alpha(z,u)=\det(1-uz^*)^{-\alpha}$ is a positive
definite kernel on $\B_{p,q}$
 if and only if $\alpha$ belongs to the set}
\begin{equation}
\alpha=0,1,2,\dots,p-1, \qquad \mbox{or} \qquad \alpha>{p-1}
\label{3.7}
\end{equation}

b) {\it  Denote by $H_\alpha$
the Hilbert space of holomorphic functions on $\B_{p,q}$ defined
by the kernel $K_\alpha$. Then the operators
 $\tau_\alpha(g)$, $g\in \U(p,q)$,
given by {\rm(\ref{3.1})} are unitary in $H_\alpha$.}

c) {\it  If $\alpha>p-1$, then the space $H_\alpha$
contains all holomorphic
 polynomials. If $\alpha=k=0,1,\dots,p-1$,
then $H_\alpha$ has the following $\U(p)\times\U(q)$-module
structure}
$$H_k
 \simeq \bigoplus\limits_{\mu_1\ge\mu_2\ge\dots\ge\mu_k\ge0}
  \rho^p_{\mu_1,\dots,\mu_k,0,\dots,0}
  \otimes\rho^q_{\mu_1,\dots,\mu_k,0,\dots,0}
$$

{\sc Proof.} A sum of positive definite kernels
is positive definite (see Proposition 1.6).
 Hence it is sufficient
to check positivity of the coefficients in the expansion (\ref{3.6}).
The last item is trivial. \kvadrat

{\sc Remark.} The space $H_0$ is one-dimensional, it contains
only constants.

\smallskip

{\sc Remark.} The space $H_q$ coincides with the Hardy space
$H^2(\B_{p,q})$. The scalar product in $H_q$ is given by
$$\langle f,g\rangle_q=
\const \int_{M_0} f(z)\overline {g(z)} \{dz\}
$$
where $M_0$ is the boundary orbit defined in 2.6
and $\{dz\}$ is the unique $\K$-invariant measure
on $M_0$.

\smallskip

{\bf \punct. Berezin scale: small values
of the parameter $\alpha$.} Consider the case
$\alpha=k=0,1,\dots,p-1$.
Let
$$\partial_{kl}=\frac{\partial}{ \partial z_{kl}}$$
Consider the matrix
\begin{equation}
 D=\begin{pmatrix}
 {\partial_{11}}&\dots &   {\partial_{1q}}\\
\vdots&\ddots&\vdots\\
 {\partial_{p1}}&\dots &  {\partial_{pq}}
\end{pmatrix}
\label{3.8}
\end{equation}
We regard the minors
of the matrix $D$ as differential operators.

\smallskip

{\sc Theorem \fact.}  {\it Let $f\in H_k$, $k=0,1,\dots,p-1$.
Then each $(k+1)\times(k+1)$ minor of $D$
annihilates $f$, i.e., for all $i_1< i_2<\dots <i_k$, $j_1< j_2<\dots <j_k$
the following identity holds}
\begin{equation}
\det\begin{pmatrix}
{\partial_{i_1 j_1}}&\dots &
          {\partial_{i_1 j_{k+1}}}\\
\vdots&\ddots&\vdots\\
{\partial_{i_{k+1}j_1}}&\dots &
       {\partial_{i_{k+1}j_{k+1}}}
\end{pmatrix} f(z)=0
\label{3.9}
\end{equation}

{\sc Proof.}  By  Corollary 1.4 and holomorphness of $f$,
 it is sufficient to prove
the statement for elements of the supercomplete system, i.e.,
for the functions
$$f(z)=\det(1-za^*)^{-k}$$
Obviously, the system of partial differential
equations (\ref{3.9})
is invariant with respect to the transformations
$$z\mapsto AzD, \qquad
\mbox{where} \quad A\in\GL(p,\C), D\in \GL(q,\C)$$
Hence it is sufficient to consider the case
in which block $p\times(p+(q-p))$ matrix  $a$ has the form
 $$a=\begin{pmatrix}1_p&0\end{pmatrix}
$$
where $1_p$ is the unit $p\times p$ matrix.
Now $f(z)=\det^{-k}(1-za^*)$
 depends only on left $p\times p$ block of $z$.
Hence, without loss of generality, we can consider
only the case $p=q$ and $f(z)=\det(1-z)^{-k}$.
But the system (\ref{3.9})
 is invariant with respect to translations
and we can change our function $f(z)$ to $\det (z)^{-k}$.

A direct (but pleasant) calculation shows that
$$\bigl({\partial_{11}   \partial_{22}} -
 {\partial_{21}
     \partial_{21}}\bigr)\det(z)^{-1}=0$$
This  implies that $3\times 3$ minors of  D
     annihilate
$\det z^{-2}=\det z^{-1}\det z^{-1}$
etc. \kvadrat

\smallskip

{\sc Remark.} Let $f$ be a solution of the system (\ref{3.9})
in $\B_{p,q}$. Generally, $f\notin H_k$, since $f$
can have too rapid growth near boundary.
Nevertheless $f$ can be approximated by
finite sums $\sum_j \det(1-za^*_j)^{-k}$ in the topology
of uniform convergence on compact sets.
\hfill $\square$

\smallskip

{\bf\punct. Gindikin--Vergne--Rossi description for $H_\alpha$.}
We consider only  the case $p=q$.
First we observe that the Berezin kernel $K_\alpha$ in
the wedge realization
$\W_q$ has
the form
\begin{equation}
K_\alpha(T,S)=\det( T^*+S)^{-\alpha}
\label{3.10}
\end{equation}

We denote by $H_\alpha(W_q)$ the space of holomorphic functions
on $W_q$ defined by the kernel (3.11). The group
$\U(q,q)$ acts in $H_\alpha(W_q)$ by the transformations
$$
f(T)\mapsto
f\bigl(A+TC)^{-1}(B+TD)\bigr)\det(A+TC)^{-\alpha}
$$
where the matrix $\bigl(\begin{smallmatrix} A&B\\C&D\end{smallmatrix}\bigr)$
satisfies  condition (2.4).

{\sc Remark.} Let  $Z\in \W_q$.
 Let $\lambda_j$ be the eigenvalues of $Z$,
and $v_j$ be the eigenvectors.
Then
$$0<\langle (Z+Z^*)v_j,v_j\rangle=  \langle Zv_j,v_j\rangle+
\langle v_j,Zv_j\rangle=(\lambda_j+\ov\lambda_j)
\langle v_j,v_j\rangle$$
Hence $\Re\lambda_j>0$ and hence the function
$\det Z^{-\alpha}:=\prod_j \lambda_j^{-\alpha}$
is well defined on $\W_q$.  Thus expression (\ref{3.10})
makes sense.

\smallskip

Consider the cone $\Pos_q$ consisting of all
positive semidefinite complex $q\times q$
Hermitian matrices $X$.
 The group $\GL(q,\C)$ acts  on $\Pos_q$
 by the transformations $h:X\mapsto h^* X h$,
 $h\in\GL(q,\C)$. Denote by $N_k$, where $k=0,1,\dots,q$,
 the set of all matrices $X\in\Pos_{q}$
 such that $\rk X=k$.
 Obviously, the sets $N_k$ are
 the  $\GL(q,\C)$-orbits on $\Pos_q$.

 Let $\chi(X)$ be a tempered
distribution supported by $\Pos_p$.
Its {\it Laplace transform} is
\begin{equation}
\widehat\chi(T)=\int_{\Pos_q} \exp(-\tr XT) \chi(X) dX
\label{3.11}
\end{equation}
The function  $\widehat\chi(T)$ is a holomorphic function of
 polynomial
growth on the Siegel wedge  $\W_q$ (see \cite{Vla}).

\smallskip

{\sc Theorem \fact.} (\cite{Gin2},\cite{RV})
{\it Let $\alpha$ satisfy
the positive definiteness conditions {\rm(\ref{3.7})}. Then
$\det (T)^{-\alpha}$ is the Laplace transform of
some positive measure $\nu_\alpha$
on $\Pos_q$. If $\alpha>q-1$, then
$\nu_\alpha=\det X^{\alpha-q}dX$.
For $\alpha=k=0,1,\dots, q-1$, the measure $\nu_k$
is the unique up to a scalar factor $\SL(q,\C)$-invariant measure
on the boundary orbit $N_k$.}

\smallskip

{\sc Proof.}
a) Let  $\alpha>q-1$.
We must check the equality
\begin{equation}
\int_{\Pos_q} \det X^{\alpha-q}\exp(-\tr XT)\,dX=\const
    \cdot \det (T)^{-\alpha}
    \label{3.12}
\end{equation}
Since the both parts are holomorphic in $T\in W_q$, it is sufficient
to consider the case $T=T^*$. Then the integral converts to
  \begin{multline}
  \int_{\Pos_q}
  \det X^{\alpha-q}\exp(-\tr T^{1/2}XT^{1/2})\,dX=\\
  =\det T^{-\alpha}
  \int_{\Pos_q} \det (T^{1/2}XT^{1/2})^{\alpha-q}
  \exp(-\tr T^{1/2}XT^{1/2})\,d (T^{1/2}XT^{1/2})=\\=
  \det T^{-\alpha} \cdot \int_{\Pos_q}
  \det Y^{\alpha-q}\exp(-\tr Y)\,dY
  \label{3.13}
  \end{multline}
as required.

  In fact, we
  only use the $\GL(q,\C)$-homogeneity of the measure
  $\nu_\alpha$:
  \begin{equation}
  \nu_\alpha(h^*Ah)=|\det(h)|^{2\alpha-2q}\nu_\alpha(A)
  \label{3.14}
  \end{equation}
  where $h\in\GL(q,\C)$, and $A$ is a subset in $\Pos_q$.

b) Let $\alpha\le q-1$.
 First we must check the existence of a $\SL(q,\C)$-invariant
volume form on $N_k$.   We have
$N_k=\SL(q,\C)/{\cal H}$, where ${\cal H}$
is
the stabilizer $\cal H$ of the point
$u_k:=
\bigl(\begin{smallmatrix} 1_k &0\\0&0\end{smallmatrix}
\bigr)\in N_k$. The subgroup $\cal H$
consists of matrices
$$g=\begin{pmatrix}R&0\\0&S\end{pmatrix},
 \qquad\text{where}\quad R\in\U(p),\,S\in\GL(q-p,\C)\quad\text{and}\quad
              \det(R)\det(S)=1
$$
The group ${\cal H}$  has no homomorphisms
to  $\R$. Hence the (unique) volume form on the tangent space
to $N_k$ at the point $u_k$ is $H$-invariant.
Then we define a volume form on $N_k$ by $\SL(n,\C)$%
-invariance.

Secondly we must  evaluate the Laplace transform
 of the measure $\nu_k$. It is not hard to
 check, that the measure $\nu_k$ satisfies the same homogeneity
 condition (\ref{3.14}), $\alpha=k$, and this gives
the formula for the Laplace transform of $\nu_k$.
 \kvadrat

  \smallskip

For the case $p=q$, this gives an independent proof
  of Theorem 3.8 about positive definiteness conditions,
i.e., we obtain the following corollary.

\smallskip

{\sc Corollary \fact.} {\it The kernel
$K_\alpha(T,S)=\det(T^*+S)^{-\alpha}$
on $\W_q$
is positive definite if $\alpha>q-1$ or
 $\alpha=0,1,\dots,q-1$.}

\smallskip

{\sc Proof.} We use directly the definition
 of positive definiteness. We must show  positivity
 of the expression
 \begin{multline*}
 \sum\limits_{1\le k,l\le N}
 \det( T_k^*+T_l)^{-\alpha}\xi_k\ov\xi_l=
 \sum\xi_k\ov\xi_l
 \int_{\Pos_q} \exp(-\tr(\ov T_k+T_l)X)\,d\nu_\alpha(X)=  \\=
 \int_{\Pos_q} |\sum\xi_k\exp(-\tr\ov T_kX)|^2d\nu_\alpha(X)
 \ge 0
 \end{multline*}
 as required. Of course, we repeated  the
  proof of the trivial
 part of the Bochner theorem (see \cite{RS}, v.2, Theorem IX.9).
 \kvadrat

{\sc Theorem \fact.} {\it The Laplace transform%
\footnote{We slightly change the notation with
respect to (\ref{3.11}) }
\begin{equation}
\widetilde f(T)=
\int_{\Pos_q} f(X)\exp(-\tr XT)d\nu_\alpha(X)
\label{3.15}
\end{equation}
 is
a unitary {\rm(}up to a scalar factor{\rm)}
 operator
from $L^2(\Pos_q, \nu_\alpha)$
to the Berezin space $H_\alpha[\W_q]$}.

\smallskip

{\sc Proof.} Consider the family of functions
$\kappa_A(X)=\exp(-\tr A^*X)$, where $A\in \W_q$,
 on $\Pos_q$.
It is sufficient to prove that
\begin{equation}
\langle \kappa_A,\kappa_B\rangle_{L^2(\Pos_q,\nu_\alpha)}=
\const\cdot\langle \widetilde\kappa_A,\widetilde\kappa_B
\rangle_{H_\alpha(\W_q)}
\label{3.16}
\end{equation}
By (\ref{3.12}), the left-hand side is
$\const\cdot\det(A^*+ B)^{-\alpha}$.
By the same formula (\ref{3.12}),
 $\widetilde\kappa_A(T)
=\det (A^*+T)^{-\alpha}$. These $\widetilde\kappa_A$
are the elements
of the supercomplete system of $H_\alpha(\W_q)$
and we know their scalar products.
\kvadrat

\smallskip

{\sc Remark.} Assign $\kappa_A$ be a supercomplete system
in $L^2(\Pos_q,\nu_\alpha)$. Then the Laplace
transform (\ref{3.15}) coincides with the transform
$\H\to \H^\circ$ described in \Subsection 1.2.
\hfill $\square$

\smallskip

{\sc Remark.}  Theorem 3.10
implies Theorem 3.9 in the case
$p=q$. For definiteness,
 assume $\alpha=q-1$. Let $D$ be given
by (\ref{3.8}). A $\delta$-distribution $f(X)\nu_{q-1}$
satisfies the condition $\det X\cdot f(X)\nu_{q-1}=0$.
After the Laplace transform, we obtain
 $D \widehat  {f(X)\nu_{q-1}}=0$ (this is not a
 complete proof). \hfill $\square$

\smallskip

{\small

{\bf \punct. Comments.} 1) The construction described
in this Section for the symmetric spaces
  $\U(p,q)/\U(p)\times\U(q)$ works for
all Hermitian noncompact symmetric spaces, i.e., also for
$\Sp(2n,\R)/\U(n)$, $\SOS(2n)/\U(n)$,
   $\SO(p,2)/\SO(p)\times \SO(2)$
   and two exceptional spaces, see \cite{Ber1},
   \cite{Gin1}, \cite{RV}, \cite{Wal}.

   \smallskip

2) The constant in formula   (\ref{3.12}) is a special case
of Gindikin's matrix $\Gamma$-function \cite{Gin1},
see also the exposition in \cite{Ter}, \cite{FK}.

 \smallskip

3)  After Theorem 3.12 there arises
a problem of transfering
of the action of the group $\U(q,q)$
to the space   $L^2(\Pos_q,\nu_\alpha)$.
The Lie algebra ${\frak u}(q,q)$ acts in $L^2(\Pos_q,\nu_\alpha)$ by
second order partial differential operators,
which can easily
 be written explicitly.
The exponents of these differential operators
(i.e elements
of the group $\U(q,q)$ itself)
are integral operators with
 kernels involving matrix Bessel functions,
see \cite{Her}, \cite{Ter}, \cite{FK}.

 Consider the case $\G=\U(1,1)$.
Then $\W_q=\W_1$ is the Lobachevskii plane
$\Re T>0$. The transformation  $T\mapsto T^{-1}$
of $W_1$ corresponds to the Hankel transform in
$L^2(\R_+,\nu_\alpha)$ (the Tricomi theorem).

 Rotations of the Lobachevskii plane
 with center at $T=1$ corresponds to the
Kepinski and Myller-Lebedeff \cite{ML}
 explicit solution
of the Cauchy problem for the Schr\"odinger equation
$$
\frac{i\partial } {\partial t} F(x,t)=
 \Bigl(- \frac{\partial^2 } {\partial x^2}+
 x^2+\frac {a}{x^2}\Bigr) F(x,t)
 $$    }

\bigskip

{\large\bf\sect  Kernel representations and spaces $\V_\alpha$}

\nopagebreak

\medskip

{\bf\punct. Definition of  kernel representations.}
Let
$$\alpha=0,1,2,\dots, p-1,\qquad \mbox{or} \quad \alpha>p-1$$
Consider   the kernel $L_\alpha(z,u)$ on $\B_{p,q}$ given by
\begin{align}
L_\alpha(z,u):=&|\det(1-zu^*)|^{-2\alpha}=|\det(1-uz^*)|^{-2\alpha}=
= \label{4.1}
\\=
  &\det(1-zu^*)^{-\alpha}
     \overline{ \det(1-zu^*)}^{\,\,-\alpha}  \nonumber
  \end{align}
 The kernel $L_\alpha$ is the product
  of two positive definite kernels
  and, by Proposition 1.6,
   $L_\alpha$ is positive definite.
  We denote by $V_\alpha$ the Hilbert space $\H^\circ[L_\alpha]$
  defined by the kernel
  $L_\alpha$.

  The group $\G=\U(p,q)$ acts in the space $V_\alpha$ by  the
  unitary operators
  \begin{equation}
  T_\alpha\begin{pmatrix}a&b\\c&d\end{pmatrix}
   f(z)=f\bigl((a+zc)^{-1}(b+zd)\bigr)|\det(a+zc)|^{-2\alpha}
\label{4.2}
   \end{equation}
  We say that the representation $T_\alpha$ is a
  {\it kernel-representation}
  of the group $\G=\U(p,q)$.

  \smallskip

  {\bf\punct.  Decomposition of
  the kernel representations into a tensor product.}
  Consider  the unitary representation $\tau_\alpha$  of $\G$
  defined by (\ref{3.1})
   and the complex conjugate
    representation $\ov\tau_\alpha$.
    Consider the tensor product
     $\tau_\alpha\otimes \ov\tau_\alpha$.
     The space $V_\alpha:= \ov H_\alpha\otimes  H_\alpha$
   of the tensor product
    is the space of holomorphic functions
   on $\B_{p,q}\times\B_{p,q}$
   defined by the positive definite kernel
   \begin{equation}
   \widetilde L(z_1,z_2;\, u_1,u_2)=\det(1-u_1z_1^*)^{-\alpha}
                            \det(1-u_2z_2^*)^{-\alpha}
\label{4.3}
\end{equation}
and the group $\G$ acts in $V_\alpha$ by
\begin{multline}
F(z_1,z_2)\mapsto\\
 F\bigl[(\ov a+z_1\ov c)^{-1}(\ov b+z_1\ov d),\,\,
            (a+z_2 c)^{-1}(b+z_2d)\bigr]
           \det (\ov a+z_1\ov c)^{-\alpha}
           \det (a+z_2c)^{-\alpha}
\label{4.4}
\end{multline}
Denote by $\Delta$ the diagonal $z_1=\ov z_2$ in
 $\B_{p,q}\times\B_{p,q}$.

 First a holomorphic function on
 $\B_{p,q}\times\B_{p,q}$ is determined by
 its restriction to $\Delta$.  Hence we
 can consider the space  $H_\alpha\otimes \ov H_\alpha$
 as a space of functions on $\Delta$.
 The reproducing kernel of the latter space
 can be obtained  by the substitution $z_1=\ov z_2$,
 $u_1=\ov u_2$
 to (\ref{4.3}), see \Subsection 1.3.

 Secondly the transformations in square brackets
 in formula (\ref{4.4}) preserve $\Delta$. Let us
  rewrite the operators (\ref{4.4})
 as operators in the space of functions on $\Delta$.
 This gives (\ref{4.2}). Thus we obtain the following
 proposition.

 \smallskip

 {\sc Proposition \fact.} {\it The operator of
  restriction of
  functions $F(z_1,z_2)$
  to $\Delta$ is a unitary
  $\G$-intertwining operator
 from   $\ov\tau_\alpha\otimes\ov\tau_\alpha$ to the
 kernel representation $T_{\alpha}$.}

 \smallskip

 {\bf \punct.  Representation in the space
  of Hilbert--Schmidt operators.}  Let $W_1$,
  $W_2$ be Hilbert spaces. Let $A$ be an operator
  $A:W_1 \to W_2$. Let $a_{ij}$ be the matrix elements
   of $A$ in some
  orthonormal bases. Recall that $A$ is called
  {\it a Hilbert-Schmidt
  operator} (see \cite{RS}, v.1, VI.6) if  $\sum  |a_{ij}|^2<\infty$.

  The space ${\cal HS}(W_1,W_2)$ of all Hilbert--Schmidt
  operators $W_1 \to W_2$ is a Hilbert space with respect to
  the scalar product
  \begin{equation}
  \langle A,B\rangle_{\cal HS}=\tr AB^*
  \label{4.5}
  \end{equation}
  There is an obvious canonical identification
  $${\cal HS}(W_1,W_2)\simeq\ov  W_1\otimes W_2$$
  Hence we can identify
  $$V_\alpha\simeq\ov  H_\alpha\otimes H_\alpha
  \simeq {\cal HS}(H_\alpha,H_\alpha)$$
  The group $\G$ acts in   ${\cal HS}(H_\alpha, H_\alpha)$
  by the conjugations
  \begin{equation}
  g: A\mapsto \tau_\alpha(g)^{-1} A \tau_\alpha(g)
 \label{4.6}
  \end{equation}

  \smallskip

{\bf \punct. Another version of the kernel $L$.}
Let us define the {\it modified Berezin kernel}
\begin{equation}
\L_\alpha(z,u)
 := \frac {\det(1-zz^*)^{\alpha}\det(1-uu^*)^{\alpha} }
       {|\det(1-zu^*)|^{2\alpha}}
\label{4.7}
       \end{equation}
We denote by  $\V_\alpha$    the space
  $\H^\circ[\L_\alpha]$.
The kernel $\L_\alpha$ is $\G$-invariant in the following sense
    \begin{equation}
    \L_\alpha\left((a+zc)^{-1}(b+zd),      (a+uc)^{-1}(b+ud)\right)
      = \L_\alpha(z,u)
      \label{4.8}
      \end{equation}
(this easily follows from (\ref{2.10})).
Hence the operators
  \begin{equation}
   {\cal T}_\alpha
   \begin{pmatrix}a&b\\c&d\end{pmatrix}
   f(z)=f\bigl((a+zc)^{-1}(b+zd)\bigr)
   \label{4.9}
   \end{equation}
 are unitary in $\V_\alpha$.

 \smallskip

  {\sc Remark.} The last formula
 does not  depend on $\alpha$. Nevertheless we shall see
(Section 7)
 that generally the unitary representations $\T_\alpha$ are nonequivalent.
 \hfill $\square$

 \smallskip

 We also define the {\it distinguished vector}
 $\Xi=\Xi_\alpha\in\V_\alpha$ by
 \begin{equation}
 \Xi_\alpha(z):=\det(1-zz^*)^{-\alpha}=\L_\alpha(z,0)
 \label{4,10}
 \end{equation}
The vector $\Xi_\alpha$ is an element of the supercomplete
basis and the $\G$-orbit of $\Xi_\alpha$ consists
of all elements of the supercomplete system.
Hence the vector $\Xi_\alpha$ is cyclic%
\footnote{Let $\rho$ be a representation of a group $G$
in a topological vector space $W$ (see \cite{Kir}, 7.2).
 A vector $w\in  W$
is called {\it cyclic} if the linear span of the vectors
$\rho(g)w$, $g\in G$, is dense in $W$. For a subset $B\subset W$
we define its {\it cyclic span} as the closure of the linear
span of all vectors $\rho(g)w$, where $g\in G$ and
 $w\in B$, see \cite{Kir}, 4.4.}.

\medskip

  By Proposition 1.7,
   the kernels $L$ and $\L$ differ unessentially,
  and the canonical  operator ${\cal U}:V_\alpha\to\V_\alpha$,
   described
  in Proposition 1.7, is given by
   $${\cal U}f(z)=f(z)\det(1-zz^*)^{\alpha}$$
 This operator intertwines  (\ref{4.2}) and (\ref{4.9}).
 Thus  we identify the unitary representation
 $\T_\alpha$ and $T_\alpha$ of the group $\G$.

 \smallskip

  {\it The last description of the kernel representation
   {\rm(}kernel $\L_\alpha$ and the representation $\T_\alpha${\rm)}
   will be main for us below.}

   \smallskip

 {\bf \punct. Limit of kernel representations
 as $\alpha\to\+\infty$.}
 Consider the space $\H^\star[\L_\alpha]$ defined
 in \Subsection 1.6.
 Consider elements of  $\H^\star[\L_\alpha]$ of the form
$$
 \phi(z)\det(1-zz^*)^{-p-q}\{dz\}
$$
 where $\phi(z)$ are compactly supported smooth
 functions on $\B_{p,q}$.
 The $\H^\star$-scalar product  of
such  distributions
is given by
\begin{multline}
\langle \phi,\psi\rangle_{(\alpha)}:=\\=
C(\alpha)^{-1}
\iint\limits_{\B_{p,q}\times \B_{p,q}} \L_\alpha(z,u) \phi(z)
\ov{\psi(u)}\det(1-zz^*)^{-p-q}\det(1-uu^*)^{-p-q}\{dz\}\{du\}
\label{4.11}
\end{multline}
We define the
normalization constant $C(\alpha)$ as
the right-hand side of
(\ref{2.18}). By (\ref{4.8}),
 this scalar product is $\G$-invariant.

\smallskip

{\sc Lemma \fact.} {\it The family of distributions
$$C(\alpha)^{-1} \L_\alpha(z,u)
\det(1-zz^*)^{-p-q}\det(1-uu^*)^{-p-q}\{dz\}\{du\}  $$
 converges
 to $\det(1-uu^*)^{-p-q}\delta(z-u)$
as $\alpha\to+\infty$.}

\smallskip

{\sc Proof.} By the invariance property (\ref{4.8}),
 it is sufficient to follow only limit behavior
 of the family of
 distributions
 $$\varOmega_\alpha(z):=C(\alpha)^{-1}L(z,0)\det(1-zz^*)^{-p-q}\{dz\}=
 C(\alpha)^{-1}\det(1-zz^*)^{\alpha-p-q}\{dz\}$$
 By the choice of the normalization constant,
 the integral of the function $\varOmega_\alpha$
  is 1. We also have
 $\det(1-zz^*)=1$ if $z=0$ and $\det(1-zz^*)<1$ otherwise.
 The last two remarks impliy the convergence
 $\varOmega_\alpha(z)\to \delta(z)$ as $\alpha\to+\infty$.
  \kvadrat

  \smallskip

  Lemma 4.2 implies the following statement.

  \smallskip

  {\sc Theorem \fact.}
  {\it The family of scalar products {\rm(\ref{4.11})}
  tends to             the following
  $L^2$-scalar product with respect to the
  $\G$-invariant measure
   $\det(1-zz^*)^{-p-q}\{dz\}$:
  $$\langle\phi , \psi\rangle=
        \int_{\B_{p,q}} \phi(z)\ov{\psi(z)}
  \det(1-zz^*)^{-p-q}\{dz\}$$
  as $\alpha\to+\infty$.}

\smallskip

We observed that
{\it the natural limit of kernel representations
${\cal T}_\alpha$ as $\alpha\to+\infty$ is the space $L^2$ on
Riemannian symmetric space $\U(p,q)/\U(p)\times\U(q)$.}

\smallskip

{\bf\punct. A canonical basis in  space of
$\K$-invariant functions.}
Denote by $\V_\alpha^\K$ the space of all
functions $f\in \V_\alpha$ that are invariant
with respect to the group $\K=\U(p)\times\U(q)$.

\smallskip

{\sc Proposition \fact.} {\it  Let $\mu$
be a collection of integers
ranging the domain
$\mu_1\ge\dots\ge\mu_p\ge 0$  if
 $\alpha>p-1$,
and
 $\mu_1\ge\dots\ge\mu_k\ge\mu_{k+1}=\dots=\mu_p= 0$
if
 $\alpha=k\le p-1$.
Then the system of functions
\begin{equation}
 \Delta_{\mu}(z)=\Delta_{\mu_1,\dots,\mu_p}(z)
:=S_{\mu}(zz^*)\det(1-zz^*)^{\alpha}
\label{4.12}
\end{equation}
where $S_\mu$ are the Schur functions,
forms an orthogonal basis in the space $\V_\alpha^\K$ and }
\begin{equation}
\|\Delta_\mu\|_{\V_\alpha^\K}^2=
\frac{\prod_{j=1}^p\Gamma^2(\alpha-j+1)}
     {\prod_{j=0}^{q-1} j!}
 \prod\limits_{j=1}^p   \frac  {(\mu_j+p-j)!(\mu_j+q-j)!}
 {\Gamma^2(\alpha+\mu_j-j+1)}
 \label{4.13}
 \end{equation}

 {\sc Proof.}  First we explain the origin of this basis.
 Consider the model ${\cal HS}(H_\alpha,H_\alpha)$
 of the kernel representation (see 4.3). The space
 $\V^\K_\alpha$ corresponds to the space of
 $\K$-intertwining operators $H_\alpha\to H_\alpha$.
 In \Subsection 3.2 we constructed the decomposition
 $$H_\alpha=\oplus_\mu \Pol^{\mu}_{p,q}$$
 into the sum of pairwise distinct representations
 of  $\K$. Let $c(\mu;\alpha)$ be given by (3.6).
 By the Schur Lemma (see \cite{Kir}, 8.2), any
  $\K$-intertwining operator is
 a scalar operator in each summand. Our basic element
 $\Delta_\mu$ corresponds to
 the operator $J_\mu$, that is
 $c(\mu;\alpha)$
 on $ \Pol^{\mu}_{p,q}$ and 0 on  $\Pol^{\nu}_{p,q}$
 for $\nu\ne\mu$.

 It remains to evaluate the Hilbert--Schmidt norm
 of the operator $J_\mu$. By (\ref{4.5}), it equals
 $$
 \dim\Pol^{\mu}_{p,q}=
 \dim \rho^p_{\mu_1,\dots,\mu_p}\cdot
 \dim \rho^q_{\mu_1,\dots,\mu_p,0,\dots,0}
 $$
 and we apply the Weyl formula (\ref{3.3}). \kvadrat

 \smallskip

 {\small

 {\punct\bf. Comments.}
 1) Formula (\ref{4.7}), (\ref{4.9}) make sense for an arbitrary
 matrix ball (see 2.10), and this defines
 {\it scalar} kernel representations
 of any classical group $G$, see \cite{Ner3},
  \cite{Ner7}.

  \smallskip

  2)   General kernel representations
  (the definition was given in \cite{Ner3},
   see also \cite{NO}) are
  realized in spaces of {\it vector-valued} functions
   on matrix balls, and analogues of the Berezin kernels
   are
 invariant matrix-valued positive definite kernels
 on matrix balls.

   We shall give the definition of the kernel
   representations more formally.
    Consider an Hermitian
 symmetric space
 $\widetilde G/\widetilde K$ (see \Subsection 3.7).
 Consider a symmetric subgroup%
 \footnote{A symmetric subgroup in a semisimple
 (resp. classical, reductive) group $\widetilde G$ is
 the set of fixed points of an automorphism of
 order 2.}
  $G\subset \widetilde G$,
 denote by $K$ the maximal compact subgroup in $G$.
 It turns out to be that there are two possibilities.

    a) {\it The first case.} $G/K$ is a complex
       submanifold in  $\widetilde G/\widetilde K$

     b){\it The second case.} $G/K$ is a totally
     real
       submanifold in  $\widetilde G/\widetilde K$
       and $\dim G/K=\frac12 \dim \widetilde G/\widetilde K$.

 Consider  a unitary highest weight representation $\tau$
 of $\widetilde G$ and its restriction to $G$.
 In the first case the spectrum of
 the restriction is discrete
 and the problem of decomposition  is reduced
 to a purely combinatorial problem (\cite{JV}).

  By definition, the restriction in the second case gives
 a kernel representation.

 For a discussion of a priori relations
 of this problems with other spectral
 problems (Howe dual pairs, $L^2$ on Stiefel manifolds,
 $L^2$ on pseudo-Riemannian symmetric spaces)
 see \cite{Ner3}.

 \smallskip

   3)  Only 4 real exceptional
    groups (from 23) have kernel representations.

\smallskip

    4) The case $\widetilde G=\OO(2,n)$
 differs from the cases
$\widetilde G=\Sp(2n,\R),\U(p,q),\SOS(2n)$
({\it "matrix balls cases"}) in many details.

   5)
  Canonical bases
 exist  for all scalar kernel representations
 (in  fact, for any irreducible representation
 of $\widetilde K$, the space of $K$-fixed vectors
 has dimension $\le 1$).
For matrix ball cases they
 consists of Jack polynomials (for a definition
 of the Jack polynomial see \cite{Mac});
  in our case $G=\U(p,q)$, the Jack
 polynomial reduces to the Schur functions.
 Existence of the canonical bases for the case in which
$G/K$ is an Hermitian space
 was observed in \cite{OZ}.

\smallskip

 }

 \bigskip

 {\large\bf \sect Index hypergeometric transform: preliminaries}

 \nopagebreak

        \medskip

 \def\F{{}_2F_1}
 \def\Fbc{\F(b+is,b-is;b+c;-x)}
 \def\dxbc{x^{b+c-1}(1+x)^{b-c}dx}
 \def\dsbc{\Bigl| \frac{\Gamma(b+is)\Gamma(c+is)}{\Gamma(2is)}\Bigr|^2ds}
  \def\dsabc{ \Bigl| \frac{\Gamma(a+is)\Gamma(b+is)\Gamma(c+is)}
                  {\Gamma(2is)}\Bigr|^2ds}

 We use the standard notation for
 the hypergeometric functions
 \begin{align*}
& \F(a,b;c;x)=\sum_{j=0}^\infty \frac{(a)_n(b)_nx^n}
                        {n! \,(c)_n}
                        \\
&{}_3F_2\Bigl[\begin{matrix} a,b,c\\d,e\end{matrix};x\Bigr]=
           \sum_{j=0}^\infty \frac{(a)_n(b)_n(c)_nx^n}
                        {n!\, (d)_n (e)_n}
\end{align*}
where $(r)_n=r(r+1)\dots(r+n-1)$ is the Pochhammer symbol.

{\bf\punct. Index hypergeometric transform}.
Fix $b,c>0$.
Consider a sufficiently decreasing function $f(x)$ on the half-line
$\R_+: \,\,x\ge 0$. We define the integral transform
$J_{b,c}f$ by
\begin{multline}
g(s)=J_{b,c} f(s)=
\frac 1{\Gamma(b+c)} \int_0^{+\infty} f(x)\,\, \F(b+is,b-is;b+c;-x)
  \dxbc
  \label{5.1}
  \end{multline}
  The inverse transform is given by
  \begin{multline}
   J_{b,c}^{-1} g(x)=
\frac 1{\Gamma(b+c)} \int_0^{+\infty}g(s)\,\,\Fbc\dsbc
\label{5.2}
\end{multline}
The integral transform $J_{b,c}$
 is called the {\it index hypergeometric transform},
or  the {\it Olevsky transform},
 or the {\it generalized Fourier transform}
or the {\it Jacobi transform},
or the  {\it Fourier--Jacobi} transform,
 see \cite{Koo1}, see also \cite{Ner8}.
For the first time, the inversion formula
 was obtained by Weyl
in 1910 (\cite{Wey}).

\smallskip

{\bf\punct. Plancherel formula.}  The following
 statement is equivalent to the inversion formula
 (\ref{5.2}).

{\sc Theorem \fact} (Weyl \cite{Wey}, \cite{Koo1}) {\it
The transform $J_{b,c}$ is a unitary operator}
$$J_{b,c}: L^2\Bigl(\R_+,\dxbc\Bigr) \to
  L^2\Bigl(\R_+,\dsbc\Bigr)
  $$

 {\bf\punct. Exotic Plancherel formulas.}
 Denote by   $W^a_{b,c}$ the Hilbert space
 of holomorphic functions
 in the disc $|z|<1$ with the scalar product
 $$\langle f, g\rangle=
   \frac 1{2\pi\Gamma(2a-1)}
    \iint\limits_{|z|<1} f(z)\overline{g(z)} (1-|z|^2)^{2a-2}
    \F(a-b,a-c;2a-1; 1-|z|^2)\{dz\}$$
    where $\{dz\}$ is the Lebesgue measure on the circle.
    A simple calculation shows that the reproducing kernel of the space
    $W^a_{b,c}$ is
 \begin{equation}
 N^a_{b,c}(z,u)=\frac{\Gamma(a+b)\Gamma(a+c)}
   {\Gamma(b+c)} \F\Bigl[ \begin{matrix} a+b,a+c\\ b+c\end{matrix}
                       ;\ov z u \Bigr]
\label{5.3}
\end{equation}

{\sc Theorem \fact.} \cite{Ner8} {\it The operator
\begin{multline}
J^a_{b,c}g(s) = \frac 1{|\Gamma(a+is)|^2 \Gamma(b+c)}
  \int\limits_0^\infty (1+x)^{-a-b} g\Bigl( \frac x{x+1} \Bigr)\times
  \\ \times
    \Fbc \dxbc
\label{5.4}
    \end{multline}
 is a unitary operator                }
 \begin{equation}
 W^a_{b,c}\to L^2\Bigl(\R_+, \dsabc\Bigr)
 \label{5.5}
 \end{equation}

 {\sc Remark.} Let $g\in  W^a_{b,c}$. Consider the new function
 $ f(x):= g\bigl(x/(x+1)\bigr)$.
 The function $f$ is defined on the half-plane
 $\Re x>-1$ and thus
  we obtain a space of holomorphic functions on
 this half-plane. The kernel (\ref{5.3}) is replaced by
$$
N^a_{b,c}(x,y)=
\frac{\Gamma(a+b)\Gamma(a+c)}
{\Gamma(b+c)} \F\Bigl[ \begin{matrix} a+b,a+c\\ b+c\end{matrix}
                       ; \frac{x\ov y}{(1+x)(1+\ov y)} \Bigr]
$$
In formula (\ref{5.4}) the factor $g(x)$ changes to $f(x)$
(and the formula
will   almost coincided with (\ref{5.1}).
 The Plancherel measure
 on the half-line $\R^*$ will be the same as in (\ref{5.5}).
 \hfill$\square$

 \smallskip

{\bf \punct. Index transform and Hahn polynomials.}
Let $a,b,c>0$. The {\it continuous dual Hahn polynomials}
(see, for instance \cite{AAR}, 6.10, \cite{KS})
are defined by the formula
\begin{equation}
\SS_n(s^2;a,b,c)=(a+b)_n(a+c)_n \,\,{}_3 F_2\Bigl[
     \begin{matrix} -n,a+is,a-is\\ a+b, a+c\end{matrix}
     ;1\Bigr]
\label{5.6}
     \end{equation}
The family of  functions $\SS_n$ is an orthogonal basis
in the space
$$      L^2\Bigl(\R_+, \dsabc\Bigr) $$
and moreover
\begin{multline}
\int_0^\infty \SS_n(s^2;a,b,c) \SS_m(s^2;a,b,c)\dsabc=
          \\=
 \Gamma(a+b+n)\Gamma(a+c+n)\Gamma(b+c+n)\, n!\,\delta_{m,n}
 \end{multline}

 {\sc Lemma \fact.} {\it The image of the function
$$
 \Bigl(\frac x{x+1}\Bigr)^n(1+x)^{-a-b}
$$
under the index hypergeometric transform $J_{b,c}$
is
$$\frac{\Gamma(a+is)\Gamma(a-is)} {\Gamma(a+b+n)\Gamma(a+c+n)}
        \SS_n\bigl(s^2;a,b,c\bigr)$$
and the image of the function $z^n\in W^a_{b,c}$
under $J^a_{b,c}$ is }
$$\frac {\SS_n\bigl(s^2;a,b,c\bigr)} {\Gamma(a+b+n)\Gamma(a+c+n)}
        $$
 The Lemma can be checked by a more or less direct calculation
 (see \cite{Ner8}).
 The second part of the Lemma implies Theorem 5.2,
 since we have an explicit correspondence of
 the orthogonal bases
 ($z^n$  and the Hahn polynomials).

 \bigskip

 {\large\bf \sect Helgason transform and spherical transform}

 \medskip

 Here we discuss the spherical transform
 (it is also is called the Harish-Chandra transform)
 and the Helgason transform.
  The latter is used only for an explanation
 of the former. In \Subsections 6.9--6.10
 we dogmatically define the spherical transform
 independently on the Helgason transform.
 Proofs of all facts
  on the spherical representations and the spherical
  transform   formulated in \Subsections
  6.1--6.8 are contained in
  Helgason \cite{Hel2}, chapter 4.

 \smallskip

 {\bf\punct.  Definition of spherical representations.}
  An irreducible representation
  in a complete
separable locally convex space (see \cite{Kir}, 7.2)
   $\rho$ of the group $\G$ in a space $W$ is called a
  {\it spherical representation}
   if $W$ contains a $\K$-invariant vector
   (this vector is called a {\it spherical vector}).

  Under some minor natural conditions on the space and
representation,
  a spherical vector  is unique up to a factor
(see \cite{Hel2}, IV.4 and references in this book).

Denote by $\xi$ the spherical vector
 of a spherical representation $\rho$.
 Consider the operator in the space $W$ given by
 \begin{equation}
 \varPi:=\int_\K\rho(k)\,dk
 \label{6.1}
 \end{equation}
where $dk$ is the Haar measure on $\K$ such that
the measure of the whole group is 1.
Obviously (see \cite{Kir}, 9.2.1)
 $\varPi$ is the $\K$-intertwining projection
 to the vector $\xi$.

The vector $\varPi\rho(g)\xi$ has the form $\Phi(g) \xi$,
where $\Phi(g)\in\C$ is a scalar.
The function $\Phi(g)$ is called the {\it spherical function}
of the representation $\rho$.

Let $k_1, k_2\in \K$. Then
$$\Phi(k_1gk_2)\xi=
\varPi\rho(k_1gk_2)\xi=\varPi\rho(k_1)\rho(g) \rho(k_2)\xi=
\rho(k_1)\varPi\rho(g) \xi =\rho(k_1) \Phi(g)\xi=\Phi(g)\xi
$$
This implies the $\K\times \K$-invariance
of the spherical function:
$$\Phi(k_1gk_2)=\Phi(g)$$
Hence we can consider a spherical function as

   a) a
    function
       on double cosets  $\K\setminus\G/\K$

   b) a $\K$-invariant function on $\B_{p,q}=\G/\K$

   As in 2.9, we shall consider
 $\K$-invariant functions on
 $\B_{p,q}=\G/\K$ as functions depending
on  variables $x_j$.
 Recall that $x_j$ are the eigenvalues
of $z^*(1-zz^*)^{-1}z$.

\smallskip

 {\bf\punct.  Construction of spherical representations.}
  As in 2.1, consider
 the space  $\C^p\oplus\C^q$ equipped with the form $Q$
 having a matrix
 $\bigl(\begin{smallmatrix}1&0\\0&-1\end{smallmatrix}\bigr)$.
  A subspace
 $R\subset \C^p\oplus\C^q$
 is called {\it isotropic} if the form
 $Q$ is the identical zero on $R$. Denote by
 $\IGr_m$ (the {\it isotropic Grassmannian})
 the set of all isotropic $m$-dimensional subspaces
 in   $\C^p\oplus\C^q$.  The space $\IGr_m$ is a
 $\G$-homogeneous space. It is  also
 $\K$-homogeneous, hence there exists a unique
 $\K$-invariant measure
 (or a volume form) on $\IGr_m$.%
 \footnote{Proof. Consider an arbitrary positive
 volume form $\Omega$ on $\IGr_m$.
  Then the  average of $\Omega$ over
  the group $K$ is an invariant volume
  form}
 For $g\in\G$ we denote by $j_m(g,R)$
 the Jacobian of the transformation $R\mapsto gR$.

 An {\it isotropic flag} in  $\C^p\oplus\C^q$ is a family
 of isotropic subspaces
 $${\cal R}: R_1\subset R_2\subset\dots\subset R_p$$
 such that $\dim R_j=j$. Denote by $\Fl$ the space of all
 isotropic flags in $\C^p\oplus\C^q$. The space $\Fl$
 is a $\G$-homogeneous space.

 The stabilizer of a flag in $\G$ is
a  {\it minimal parabolic subgroup} in $\U(p,q)$.
 Denote by  $e_1,\dots , e_{p+q}$ the standard basis in
$\C^p\oplus\C^q$. Consider the isotropic subspace $S_j$ spanned
by the vectors $e_i+e_{p+i}$ for $i\le j$.
Consider the flag
 $${\cal S}:\quad S_1\subset S_2\subset \dots\subset S_p$$
We define the {\it standard parabolic subgroup}
 $\P\subset \G$ as the stabilizer of the flag
$\cal S$ in $\G$.

The space $\Fl$
  is also $\K$-homogeneous,
  \begin{equation}
  \Fl=\K/\M,\qquad \mbox{where}\,\,
  \M\simeq\U(q-p)\times
  \underbrace
  {\U(1)\times \dots\times\U(1)}_{p\,\,\,\mbox{times}}
  \label{6.2}
  \end{equation}
   hence
  there exists a unique (up to a factor)
 $\K$-invariant measure on $\Fl$.
 For $g\in\G$ we denote by $J(g,{\cal R})$
 the Jacobian of the transformation ${\cal R}\mapsto g{\cal  R}$.

 A {\it multiplier} $\omega(g,x)$ on a homogeneous
 space $G/H$ is a function on $G\times\, (G/H)$ satisfying
 the condition
 $$\omega(g_1g_2,x)=\omega(g_1,g_2x)\omega(g_2,x)$$

{\sc Example.} The Jacobian is a multiplier. In particular,
$J(g,{\cal R})$ and $j_m(g,R_m)$, $m=1,\dots,p$
are multipliers on $\Fl$.

\smallskip

 Let $\omega$  be a multiplier. Then the formula
 $$\rho_\omega(g)f(x)=f(g(x))\omega(g,x)$$
 defines a representation of $G$ in the space of functions
 on $G/H$.
 For a description of all multipliers in general case see
 \cite{Kir}, 13.2, see also \cite{Kir}, 13.5 for geometric explanations.

Consider the action of $\G$ on the homogeneous space
$\Fl$. Obviously, any function of the form
$$\omega(g,{\cal R})=\prod_{m=1}^p j_m(g,R_m)^{\beta_m}$$
is a multiplier on $\Fl$.   It is more convenient
to define the "basic" multipliers by the formula
\begin{align*}
\omega_1(g,{\cal P})&=j_1(g,P_1)^{\frac{1}{2(q-p+1)}}; \\
\omega_m(g, {\cal P})&=
j_m(g,P_m)^{\frac{1}{2(q-p+m)}}
j_{m-1}(g,P_m)^{-\frac{1}{2(q-p+m-1)}}
;\qquad \mbox{for} \,\,\, m=2,\dots,p
\end{align*}
In this notation,
$$J(g,{\cal P})=\prod_{m}\omega_m(g,P_m)^{2(q-p+1-2m)}$$
Fix $s_1,\dots,s_p\in\C$. We define the representation
$\widetilde \rho_s$ of the group $\G$ in the space
of $C^\infty$-smooth
functions on $\Fl$ by the formula
\begin{equation}
\widetilde \rho_s(g)f({\cal P})=f(g{\cal P})J(g,{\cal P})^{1/2}
      \prod \omega_j(g,P_j)^{s_j}
\label{6.3}
\end{equation}
For  $(s_1,\dots,s_p)$ in a general position,
 the representation
$\widetilde \rho_s$  is irreducible.
Otherwise there exists a unique irreducible
spherical subquotient in $\widetilde\rho_s$.
Let us describe it.

Denote by $M$ the minimal closed
$\widetilde \rho_s$-invariant subspace of $C^\infty(\Fl)$
containing the function $f({\cal P})=1$. Denote by $N$
the maximal proper closed
$\widetilde \rho_s$-invariant subspace in $M$
\footnote{The space $M$ is the cyclic span of
the vector $1$. Consider a proper $\G$-invariant subspace
$L\subset M$. Then $1\notin L$. Hence $L$
is contained in the kernel of the
$\K$-invariant projection $\varPi$ (see (6.1)).
Thus a sum of all proper subspaces in $M$
is contained in $\ker \varPi$.}
 and
consider the quotient $M/N$. We denote
by $\rho_s$ the representation of $\G$ in the space
$M/N$ (if $\widetilde\rho_s$   is reducible,
then $\rho_s\simeq\widetilde\rho_s$).

\smallskip

Consider the space $\C^p$ with coordinates
$s_1,\dots,s_p$.
Consider the {\it hyperoctahedral group} $D_p$, i.e.,
 the group
generated by all permutations of coordinates and by
the reflections
 $(s_1,\dots,s_p)\to(\sigma_1 s_1,\dots,\sigma_p s_p)$,
 where $\sigma_j=\pm 1$.  This group also coincides with
the so-called {\it restricted Weyl group} of $\U(p,q)$.

 \smallskip

{\sc Theorem \fact.}
a) {\it The representations $\rho_s$ are exactly
all spherical representations of $\G$. More precisely,
the representations
$\rho_s$ are the all spherical Harish-Chandra modules up to
equivalence of Harish-Chandra modules.}

b) {\it The representations $\rho_s$ and $\rho_{s'}$
are equivalent iff there exists $\gamma\in D_p$
such that $\gamma s=s'$.}

 Denote by $\dual$ the set of all
 unitary spherical representations.

 Explicit description of the set $\dual$ is unknown.
 We shall formulate two facts about $\dual$ that are necessary
 for understanding of the subsequent text.

 \smallskip

1. Obviously, if all the coordinates $s_j$ are pure imaginary,
then the representation
$\widetilde \rho_s$ is unitary in $L^2(\Fl)$.
These representations are called {\it representations of
the spherical unitary principal nondegenerate series.}

\smallskip

2. If a representation $\rho_s$ is unitary,
 then  $s_j\in \R\cup i\R$
 (this a corollary of Theorem 6.1.b; indeed
a unitary representation $\rho$ is equivalent to
 its contragredient (= dual) representation,
 and the representation dual to $\rho_s$
 is $\rho_{-\bar s}$).

 \smallskip

 {\bf\punct. Another realization of
 spherical representations.}
 Let $h$ be an element of the standard parabolic subgroup $\P$.
 Then $h$ induces a linear transformation
 in each $1$-dimensional quotient $S_j/S_{j-1}$,
 clearly, it is the multiplication by some
 complex number $\chi_j(g)$\footnote{certainly
 $\chi_j(g)$ are the diagonal elements of the matrix $h$ in
 the basis  $e_1+e_{p+1}$, \dots, $ e_p+e_{2p}$, $e_{2p+1},\dots$,
 $e_q$, $ e_p-e_{2p}$, \dots, $ e_1-e_{p+1}$.}.

 Let $s_j\in\C$.
 Consider the space $L_s$ of all smooth functions on $\G$
 satisfying the condition

 \begin{equation}
 F(rh^{-1})=F(r)\prod_{m=1}^p |\chi_m(h)|^{(q-p+1-2m)+s_m}
 \qquad r\in \G, h\in \P
 \label{6.4}
 \end{equation}
 Denote by $\widehat\rho$ the representation of $\G$
 in $L_s$ given by
 \begin{equation}
 \widehat\rho(g)F(r)=F(gr)
 \label{6.5}
 \end{equation}

 {\sc Lemma \fact.} {\it Any element of $\G$ has
 the decomposition $g=hk$, where $h\in\P$, $k\in\K$.}

 \smallskip

 {\sc Proof.} This is equivalent to the transitivity
 of $\K$ on $\Fl=\G/\P$.\kvadrat

 \smallskip

 Thus any function
 $F\in L_s$ is determined by its restriction
  to the submanifold $\K\subset\G$.
 In (\ref{6.2}) we defined the subgroup $\M:=\K\cap\P$.
For $h\in\M$, we have $|\chi(h)|=1$. Thus
$F$ can be regarded as a function on $\K/\M=\Fl$.

Thus we obtain a canonical operator
from $L_s$ to the space of  $C^\infty(\Fl)$,
it is not hard to check, that this operator intertwines
$\widehat\rho$ with $\widetilde\rho$.

\smallskip

{\sc Remark.}
This construction also explains the appearance of  the
"canonical multipliers" $\omega_m$, they corresponds to
the $1$-dimensional
characters $|\chi_m(h)|$ of $\P$.
\hfill $\square$

 {\bf\punct. Preliminary remarks on
 Plancherel formula.}

     By a general abstract
 theorem (see \cite{Kir}, 4.5, 8.4), any unitary
 representation of a
 locally compact group is a direct integral of irreducible
 representations. For some types of groups
 (including semisimple groups and, hence, $\U(p,q)=\U(1)\times{\rm SU}(p,q)$)
 this decomposition is unique in a natural sense.

 \smallskip

 {\sc Theorem \fact.}
a) {\it Representation of $\G$
 in $L^2(\G/\K)$ is a direct integral over  $\dual$
  with multiplicities
 $\le 1$.}

b) {\it A kernel representation $\T_\alpha$ is a direct
integral over $\dual$ with multiplicities $\le 1$.}

 \smallskip

 {\bf\punct. Proof of Theorem 6.3.}  Denote by
 ${\cal M}(\G)$ the algebra of   compactly
 supported (complex-valued) measures on $\G$;
 the multiplication
  in  ${\cal M}(\G)$  is the usual convolution $*$.
    We define by $\mu^\square$ the pushforward of
    a measure $\mu$  under the map $g\mapsto g^{-1}$.
   Obviously,  $\mu\mapsto\mu^\square$ is an involution
   on ${\cal M}(\G)$:
  $(\mu*\nu)^\square=\nu^\square * \mu^\square$

Denote by ${\cal C}(\G)$ the  subalgebra of ${\cal M}(\G)$
consisting of measures invariant with
respect to the transformations
$g\mapsto k_1gk_2$, where $k_1,k_2\in \K$.

  \smallskip

 {\sc Theorem \fact.} (Gelfand, see \cite{Hel2}, 5.1)

   a) {\it $\mu^\square=\mu$ for all $\mu\in {\cal C}(\G)$.}

    b) {\it The  algebra ${\cal C}(\G)$ is commutative.}

    \smallskip

 {\sc Proof.}
a) It easily can be checked that $g^{-1}$ is contained
in the double coset $\K g\K$. Hence $\mu=\mu^\square$.

b) By a), the identity
$(\mu*\nu)^\square=\nu^\square * \mu^\square$
  coincides with
$\mu*\nu=\nu*\mu$.                        \kvadrat

\smallskip

 Let $\kappa$
 be a unitary representation of $\G$
 in a Hilbert space $H$.
 For $\mu\in{\cal M}(\G)$ we define the operator
 $$\kappa(\mu)=\int_\G \kappa(g)\,d\mu(g)$$
 It is easily shown that
 $\kappa(\mu*\nu)=\kappa(\mu)\kappa(\nu)$ (see \cite{Kir}, 10.2).

 Denote by $H^\K$ the space of all $\K$-invariant
 vectors in $H$. The projection $\varPi$ to
 $H^\K$ is given by
 $$\varPi=\int_\K \kappa(k)dk$$
 Obviously, $\varPi$
  has the form $\kappa(\delta_\K)$,
 where $\delta_\K$ is the  Haar measure on $\K$
normalized by the condition: the measure of the whole group
is 1.
 We consider this measure  as a $\delta$-measure
 on $\G$
 supported by $\K$).

 \smallskip

{\sc Lemma \fact.}
a)  {\it The operators $\kappa(\mu)$, where $\mu\in{\cal C}(\G)$,
 are zero on the
orthocomplement to the  space $H^\K$.}

b) {\it The operators $\kappa(\mu)$, where $\mu\in{\cal C}(\G)$,
    are self-adjoint.}

c) {\it Let $R$ be an ${\cal C}(\G)$-invariant subspace
in  $H^\K$. Let $Z$ be the $\G$-cyclic span of
$R$. Then the projection of $Z$ to $H^\K$ coincides
with $R$.}

d) {\it If $\kappa$ is irreducible, then $\dim H^\K\le 1$.}

e) {\it If there exists a cyclic $\K$-fixed   vector $\Xi$ in $H$,
then the spectrum of $\kappa$  contains
only spherical representations and
their multiplicities are $\le 1$.}

\smallskip

{\sc Proof.}
a) $\kappa(\mu)=
\kappa(\mu*\delta_\K)=\kappa(\mu)\varPi$ for $\mu\in {\cal C}(\G)$.

\smallskip

b) We have $\kappa(\mu)^*=\kappa(\mu^\square)$.
 By Theorem 6.4.a, this coincides with $\kappa(\mu)$.

\smallskip

c) Let $h\in R$, then
$$\varPi \kappa(g)h=\varPi \kappa(g)\varPi h=
\kappa(\delta_\K)\kappa(g)\kappa(\delta_\K)h=
\kappa(\delta_\K*\delta_g*\delta_\K)h\in R$$

d) The algebra ${\cal C}(\G)$ is commutative.
Any commutative family of self-adjoint bounded
operators in a Hilbert space of dimension $>1$ has
a proper invariant subspace. Now apply c).

e) First let $L$ be a $\G$-invariant subspace without
$\K$-invariant vectors. Then $\Xi$ is contained in
its orthocomplement $L^\perp$.
But $\Xi$ is cyclic and hence $H=L^\perp$.

Second, a cyclic representation of a commutative
$*$-algebra has multiplicities $\le 1$
(see \cite{Kir}, 4.4, Problem 3).
\kvadrat

\smallskip

{\sc Proof of Theorem 6.3} b)
The distinguished vector $\Xi_\alpha$
is cyclic (see 4.4). It remains to apply
  Lemma 6.5.e.
 \kvadrat

\smallskip

{\sc Proof of Theorem 6.3.a.}
There exists a unique $\G$-invariant differential operator  $\Delta$
of order 2 on $\G/\K$. Consider the heat equation
$$\frac d{dt}f(t,z)=\Delta f(t,z); \qquad t\ge 0$$
 Let $N(t;z,u)$ be the corresponding
heat kernel and let $A_t$ be the corresponding evolution operator, i.e.,
$$
f(t,z)=A_t f(0,z)= \int_{\B_{p,q}} N(t;z,u)f(u)\det(1-uu^*)^{-p-q}\{du\}
$$

Fix $t>0$.
Obviously, the function $N(t;z,0)$ is a cyclic vector
in the (closed) image of the evolution operator $A_t$.
Thus the representation of $\G$ in the subspace $\Im A_t$
has a multiplicity free spectrum.

We have $A_tA_{t'}=A_{t+t'}$. Hence $\Im A_t\supset \Im A_\tau$
if $t <\tau$.  Thus $L^2(\G/\K)$ is the closure of the union
of increasing family of multiplicity free subrepresentations
$$\Im A_1\subset\Im A_{1/2}\subset\Im A_{1/3}\subset\dots$$
Thus $L^2(\G/\K)$ itself is multiplicity free.



 \smallskip

 {\bf\punct. Normalization of  Plancherel measure.}
 Let $\rho_s\in\dual$ be a spherical representation.
  Denote by $W_s$ the space
 of the representation $\rho_s$,
   denote by
 $\langle\cdot,\cdot \rangle_s$  the scalar product
 in $W_s$,
  denote by $\xi_s$
 a spherical vector in $W_s$, such that
 $\langle\xi_s, \xi_s\rangle_s=1$.

We shall consider Borel measurable functions
$f$ on  $\dual$
such that a value of $f$ at a point $s$ is an element of the
space $W_s$ (we omit a definition of measurability).

Consider a Borel measure  $\nu$  on the set
$\dual$ of spherical unitary representations
of $\G$, let $M$ be the  support of $\nu$.
We shall consider the Hilbert space
$\int W_s d\nu(s)$
of all Borel measurable functions
$f:s\mapsto f(s)\in W_s$
satisfying the condition
$$\int_M\langle f(s), f(s)\rangle_s d\nu(s)<\infty$$
The scalar product in this space is defined by
$$[f,g]=\int_M\langle f(s), g(s)\rangle_s d\nu(s)$$
The group $\G$ acts in our space pointwise.
We denote this representation
({\it direct integral})
by
$\int \rho_sd\nu(s)$ (for details of the definition of  direct integral
see \cite{Dix}).

In this Section we shall
present the classical description
of the measure $\nu_\infty(s)$ on  $\dual$ such that
$\int \rho_s(g)d\nu_\infty(s)$
 is equivalent to $L^2(\G/\K)$
and of a canonical unitary operator $U$ from
 $L^2(\G/\K)$ to the direct integral  $\int W_s d\nu(s)$.

In the next Section we  solve the same
problem for the kernel representations
of $\G=\U(p,q)$.
These measures
 are called the {\it Plancherel measures.}

\smallskip

{\sc Remark.}
A Plancherel measure  is not canonically defined
(see \cite{Kir}, 4.5, 8.4).
Let $\beta(s)$ be a positive function on
the support $M$ of the measure $\nu$.
Clearly, the direct integrals
$\int \rho_s(g)d\nu(s)$  and
 $\int \rho_s(g)\bigl[\beta(s)d\nu(s)\bigr]$
are equivalent. \hfill $\square$

\smallskip

For the case $L^2(\G/\K)$, we define
 a canonical normalization  of a Plancherel measure
 by the following rule: {\it the element
 $s\mapsto \xi_s$ of the direct integral
  corresponds to the  $\delta$-function on $\G/\K$
  supported by $z=0$.}

   Let us repeate this more formally.
  Consider  an element $s\mapsto a(s)\xi_s$
  of the direct integral. Let $F_a$ be the  image
  of  $s\mapsto a(s)\xi_s$
  under the operator $U^{-1}$. Then we require
  $F_a(0)=\int a(s)\,d\nu_\infty (s)$ for all $a(s)$.

  This normalization is not  an arbitrary rule.
  If we shall change it, then formulas for the measure $\nu_\infty$
  and for the operator $U$ will change. It is easy to
  obtain a more complicated formula in this way, but
  it is impossible to obtain a simpler formula.

  Now we shall describe the map $U$ (the Helgason transform)
  and the measure $\nu_\infty$ (the Gindikin--Karpelevich measure).

  \smallskip

{\bf\punct. Helgason transform and Gindikin--Karpelevich
measure.}
 Consider the "section of  wedge"
model $\SW_{p,q}$ of the space  $\G/\K=\U(p,q)/\U(p)\times \U(q)$,
see \Subsection 2.8.
Fix  $s:=(s_1,\dots, s_p)\in \C^p$, assume $s_{p+1}=0$.
 Let
 $Z=\bigl(
 \begin{smallmatrix}1&0\\2K&L\end{smallmatrix}\bigr)
 \in\SW_{p,q}$.   Denote by $[Z]_k$ the left upper $k\times k$ block
of $Z$.
 We define the functions $\Psi_s(Z)$ by
\begin{equation}
\Psi_s\begin{pmatrix}1&0\\2K&L\end{pmatrix}=
\prod_{j=1}^p
\det \left[\begin{array}{cc}
1&K^*\\K&\frac12(L+L^*)
\end{array}\right]_{q-p+j}^%
{-\sigma_j+(s_j-s_{j+1})/2}
\label{6.6}
\end{equation}
where
\begin{equation}
\sigma_1=\dots=\sigma_{p-1}=1;\qquad \sigma_p=q-p+1
\end{equation}

\smallskip

{\sc Lemma \fact.} {\it
For any $h$ in the standard parabolic subgroup $\P$, }
\begin{equation}
\Psi_s(Z^{[h]})=\Psi_s(Z)%
\prod_{m=1}^p |\chi_m(h)|^{(q-p+1-2m)+s_j}
\label{6.8}
\end{equation}

{\sc Proof.} A calculation, see \cite{Ner5}. \hfill $\boxtimes$

\smallskip

Let $f$ be a $C^\infty$-function on $\SW_{p,q}$
with a compact support.
Its {\it Helgason transform} is the function
$$
F(r;\, s_1,\dots, s_p);
      \qquad r\in\G,\,s_j\in\C
$$
defined by
\begin{equation}
F(r;\, s_1,\dots, s_p)=
\int_{\SW_{p,q}}f(Z)\Psi_s(Z^{[r]})
 \det \left[\begin{array}{cc}
1&K^*\\K&\frac12(L+L^*)
\end{array}\right]^{-p-q} \{dZ\}
\label{6.9}
\end{equation}

By (\ref{6.8}), the functions $F(r;s_1,\dots,s_p)$
satisfy (\ref{6.4}), i.e., for a fixed $s\in\C^n$
a function $F(r,s)$ is an element of the space
$L_s$ of the representation $\widehat\rho_s$ (see (\ref{6.3}).

Obviously, the
 pushforward of the transformation
$f(Z)\mapsto  f(Z^{[g]})$
under the Helgason transform is given by
(\ref{6.5}).

Above we identified the representations $\widehat \rho_s$
and $\widetilde \rho_s$. Thus we can consider a function $F$
as a function on $\Fl\times\C^p$.

\smallskip

 {\sc Remarks.} 1) Functions $F(\dots)$ are holomorphic
 in $s_1,\dots,s_p\in\C^n$.\hfill $\square$

  2)
  For each element $\gamma$ of the hyperoctahedral
  group $D_p$ there exists an integral
  operator $A_\gamma$ on $\Fl$ (which can be written  explicitly)
   such
  that $F({\cal P},\gamma s)=A_\gamma F( {\cal P},s)$.
                      \hfill $\square$

  Denote by $\Sigma_p$ the simplicial cone
  $$\Sigma_p: s_1\ge s_1\ge \dots \ge s_p\ge 0$$
  Denote by $i\Sigma_p\subset \C^p$
   its image under the multiplication by $i$.

   \smallskip

{\sc Theorem \fact.} {\it The Helgason transform
 is a unitary operator
\begin{multline*}
L^2\Bigl(\SW_{p,q},
\det \left[\begin{array}{cc}
1&K^*\\K&\frac12(L+L^*)
\end{array}\right]^{-p-q} \{dZ\}\Bigr)
\to
L^2\Bigl(\Fl\times (i\Sigma_p),\,\,\, {\frak R}(s)\,
 \,\{ds\}\, d{\cal P}\Bigr)
 \end{multline*}
 where $\{dZ\}$ is the Lebesgue measure
on $\SW_{p,q}$, $d\cal P$ is the $\K$ -invariant
 measure on $\Fl$, $\{ds\}$ is the Lebesgue measure
 on the set $i\Sigma_p$
 and ${\frak R}(s)$ the Gindikin--Karpelevich density
 \begin{equation}
 \prod_{k=1}^p\frac{|\Gamma (\tfrac12 (q-p+1)+s_k)|^4}
             {|\Gamma(2s_k)|^2}
                  \prod_{1\le k<l\le p}(s_k^2-s_l^2)^2
 \label{6.10}
 \end{equation}
 The pushforward of the transformation $f(Z)\mapsto f(Z^{[g]})$
 under the Helgason transform is given by}
 $$F({\cal P};s_1,\dots,s_p)\mapsto
   F(g{\cal P}; s_1,\dots,s_p)J(g,{\cal P})^{1/2}
      \prod \omega_j(g,P_j)^{s_j}
 $$

{\bf \punct. Spherical transform: preliminary remarks.}
 Consider the space of  $\K$-invariant functions $f$ on $\G/\K$.
As above, we can consider elements of this space as functions
in the variables $x_j>0$
(see \Subsection 2.9) symmetric with respect to
the permutations.

Consider the Helgason transform $F$
of a $\K$-invariant function $f$.
 Obviously, $F$ is a $\K$-invariant function.
Hence $F$ depends only in the variables $s_1$,\dots, $s_p$.
Thus we obtain the transform
from the space of symmetric functions
in the variables $x_1$,\dots,$x_p$ to the space
 of $D_p$-symmetric functions in the variables
 $s_1,\dots,s_p$.

 If a function $f$ in formula (\ref{6.9})
  is $\K$-invariant, then
 we can replace the factor $\Psi_s$ by its average over
 the group $\K$. But this average is the spherical function
 $\Phi_s$
 of the representation $\rho_s$ (see \cite{Hel2}, 4.4.3).
Thus the Helgason transform on the space
of $\K$-invariant functions is given by
 \begin{equation}
 F(s)=\const\cdot\int_{\R_+^p}
 f(x)\Phi_s(x)w(x)
  dx
 \label{6.11}
 \end{equation}
where
\begin{equation}
w(x)\,dx:=
\prod_{j=1}^px^{q-p}
 \prod_{1\le k<l \le p}(x_k-x_l)^2  \prod_{j=1}^p dx_j
\label{6.12}
\end{equation}
  The map $f\mapsto F$ given by (\ref{6.11})
   is called the {\it spherical transform}.

 {\sc Remark.} In some sense, the spherical transform is
 simpler than the Helgason transform.

   1. We have only $p$ variables instead of $2pq$ variables.

   2. The function $s\mapsto \Phi_s(x)$ is holomorphic
   in $s$  and $D_p$-symmetric with respect to $s$
   (by Theorem 6.1b). Hence for a compactly supported $f$
   its spherical transform $F$ is a $D_p$-symmetric
   holomorphic function.

   3. The index hypergeometric transform $J_{b,c}$
    for some special
   values of the parameters $b,c$ gives
   the spherical transforms for
   all rank 1 groups ($\OO(p,1)$, $\U(p,1)$, $\Sp(p,1)$).
   Heckman and Opdam \cite{HO} constructed the
   family of integral transforms
   interpolating the spherical transforms for
   simple Lie
   groups of arbitrary rank. \hfill $\square$

   \smallskip

{\bf\punct.  Spherical functions
and Berezin--Karpelevich formula.} First we fix the notation
$$\det\limits_{k,l}{c_{k,l}}:=\det
  \begin{pmatrix} c_{11}&\dots&c_{1p}\\
                  \vdots&\ddots&\vdots\\
                  c_{p1}&\dots&c_{pp}
  \end{pmatrix}$$
{\it Below all determinants are of the size $p\times p$.}

\smallskip

{\sc Theorem \fact.}
 (Berezin--Karpelevich \cite{BK}, Hoogenboom\cite{Hoo})
{\it The spherical functions of the group
$\G=\U(p,q)$ are given by}
\begin{equation}
\Phi_s(x)=\const\cdot\frac
{\det\limits_{k,j}
\bigl\{
 \F\Bigl[\begin{array}{c}
\tfrac12(q-p+1)+s_j,\tfrac12(q-p+1)-s_j\\ q-p+1
\end{array};-x_k\Bigr]\bigr\}}
{\prod_{1\le k < l \le p}(s_k^2-s_l^2)
 \prod_{1\le k < l \le p}(x_k-x_l)}
\label{6.13}
\end{equation}
In order to simplify this expression, we introduce the notation
$$\boxed{ r=(q-p+1)/2}$$

{\bf\punct. Spherical transform.}
Thus the spherical transform for the group $\G$
is given by
\begin{equation}
F(s)=\overset{\curvearrowright} f(s)=
\const \int\limits_{x_1\ge 0,\dots,x_p\ge 0}
       f(x) \Phi_s(x)w(x)\, dx
\label{6.14}
\end{equation}
where $w(x)$ is given by (\ref{6.12}) and $\Phi_s$ is given by
 the Berezin--Karpelevich formula.

 \smallskip

 {\sc  Theorem \fact.}  {\it Spherical transform is a unitary
 {\rm(}up to a factor{\rm)} operator from the space of symmetric
 functions with the scalar product
 $$\langle f,g\rangle=\int_{\R^p_+}f(x)\ov{g(x)}w(x)\,dx $$
 to the space of $D_p$-symmetric functions with the scalar
 product
 $$\langle F,G\rangle=\int_{i\R^p}  F(s)\ov{G(s)}
   {\frak R}(s)\,\{ds\}$$
 where ${\frak R}(s)$ is the Gindikin--Karpelevich density
 {\rm(\ref{6.10})} and $\{ds\}$ is the Lebesgue measure on
 $i\R^p$.
 The inversion formula is given by                }
 \begin{equation}
 f(x)= \overset{\curvearrowleft}F(s)
  =\const\int_{i\R^p} F(s) \Phi_s(x) {\frak R}(s)\,\{ds\}
  \label{6.15}
  \end{equation}

 {\sc Remark.}  Theorems 6.7 and 6.9 coincide.
 The implicator 6.7 $\Rightarrow$ 6.9 is obvious.
 Let us explain  $\Leftarrow$. The Helgason
 transform is an
 operator from $L^2(\G/\K)$ to the direct integral of
 the principal nondegenerate series over the Gindikin--Karpelevich
 measure. We must show that this operator is unitary.
 Assume that the spherical transform is unitary. Then
 the Helgason transform preserves the scalar products
 $\langle g v,w\rangle$, where $v,w$ range
 the space of $\K$-fixed vectors and $gv$ denotes
 an action of the group on the Hilbert space.
 Thus the Helgason transform preserves  the
 scalar products
 $\langle g_1v,g_2w\rangle=\langle g_2^{-1} g_1v,w\rangle$
 and hence the Helgason transform is unitary
 (since the vectors $gv$ span the both Hilbert spaces).
 \hfill $\square$

\smallskip

 {\bf \punct. Some integrals with determinants.}

 \smallskip

 {\sc Lemma \fact.} {\it Let $\mu$  be a measure on $\R$.
 Then
 \begin{multline}
 \int_{\R^p}\det\limits_{k,l}\{f_k(x_l)\}
            \det\limits_{k,l}\{g_k(x_l)\}
            d\mu(x_1)\dots d\mu(x_p) 
=n!\det\limits_{k,m}
\Bigl\{\int_\R f_k(x)g_m(x)\,d\mu(x)\Bigr\}
\label{6.16}
\end{multline}
if the right-hand side of the equation has sense.}

{\sc Proof.} Obvious.                                \kvadrat

\smallskip

{\bf\punct. Some spherical transforms.}
{\sc Lemma \fact.} {\it For  functions $\beta_1$, \dots, $\beta_p$
 on $\R_+$ we define
the function
\begin{equation}
\Theta_\beta(x):=
\frac{\det\limits_{k,l}\{\beta_k(x_l)\} }
                  {\prod_{k<l} (x_l-x_k)}
\label{6.17}
\end{equation}
                  Then its spherical transform is   }
$$\overset{\curvearrowright}
\Theta_\beta(x) =
\frac {1}{\prod\limits_{k<l}(s_k^2-s_l^2)}
\det\limits_{k,l}\Bigl\{
 \int\limits_0^\infty \beta_k(x)
 \F(r+s_l,r-s_l;
    2r;-x)x^{q-p}dx\Bigr\}$$

We see that the integral in the  curly brackets  is
the index hypergeometric transform  $J_{r,r}\beta(s)$.

\smallskip

{\sc Proof.} We evaluate the integral (\ref{6.14})
by Lemma 6.10.\kvadrat

\smallskip

{\sc Corollary \fact.} {\it The spherical transform of the
 function $\prod_j b(x_j)$ is}
$$\frac 1{\prod_{k<l}(s_k^2-s_l^2)}
\det\limits_{k,l}\Bigl\{
 \int_0^\infty
 x^{k-1}b(x) \F(r+s_l,r-s_l;
    2r;-x)x^{q-p}dx\Bigr\}$$

Our next Section is based on the following simple
formula

\smallskip

{\sc Theorem \fact.} {\it The spherical transform of the function
 $$\det(1-zz^*)^\alpha=\prod\nolimits_{j=1}^p(1+x_j)^{-\alpha}$$
is}
\begin{equation}
\frac{\prod_{k=1}^p \Gamma(\alpha-\tfrac12(q+p-1)+s_k)
                      \Gamma(\alpha-\tfrac12(q+p-1)-s_k)}
       {\prod_{j=0}^{p-1}\Gamma^2(\alpha-j)}
\label{6.18}
\end{equation}

Below we shall use the notation
$$\boxed{h:=(q-p+1)/2}$$

       {\sc Proof.} We must evaluate
\begin{multline*}
\int\limits_{\R^p_+}
\prod_{j=1}^p(1+x_j)^{-\alpha}
\frac{\det\limits_{k,j}
\bigl\{\F(r+s_k,r-s_k; 2r;-x_j)\bigr\}}
{\prod_{1\le k < l \le p}(s_k^2-s_l^2)
 \prod_{1\le m < n \le p}(x_m-x_n)} \times\\ \times
 \prod_{1\le m < n \le p}(x_m-x_n)^2
 \prod_{k=1}^p x_k^{q-p} dx_1\dots dx_p
 \end{multline*}
 It is possible to apply directly  Corollary 6.12,
 but it is more convenient to write
 $$\prod\limits_{1\le k<l\le p}(x_k-x_l)=\prod\limits_k(1+x_k)^{p-1}
 \prod\limits_{1\le k<l\le p}
 \Bigl(\frac {x_k}{1+x_k}- \frac {x_l}{1+x_l}\Bigr)=
  \prod_k(1+x_k)^{p-1} \det\limits_{k,m}
  \Bigl(\frac {x_k}{1+x_k}\Bigr)^{m-1}$$
 By Lemma 6.11,  we obtain
 $$\frac{1}  {\prod(s_k^2-s_l^2) }
 \det\limits_{m,l}\Bigl\{\int\limits_0^\infty
   (1+x)^{-\alpha+p-1} \Bigl(\frac {x}{1+x}\Bigr)^{m-1}
   \F(r+s_l,r-s_l;2r;-x)x^{q-p}dx\Bigr\}
  $$
 By Lemma 5.3, the integrals under the determinant
  are the continuous dual Hahn polynomials $\SS_{m-1}$.
 We obtain
 $$
 B(\alpha)\cdot \frac {\det\limits_{k,m}
  \bigl\{\SS_{m-1}(s_k^2; \alpha-h,r,r)\bigr\} }
   {\prod(s_k^2-s_l^2) }
   $$
   where $B(\alpha)$ is given by (\ref{6.18}).
By (\ref{5.6}), $\SS_{m-1}(s^2)=(s^2)^{m-1} +\dots $
 and hence the determinant in the numerator is
 the Vandermonde determinant.\kvadrat

 \smallskip

 {\bf \punct. Image of the canonical basis
  under the spherical transform.}
The canonical orthogonal basis
 $\Delta_\mu$ in the space $\V^\K_\alpha$
 (see \Subsection 4.6) in our
 coordinates $x_k$ is given by
 $$\Delta_\mu(x)=\prod_{k=1}^p (1+x_k)^{-\alpha+p-1} \,\,\cdot\,\,
 \frac
 {\det\limits_{k,j}
 \Bigl\{\Bigl(\frac {x_k}{x_k+1}\Bigr)^{\mu_j+p-j}\Bigr\} }
 {\prod\limits_{1\le k< m\le p} (x_k-x_m)}
 $$

 {\sc Theorem \fact.} {\it The image of the function
 $\Delta_\mu$ under the spherical transform is given by
 \begin{equation}
 \const  \cdot
 \frac{\prod\limits_{k=1}^p \Gamma(\alpha-h+s_k)
                      \Gamma(\alpha-h-s_k)}
       {\prod\limits_{j=0}^{p-1}\Gamma^2(\alpha-j)}\,\cdot \,
\frac{\det\limits_{k,l}
\bigl\{\SS_{\mu_l+p-l}(s_k^2;\alpha-h,r,r)\}}
     {\prod\limits_{1\le k < l \le p}(s_k^2-s_l^2)}
\label{6.19}
\end{equation}
 where $\SS_n$ are the continuous dual Hahn polynomials} (\ref{5.6}).

 \smallskip

 {\sc Proof.} We apply  Lemma 6.10 and Lemma 5.3.\kvadrat

 \smallskip

 {\bf\punct. Reduction
 of the Gindikin--Karpelevich inversion
 formula to the Berezin--Karpelevich formula.}
  Let
 $\Theta_f(x)$ be the same as above (\ref{6.17}).
 Then
 $$\overset{\curvearrowright}\Theta_f(s)=
 \det\limits_{k,m}\bigl\{J_{r,r}f_m(s_k)\bigr\}
 /\prod(s_k^2-s_l^2) $$
 We must check the equality
 $$
 \langle  \Theta_f, \Theta_g \rangle=
 \langle  \overset{\curvearrowright}\Theta_f ,
    \overset{\curvearrowright}\Theta_g  \rangle
    $$
  By Lemma 6.10, this reduces to
 $$   \det\limits_{k,m}
 \bigl\{\int_0^\infty f_m(x)\ov {g_k(x)}x^{q-p}dx\bigr\}=
      \det\limits_{k,m}\bigl\{\int_0^\infty (J_{r,r}f_m)(s)
                            \ov{  (J_{r,r}g_k)(s)}
                              \frac{|\Gamma(r+is)|^4}
                          {|\Gamma(2is)|^2}ds\bigr\}
  $$
  By Theorem 5.1, the matrix elements of these two matrices
  coincide.

 {\small
 {\bf\punct.  Comments.} 1)
 The main facts concerning the Helgason transform,
the spherical functions and the spherical transform are the same
 for all semisimple groups.
 But generally,
the spherical functions of the real semisimple groups
 are certain multivariate special functions,
 they are one of the natural multivariate analogues of
 the
 hypergeometric functions (see \cite{HO}).
  Simple determinant
 formulas for spherical functions
  exist only  for the complex   groups \cite{GN}
 and $\U(p,q)$.

 \smallskip

 2) In general, the image of the canonical basis under
 the spherical transform consists of multivariate
 continuous  dual Hahn polynomials or multivariate
  Meixner--Pollachek (for series $\GL_n$)
 polynomials. In our case, the multivariate
 Hahn polynomials admit a simple expression.

   Koornwinder \cite{Koo2} constructed
a multivariate analogue of the Askey--Wilson polynomials. All classical
 and neoclassical orthogonal polynomials in one variable
 can be obtained by a degeneration of the Askey--Wilson polynomials
 (see the treatise \cite{KS}). Multivariate versions
 of (neo)classical orthogonal polynomials can be obtained
 by a degeneration of the Koornwinder construction
 (the basic steps were done in \cite{DS}, \cite{Die}).
 In particular, the multivariate
 continuous  dual Hahn polynomials or Meixner--Pollachek
 polynomials can be constructed in this way.
        }

 \bigskip

 {\large\bf\sect Plancherel formula for kernel representations}

 \nopagebreak

 \medskip
 We preserve the notation
 $$r=(q-p+1)/2;\qquad h=(q+p-1)/2$$

We preserve the notation $\rho_s$ for a spherical
 representation, $W_s$ for its space,
  $\langle\cdot,\cdot\rangle_s$
 for the scalar product in $W_s$, $\Phi_s(g)=\Phi_s(x)$ for
 the spherical function of $\rho_s$.

 \smallskip

 {\bf\punct. Normalization of the Plancherel
 measure.}
By Theorem 6.3b, any kernel representation
$\T_\alpha$ is equivalent
to a multiplicity-free
 direct integral  $\int\rho_sd\nu_\alpha(s)$
of spherical representations
 over some measure $\nu=\nu_\alpha$ on the space $\dual$
 of all spherical representations
(this measure is called  the {\it Plancherel measure}).

As we have seen in 6.6,
the measure $\nu_\alpha$  is defined up to  a multiplication
by a positive function.
Now we shall define  a natural normalization of $\nu_\alpha$
and of
a unitary intertwining operator
$U_\alpha$ from $\T_\alpha$
to $\int \rho_sd\nu_\alpha(s)$.
 We require  {\it the image of the distinguished vector
$\Xi_\alpha\in \V_\alpha$ under     $U_\alpha$
 to be the function}
$s\mapsto \xi_s$\footnote{This normalization is
consistent with the normalization of the Plancherel measure
for $L^2$ defined in \Subsection 6.6, since the limit
of $\Xi_\alpha\in\V_\alpha$ as $\alpha\to+\infty$
is the $\delta$-function $\delta(z)$}.

We want to find the measure $\nu_\alpha$
normalized in this way. To do this,
 we calculate the matrix element
$\langle \T_\alpha (g)\Xi_\alpha,\Xi_\alpha\rangle$ in two ways.
We recall that this matrix element can be regarded
as $\K$-invariant  function on $\G/\K$.
A calculation in the kernel representation $\T_\alpha$ gives
$$ \langle \T_\alpha (g)\Xi_\alpha,\Xi_\alpha\rangle=
\det(1-zz^*)^\alpha=\prod_{k=1}^p(1+x_k)^{-\alpha}$$
A calculation in  the direct integral
  gives
\begin{equation}
\langle \T_\alpha (g)\Xi_\alpha,\Xi_\alpha\rangle=
\int\langle \rho_s(g) \xi_s,\xi_s\rangle_{s}d\nu_\alpha(s)
=\int \Phi_s(g) d\nu_\alpha(s)
\label{7.1}
\end{equation}
Thus we must find the measure $\nu_\alpha$
on $\dual$ such that
\begin{equation}
\int \Phi_s(x) d\nu_\alpha(s)=\prod_{k=1}^p(1+x_k)^{-\alpha}
\label{7.2}
\end{equation}

Conversely, if we have a measure $\nu_\alpha$ on $\dual$
satisfying
(\ref{7.2}), then  we have equality (\ref{7.1})
 for matrix elements.
Hence
 the direct integral is equivalent to $\T_\alpha$.%
\footnote{Let $\xi$ be a cyclic vector of a unitary representation
of a group $G$ in a Hilbert space $H$. Assume we know the matrix element
$\gamma(g)=\langle \rho(g)\xi,\xi\rangle$.
Let us explain why we know the representation $\rho$ itself.
Then
$$\gamma (g_2^{-1}g_1)= \langle \rho(g_2)^{-1}\rho(g_1)\xi,\xi\rangle
=\langle \rho(g_1)\xi,\rho(g_2) \xi\rangle$$
is a positive definite kernel (see 1.1) on $G$ and after this
we can reconstruct the Hilbert space $H$ with the distinguished
system of vectors $\rho(g)\xi$ in the usual way.}

\smallskip

 {\bf \punct. Plancherel formula for large values of $\alpha$.}

 \nopagebreak

 {\sc Theorem \fact.} (Berezin \cite{Ber2}) {\it Let
 $\alpha>q-p+1$. Then the Plancherel measure is
 supported by the pure imaginary $s$ and its density
 with respect to the Lebesgue measure
 on $i\R^p$ is}
\begin{multline}
   \frac {1}
   {\prod_{j=0}^{p-1}\Gamma^2(\alpha-j)}
   {\prod_{k=1}^p \Gamma(\alpha-h+s_k)
                      \Gamma(\alpha-h-s_k)}
       \times\\ \times
       \prod_{k=1}^p\frac{\Gamma^2 (r+s_k)\Gamma^2 (r-s_k)}
                 {\Gamma(2s_k)\Gamma(-2s_k)}
                 \prod_{1\le l <m\le p} (s_l^2-s_m^2)^2
 \label{7.3}
 \end{multline}

  In fact,  this is the product of (\ref{6.18}) and the
  Gindikin--Karpelevich density (\ref{6.10}).

  \smallskip

  {\sc Proof.}
   We want to represent $f(x):=\prod(1+x_j)^{-\alpha}$
  as the inverse spherical transform (\ref{6.15})
   of some function $F(s)$.
  It is
  sufficient to evaluate the direct spherical transform
  (\ref{6.14}) for $f(x)$. This was done in Theorem 6.13.
  This operation is correct  (see \cite{FK})
  if $f\in L^2\cap L^1$
  (i.e., $\alpha>q+p-1$).
   Then we consider analytic continuation.
  For $\alpha<h$ we have a pole of integrand
  and for $\alpha< h$ our Theorem 7.1
  becomes incorrect, see below.
  \kvadrat

  \smallskip

   {\bf \punct. Analytic continuation of the Plancherel formula.}
 In fact, the Plancherel formula obtained in Theorem 7.1
 is the following identity  for the hypergeometric functions
 \begin{multline}
 \prod_{k=1}^p(1+x_k)^{-\alpha}
 \prod_{1\le k<l\le p}(x_k-x_l)=
 \\=\const\cdot
  \frac {1}
   {\prod_{j=0}^{p-1}\Gamma^2(\alpha-j)}
\int\limits_{i\R^p}
\det_{k,m}\bigl\{ \F(r+s_m,r-s_m;2r;-x_k)\Bigr\}
\times \\ \times
   {\prod_{k=1}^p \Gamma(\alpha-h+s_k)
                      \Gamma(\alpha-h-s_k)}
       \prod_{k=1}^p\frac
{\Gamma^2 (r+s_k) \Gamma^2 (r-s_k)}
                 {\Gamma(2s_k)\Gamma(-2s_k)}
                  \prod_{1\le l <m\le p} (s_l^2-s_m^2)
 ds_1\dots ds_p
 \label{7.4}
 \end{multline}
 where $r=(q-p+1)/2,\,h=(q+p-1)/2$
(a direct verification of this formula is a nice exercise).

The left-hand side of the equation is holomorphic in $\alpha\in\C$.
Let us discuss the right-hand side.
The integrand has singularities if
\begin{equation}
\Re \alpha=h-k;\qquad k=0,1,2,\dots
\label{7.5}
\end{equation}
The $\Gamma$-function exponentially decreases in the
imaginary direction
(see, for instance, \cite{HTF}, v.1, (1.18.6))
$$|\Gamma(a+is)|=
 (2\pi)^{\frac 12} |s|^{\frac 12-a}e^{-\frac\pi2|s|}
 (1+o(1)),\qquad |s|\to\infty$$
and the spherical functions of unitary representations
are bounded by $1$. Hence the right-hand side
of (\ref{7.4}) is holomorphic
in $\alpha$
except for  the lines  (\ref{7.5}).

Let us construct the analytic continuation of the right-hand side
from the domain $\Re\alpha>\tfrac12(q+p-1)$ to
the whole  $\C$.

By Lemma 6.10, we can represent the identity
 (\ref{7.4}) in the form
\begin{multline}
 \prod_{k=1}^p(1+x_k)^{-\alpha}\prod_{1\le k<l\le p}
 (x_k-x_l)=\const\cdot
  \frac {1}
   {\prod_{j=0}^{p-1}\Gamma^2(\alpha-j)} \times \\
   \times
\det\limits_{k,m} \Bigl\{\int_{-i\infty}^{i\infty}
s^{2(m-1)}
\Gamma(\alpha-h+s)
                      \Gamma(\alpha-h-s)
\,\, \F(r+s,r-s;2r;-x_k)   \times \\
  \qquad\times
\frac {\Gamma^2 (r+s) \Gamma^2 (r-s)}
                 {\Gamma(2s)\Gamma(-2s)}
                 ds \Bigr\}
                 \label{7.6}
                 \end{multline}

{\sc  Lemma \fact.} {\it Denote by
$\int_{-i\infty}^{i\infty}I(\alpha,s)ds$ the integral
in the curly brackets in {\rm(\ref{7.6})}.
 Then the meromorphic continuation of
$\int_{-i\infty}^{i\infty} I(\alpha,s)$
 to the whole $\C$  is given by
\begin{equation}  \int\limits_{-i\infty}^{i\infty}
I(\alpha,s)\,ds+4\pi
\!\!\! \!\!\! \!
\sum\limits_{0\le k< h- \alpha}
\!\!\! \!\!\! \!
c_k(\alpha)
\F\Bigl[\begin{matrix}\alpha-p+1+k,-\alpha +q-k\\  2r
\end{matrix};-x\Bigr]
(\alpha-h+k)^{2(m-1)}
\label{7.7}
\end{equation}

where the coefficients $c_k(\alpha)$ are given by }
\begin{equation}
c_k(\alpha)= \frac{
 \Gamma(2\alpha-2h+k)\Gamma^2(-p+1+\alpha+k)
    \Gamma^2(q-\alpha-k) (-1)^k
                    }{
   \Gamma(2\alpha-2h+2k)\Gamma(-2\alpha+2h-2k) k!}
\label{7.8}
\end{equation}

{\sc Remark.}
 The expression (\ref{7.8}) has poles
at the points $\alpha=p-1-k,p-2-k,\dots$.
 Thus  expression (\ref{7.7}) has poles
at the points $\alpha=p-1,p-2,\dots$.

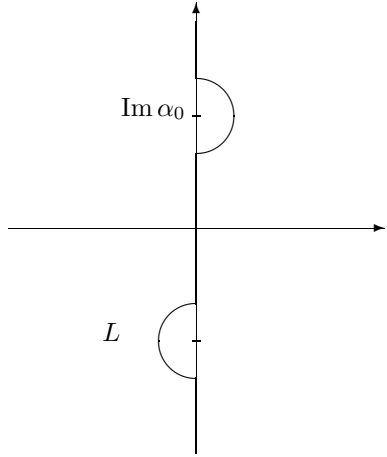
\hfill $\square$

\def\ris{
\unitlength=0.5mm
\begin{picture}(100,150)(-50,-80)
{\linethickness{0.005mm}
\put(0,-60){\vector(0,1){120}}
\put(-50,0){\vector(1,0){100}} }
\put(0,-20){\line(0,1){40}}
\put(0,40){\line(0,1){15}}
\put(0,-40){\line(0,-1){20}}
\put(0,30){\oval(20,20)[r]}
\put(0,-30){\oval(20,20)[l]}
\put(-1,30){\line(1,0){2}}
\put(-1,-30){\line(1,0){2}}
\put(-20,30){$\Im\alpha_0$}
\put(-25,-30){$L$}
\end{picture}}

\ris

 {\sc Proof.} We shall  obtain the first
 summand of the formula.
Let $\Re\alpha_0=h$, assume $\Im\alpha_0>0$.
 We want to construct an analytic continuation of the integral
$\int I(\alpha,s)ds$ to a small neighborhood of the point
$\alpha_0$.
Our integrand has poles at the points $s=\pm\Im\alpha_0$.
Consider the contour $L$   shown on the Picture.
In a small neighborhood of $\alpha_0$, the expression
$\int_L I(\alpha,s)ds$ depends on $\alpha$ holomorphically
and
$$
\int_{-i\infty}^{i\infty} I(\alpha,s)ds
- \int_L I(\alpha,s)ds
$$
is the sum of the residues. \kvadrat

\smallskip

Thus we obtain the analytic continuation
of  the right-hand side of (\ref{7.6}):
\begin{equation}
 \const\cdot
  \frac {1}
   {\prod_{j=0}^{p-1}\Gamma^2(\alpha-j)}
\det\limits_{k,m} \Bigl\{\int_\C
s^{2(m-1)}
 \F(r+s,r-s;2r;-x_k)
           d\mu_\alpha(s) \Bigr\}
\label{7.9}
                 \end{equation}
where the measure  $\mu_\alpha$ on the complex plane $\C$
 is the sum
of a continuous measure on the imaginary axis
and  $\delta$-measures
 supported by
the points $\pm(\alpha-h+k)$:
 \begin{multline*}
 d\mu_\alpha(s)=
  \Gamma(\alpha-h+s)
                      \Gamma(\alpha-h-s)
 \frac {\Gamma^2 (r+s) \Gamma^2 (r-s)}
                 {\Gamma(2s)\Gamma(-2s)}
\{ds\}+    \\ + 2\pi
\sum_{0\le k<h-\alpha}
c_k(\alpha)\delta\bigl(s\pm(\alpha+h-k)\bigr)
\end{multline*}
where $\{ds\}$ denotes the Lebesgue measure on
the imaginary axis and $c_k(\alpha)$ are the same as above
(\ref{7.8}).

\smallskip

{\sc Remark.}  The measures $\nu_\alpha$ are complex-valued
if $\alpha\in\C$.
They are real-valued for $\alpha\in\R$.
\hfill $\square$

\smallskip

Further, we apply Lemma 6.10 (from the right-hand side to the left-hand side)
and obtain the equality
\begin{equation}
\prod_{k=1}^p(1+x_k)^{-\alpha}=
\int_{\C^n}\Phi_s(x)d\nu_\alpha(s)
\label{7.10}
\end{equation}
where
\begin{equation}
d\nu_\alpha(s)=
\frac{1}
{\prod_{j=0}^p\Gamma^2(\alpha-j)}
\prod\limits_{1\le k<l\le p}
(s_k^2-s_l^2)^2\,\, d\mu_\alpha(s_1)\dots d\mu_\alpha(s_p)
\label{7.11}
\end{equation}
and the measures $\mu_\alpha$  are the same as above.

The identity  (7.10)
 is the required expansion of the distinguished matrix element
 in the spherical functions.

In fact, in this formula we have integration over
 the family
of planes having (up to the action
of the group $D_p$) the form
\begin{equation}
\Pi^\alpha_{k_1,\dots,k_\sigma}:\qquad
s_1=\alpha-h+k_1,\,\dots, s_\sigma=\alpha-h+k_\sigma,
\qquad s_{\sigma+1},\dots, s_p\in i\R
\label{7.12}
\end{equation}
where $\sigma=0,1,\dots, p$ and $k_j$ are nonnegative integers.
If some $k_j$ coincide, then
(due to the factor $\prod(s_k^2-s_l^2)^2$) the density
of measure on the plane (\ref{7.12}) is 0.
Hence only the case
$$k_1>k_2>\dots > k_\sigma\ge 0$$
 really exists.

By construction, our family of measures
is meromorphic in $\alpha\in\C$. Nevertheless
we have a factor $\prod\Gamma^{-2}(\alpha-j)$,
which is a zero at the poles of (7.8).
 It is easy to show that
our family of measures $\nu_\alpha$
 is holomorphic in $\alpha\in\C$.
The final formula without poles can be easily obtained
and hence we omit them.

\smallskip

{\bf\punct. Positive definiteness.}

\nopagebreak

\smallskip

{\sc Theorem \fact.} (\cite{Ner7}\footnote{a partial result
was obtained in \cite{Hil}})
{\it The formula {\rm (\ref{7.10})--(\ref{7.11})}
 is the Plancherel formula.}

 \smallskip

Identity (\ref{7.10}) has the form (\ref{7.2})
 but our considerations do not
imply that our  measure $\nu_\alpha$
  is supported by {\it unitary}
spherical representations.
 Hence some proof is necessary.
 An a priori proof
of positive definiteness of all representations
at
the support of the Plancherel  measure
is given in \cite{Ner7}.
For our case  $\G=\U(p,q)$, also it is possible
to apply Molev unitarizability results \cite{Mol}.

\smallskip

 {\sc Corollary \fact.} {\it The  Helgason transform
 is a unitary operator from the space of functions with the scalar
product (4.11) to
 $\int\rho_s d\nu_\alpha(s)$.}

\smallskip

 {\bf\punct. Discrete spectrum.}
  For the case $\sigma=p$,
 the plane  $\Pi^\alpha_{k_1,\dots,k_p}$  is a one-point
 set. Hence the representation
 $\rho_{\alpha-h+k_1,\dots,\alpha-h+k_p}$
 is a direct summand in $\T_\alpha$.
 These representations were subject of papers
 \cite{Ols1}, \cite{Mol}, \cite{NO}.

 All other planes  $\Pi^\alpha_{k_1,\dots,k_\sigma}$
 correspond to direct integrals of some spherical series
of unitary representations.

\smallskip

{\small
{\bf \punct. Comments.}
1)  Theorem 7.1  (large values of $\alpha$)
 was announced by Berezin \cite{Ber2}
for the series $G=\U(p,q)$, $\Sp(2n,\R)$ , $\SOS(2n)$, $\SO(n,2)$,
a proof was published by Unterberger and Upmeier
in \cite{UU}. For other series, the problem
was solved in \cite{Ner6}, the construction
 is based on the matrix
$\B$-function, which was constructed by Gindikin \cite{Gin1}
for the groups $\GL(n,\R)$, $\GL(n,\C)$, $\GL(n,\HH)$
and by the author for other series
of groups.
The general Plancherel formula (for arbitrary $\alpha$)
was obtained in \cite{Ner8}.

2) The simple  proof of the Plancherel formula given above
works also for series $\GL(n,\C)$, $\OO(n,\C)$, $\Sp(2n,\C)$.

3) There exists a counterpart of our analytic
continuation construction on the level of orthogonal
polynomials. It goes at least to Wilson \cite{Wil},
the most general construction
is contained in \cite{DS}.
}

\bigskip

{\large\bf \sect Boundary behavior of
holomorphic functions and separation of spectra.}

\medskip

Consider the planes $\Pi_{k_1,\dots,k_\sigma}^\sigma$
defined by (\ref{7.12}).
We have  a canonical decomposition
$$
\T_\alpha\simeq
\int \rho_s d\nu_\alpha(s)=
\bigoplus\limits_{\begin{matrix}\sigma; k_1,\dots,k_\sigma \\
               \sigma=0,\dots,p;\,\, h-\alpha>k_1>\dots> k_\sigma\ge 0
               \end{matrix}}
\int\limits_{s\in\Pi^\alpha_{k_1,\dots,k_\sigma}}
 \rho_s \,d\nu_\alpha (s)
 $$
 Our purpose is to obtain a natural realization
 of the summands of this decomposition.

 \smallskip

 {\bf\punct. Restriction of holomorphic functions to
 submanifolds in boundary.}
Let $\Omega\subset\C^N$ be an open domain,
and let $\partial \Omega$ be its boundary.
Suppose  $\Omega$ satisfies the conditions

 1) If $z\in \Omega$, $\lambda\in\C$, and $|\lambda|\le 1$,
 then $\lambda z\in\Omega$, i.e., $\Omega$ is a {\it circle domain}.

 2) If $z\in \partial\Omega$  and $|\lambda|< 1$,
 then $\lambda z\in\Omega$

 Let $K(z,u)$ be a positive definite kernel on $\Omega$.
Let $K(z,u)$ be holomorphic in $u$ and antiholomorphic in
$z$ (and hence the space $\H^\circ[K;\Omega]$
consists of holomorphic functions).

Assume that the kernel $K(z,u)$ is invariant
 with respect to the rotations
$$K(e^{i\phi}z,e^{-i\phi}u)=K(z,u)$$

{\sc Theorem \fact.} (\cite{NO})
 {\it Let $\mu$ be a positive measure
supported by a subset $M\subset\partial\Omega$.
Suppose  that

 {\rm 1)} The limit
  \begin{equation}
 K^*(z,u)= \lim_{\epsilon\to+0}
  K\bigl((1-\epsilon)z,(1-\epsilon)u\bigr)
  \label{8.1}
  \end{equation}
  exists almost everywhere on $M\times M$ with
  respect to $\mu\times\mu$.

  {\rm 2)} $K^*\in L^1(M\times M,\mu\times\mu)$
  and the limit {\rm(\ref{8.1})} is dominated, i.e.,
  there exists
  a positive function
  $\gamma(z,u)\in L^1(M\times M,\mu\times\mu)$
    such that
    $$|K\bigl((1-\epsilon)z,(1-\epsilon)u\bigr)|\le \gamma(z,u)$$

  Then

   {\rm a)} For any $f\in\H^\circ[K]$ the limit
  \begin{equation}
  R_\mu f(z):=\lim_{\epsilon\to+0} f((1-\epsilon)z)
  \label{8.2}
  \end{equation}
  exists almost everywhere   on $M$ with respect to $\mu$,
  and the restriction operator $R_\mu$ is a bounded operator
  $\H^\circ[K]\to L^1(M,\mu)$.

  {\rm b)} Let $\chi$ be a bounded
  measurable function on $M$.
   Then the limit
   $$l_{\chi}(f):=\lim_{\epsilon\to+0}
   \int_M f\bigl((1-\epsilon)z\bigr)\chi(z)d\mu(z)
   $$
   exists for all $f\in \H^\circ[L]$
    and $l_\chi$ is a bounded linear functional
   on $\H^\circ[L]$. Moreover, the map $\chi\mapsto\l_\chi$
   is a bounded operator from $L^\infty(M,\mu)$
   to the space of linear functionals on $\H^\circ[L]$.}

\smallskip

{\sc Remark.} Generally, functions $f\in\H^\circ[K]$ are discontinuous
on the boundary. Hence they have no values at an individual point
 $z\in\partial \Omega$. Theorem 8.1 claims that under some
conditions there
exists  the operator of restriction of function to a submanifold
in the boundary.

\smallskip

{\sc Remark.}
Statement b) of the Theorem means that the space $\H^\star[K]$
contains distributions supported by the subset
 $M\subset \partial \Omega$.
 \hfill $\square$

 \smallskip

In fact, the Theorem defines two
following  spaces of functions
 on
$M$.

\smallskip

The first  space
${\frak E}^\circ(M,\mu)$  consists
 of all functions on $M$ that can be obtained
by the restriction of $f\in \H^\circ[L]$.
This space  ${\frak E}^\circ(M,\mu)$   is a quotient space of $\H^\circ[L]$

\smallskip

The second space ${\frak E}^\star(M,\mu)$
 is the subspace  in $\H^\star[L]$
spanned by the (complex-valued) measures
$\chi(z)\mu(z)$, where $\chi(z)$ in $L^\infty(M,\mu)$.
 The space   ${\frak E}^\star(M,\mu)$
can be described more directly in the following way.
We consider the scalar product
$$
\langle \chi_1,\chi_2 \rangle:=
\iint_{M\times M}K^*(z,u)
\chi_1(u)\ov{\chi_2(z)}\,d\mu(z)\,d\mu(u)
$$
in $L^\infty(M,\mu)$. Then  ${\frak E}^\star(M,\mu)$
is a Hilbert space associated with the pre-Hilbert
space $L^\infty(M,\mu)$.

The space   ${\frak E}^\star(M,\mu)$
is a subspace of $\H^\star[K]$.

\smallskip

{\bf\punct. Restriction operators  in the spaces $V_\alpha$.}
We apply Theorem 8.1 to the objects described in \Subsection
4.2. The domain $\Omega$ is $\B_{p,q}\times \B_{p,q}$,
the kernel $K$ is the kernel $L_\alpha$ given by
(\ref{4.3})). The set  $M$ is a submanifold lying   in the boundary
of the diagonal $\Delta$: $z_1=\ov z_2$.

We obtain the following statement.

\smallskip

{\sc Corollary \fact.} {\it    Consider the space
$V_\alpha=\H^\circ[L_\alpha]$ of functions on $\B_{p,q}$ described in
\Subsection {\rm 4.1}. Let $M$ be a subset in the boundary
of $\B_{p,q}$,
and let $\mu$ be a measure supported by $M$.
Assume  $\mu$ satisfies   conditions {\rm 1--2}
of  Theorem {\rm 8.1}. Then the limit {\rm(\ref{8.2})} exists
and the restriction operator $R_\mu$ is a well-defined
operator $\H^\circ[L_\alpha]\to L^1(M,\mu)$.}

Now let us consider the $\G$-orbits $M_k$ in the boundary
of $\B_{p,q}$ (see \Subsection  2.6).

\smallskip

{\sc Theorem  \fact.} (\cite{Ner})
 {\it Let  $\alpha<(q-p+1)/2+k$.
Then the restriction operator is a well-defined operator
 from $V_\alpha= \H^\circ[L_\alpha]$ to
 the space $L^1_{loc}(M_k)$ of locally integrable functions
 on $M_k$.}

The question is
 reduced to an estimation of convergence of the integral
$$\iint_{A\times A\subset M_k\times M_k}
L_\alpha(z,u) dl(z)\,dl(u)$$
where $A$ is a compact subset in $M_k$ and $l(z)$
is the surface Lebesgue measure on $M_k$.
  I cannot simplify
estimates of \cite{Ner7} for our case $\G=\U(p,q)$
 and hence I omit a proof.

 Thus, we obtain the family of
 $\G$-invariant subspaces
 ${\frak E}^\star(M_k)\subset \H^\star[L_\alpha]$
 associated with the orbits
 $M_k$.

 \smallskip

 {\bf \punct. Restrictions of derivatives.}

 \nopagebreak

 {\sc  Theorem \fact.} (\cite{Ner6})
 {\it  Let  $\alpha<(q-p+1)/2+k-l$. Then the operator
 of restriction of partial derivatives of order $l$
 to $M_k$ is a well-defined operator}.

 {\sc Remark.} Discrete part of the spectrum
 corresponds to the compact $\G$-orbit $M_0$.
 \hfill $\square$

 For further discussion see \cite{NO}, \cite{Ner6}.

 \smallskip

 {\small
 {\bf\punct. Comments.} 1)The problem of separation of
 spectra in noncommutative harmonic analysis
 goes back to  Gelfand and Gindikin \cite{GG}.
  Olshanskii \cite{hardy} proposed  a way, which in
  some cases separates one of the pieces of spectra.
  Our way differs from \cite{GG}, \cite{hardy}, see also
  \cite{Mol2}.

   2. The condition $K^*(z,u)\in L^1$ is not necessary
for existence of
the restriction operator, some phenomena
related to the restriction problem are discussed
in \cite{Ner4}

   3) The discrete part of the spectrum of $\T_\alpha$
   corresponds to the minimal boundary orbit
   $M_0$. This part of the spectrum was the subject of the work
   \cite{NO}.
        }

 \bigskip

{\large\bf\sect Interpolation between
  $L^2\bigl(\U(p,q)/\U(p)\times\U(q)\bigr)$
  and
  $L^2\bigl(\U(p+q)/\U(p)\times\U(q)\bigr)$.
   Pickrell formula}

   \nopagebreak

  \medskip

{\bf\punct. Analytic continuation
of the Plancherel formula to negative integer $\alpha$.}
Assume that $\alpha$ in the Plancherel
formula (\ref{7.10})--(\ref{7.11})
 is a nonpositive integer, $\alpha=-N$.

 \smallskip

{\sc Theorem \fact.}
\begin{multline}
\prod_{j=1}^p(1+x_j)^N\prod\limits_{1\le k<l\le p}(x_k-x_l)=
 \prod_{j=1}^p\Gamma^2(N+j) \times  \\   \times
\sum\limits_{\begin{smallmatrix}m_1,\dots,m_p:\\
             N+p-1\ge m_1>m_2>\dots>m_p\ge0\end{smallmatrix}}
\Biggl\{\prod\limits_{j=1}^p   \frac
{(2m_j-p+q+1)\bigl((q-p+m_j)!\bigr)^2}
{\bigl(m_j!\bigr)^2\Gamma(N+q+m_j+1)\Gamma(N+p-m_j)}
\times  \\  \times
\prod\limits_{1\le k<l\le p}(m_k-m_l)(m_k+m_l+q-p+1)
\Biggr\}
\cdot
\det\limits_{j,l} \Bigl\{
\F(-m_j,q-p+1+m_j;q-p+1;-x_l)\Bigr\}
\label{9.1}
\end{multline}
{\it or}
\begin{multline}
\prod_{j=1}^p(1+x_j)^N=
 \prod_{j=1}^p\Gamma^2(N+j) \times  \\   \times
\sum\limits_{\begin{smallmatrix}m_1,\dots,m_p:\\
             N+p-1\ge m_1>m_2>\dots>m_p\ge0\end{smallmatrix}}
\Biggl\{\prod\limits_{j=1}^p   \frac
{(2m_j-p+q+1)\bigl((q-p+m_j)!\bigr)^2}
{\bigl(m_j!\bigr)^2\Gamma(N+q+m_j+1)\Gamma(N+p-m_j)}
\times  \\  \times
\prod\limits_{1\le k<l\le p}(m_k-m_l)^2(m_k+m_l+q-p+1)^2
\Biggr\}        \cdot  \,
\Phi_{p-h-m_1-1,\dots, p-h-m_p-1}(x_1,\dots,x_p)
\label{9.2}
\end{multline}

{\sc Remark.}  The hypergeometric functions in the right-hand
side of (\ref{9.1}) are the {\it Jacobi polynomials}
(see \cite{HTF}, v.2, or \cite{AAR}, Chapters 2,6)
$$
\F(-m,q-p+1+m_j;q-p+1;-x_l)=
\frac{m!\Gamma(q-p+1)}{\Gamma(q-p+m+1)}
P^{q-p,0}_m(1-2x)
$$
Thus the right-hand side of (9.1) can be represented in the form
\begin{multline*}
\Gamma^p(q-p+1) \prod_{j=1}^p\Gamma^2(N+j) \times  \\   \times
\sum\limits_{\begin{smallmatrix}m_1,\dots,m_p:\\
             N+p-1\ge m_1>m_2>\dots>m_p\ge0\end{smallmatrix}}
\Biggl\{\prod\limits_{j=1}^p   \frac
{(2m_j-p+q+1)(q-p+m_j)!}
{m_j!\Gamma(N+q+m_j+1)\Gamma(N+p-m_j)}
\times  \\  \times
\prod\limits_{1\le k<l\le p}(m_k-m_l)(m_k+m_l+q-p+1)
\Biggr\}
\cdot
\det\limits_{j,l} \Bigl\{
P^{q-p,0}_{m_j}(1-2x)\Bigr\}
\end{multline*}

{\sc Proof of Theorem 9.1.}
The factor
$\prod_{j=0}^p\Gamma^{-2}(\alpha-j)$ of (\ref{7.11})
has zero of order $2p$ at $\alpha=-N$.
   Hence a summand of the Plancherel
formula can be nonzero only in the case than the product
$\prod_j c_{k_j}(\alpha)$ has a pole of order $2p$
at $\alpha=-N$.
But only the factors
\begin{equation}
\Gamma^2(-p+1+\alpha+k_j)
\label{9.3}
\end{equation}
of $c_{k_j}(\alpha)$ (see (\ref{7.8})) can give poles
(the poles of the first factor in the numerator in (\ref{7.8})
are annihilated by the poles of the first factor of the denominator).
Hence all $p$ factors (\ref{9.3}) have to be present, and hence
the continuous components of the Plancherel measure $\nu_\alpha$
are absent.
Now it remains to write the formula.\kvadrat

\smallskip

Our next aim is to explain the group-theoretical meaning
of Theorem 9.1.

\smallskip

{\bf\punct. Symmetric spaces $\U(p+q)/\U(p)\times\U(q)$.}
Consider the space $\C^p\oplus\C^q$ equipped with
the standard scalar product. We denote by $\U(p+q)$
the unitary group of this space. We write elements of
$\U(p+q)$ as $(p+q)\times(p+q)$ matrices
$g=\bigl(\begin{smallmatrix}
a&b\\c&d
\end{smallmatrix}
\bigr)$.

We denote by $\Gr_{p,q}$ the space
of all $p$-dimensional subspaces
in $\C^p\oplus\C^q$. Obviously, $\Gr_{p,q}$
is a $\U(p+q)$ homogeneous space
$$\U(p+q)/\U(p)\times\U(q)$$

The graph of a linear operator   $\C^p\to\C^q$
is an element of $\Gr_{p,q}$, and
 elements of $\Gr_{p,q}$ in a general position
  have this form. Hence we obtain
 a parametrization of an open dense subset in $\Gr_{p,q}$ by
 points of the space $\Mat_{p,q}$ of all
 $p\times q$ matrices.

 The action of the group $\U(p+q)$ on
 the space $\Mat_{p,q}$ is given by
 the same linear-fractional formula as above (\ref{2.7}).

 The $\U(p+q)$-invariant measure on the space
 $\Mat_{p,q}$  is given by
 \begin{equation}
 \det(1+zz^*)^{-p-q} dz
 \label{9.4}
 \end{equation}

 \smallskip

 {\bf\punct. Representations $\tau_{-N}$ of the group
 $\U(p+q)$.} We define these representations in three ways.

 1) $\tau_{-N}$ is the irreducible representation
 of $\U(p+q)$ with
 the signature
 $(N,\dots,N,0\dots,0)$, where $N$ is on the first
 $p$ places (see \cite{Zhe}, \S49).

 2)  Denote by $Det$ the determinant line bundle on
 $\Gr_{p,q}$. Denote by $Det^{\otimes N}$ the $N$-the
 power of $Det$. Then
 $\tau_{-N}$ is the representation of $\U(p+q)$ in
 the space of holomorphic sections of
 $Det^{\otimes N}$.

 3) Consider the space $\H^\circ[L_{-N}]$
  of functions on $\Mat_{p,q}$
 defined by the positive definite kernel
 $$L_{-N}(z,u)=\det(1+uz^*)^N$$

{\sc Lemma \fact.} {\it The kernel $L_{-N}$ is positive definite.}

\smallskip

{\sc Proof.} By Proposition 1.6, it is sufficient to prove the
statement for $N=1$. Consider Euclidean space
$\C^{p+q}$ with an orthonormal basis
$e_1,\dots, e_p, h_1,\dots,h_q$. Consider the $p$-the
exterior power $\Lambda^p\C^{p+q}$ and the system of vectors
$$v_z:=(e_1+\sum \overline z_{1j}h_j)\wedge
       (e_2+\sum \overline z_{2j}h_j)\wedge \dots
  \wedge
     (e_p+\sum \overline z_{pj}h_j)
\,\,\,\in \Lambda^p\C^{p+q}
$$
Then $\langle v_z,v_u)=\det(1+zu^*)$.  \hfill $\boxtimes$

\smallskip

 The functions
 $$\theta_z(u):=L_{-N}(z,u)$$ are  holomorphic
 polynomials  of degree $\le pN$. By  Corollary 1.4,
 all elements of the space  $\H^\circ[L_{-N}]$
 are holomorphic polynomials of degree $\le pN$.

 The group $\U(p+q)$ acts in   $\H^\circ[L_{-N}]$
 by the operators
 $${ \tau}_{-N}(g)f(z)=f\bigl((a+zc)^{-1}(b+zd)\bigr)
   \det(a+zc)^N$$

{\sc Remark.} The image of an element
 $\theta_u(z):=\det(1+zu^*)^N$ of the supercomplete system
under the transformation ${ \tau}_{-N}(g)$
 is an element of the supercomplete system (up to a factor).
 Hence the operators $\tau_{-N}(g)$
 preserve the space $\H^\circ(L_{-N})$.
 \hfill$\square$

 \smallskip

   {\bf \punct. Kernel representations of $\U(p+q)$.}
 Consider the positive definite kernel
 $$
 {\cal L}_{-N}(z,u)=\frac {|\det(1+uz^*)|^{2N} }
    {\det(1+zz^*)^N\det(1+uu^*)^N}
 $$
and the space   $\V_{-N}:=\H^\circ[{\cal L}_{-N}]$
defined by this kernel.
We define the {\it kernel representation} $\T_{-N}$ of the group
$\U(p+q)$ in the space  $\H^\circ[L_{-N}]$
by the formula
$${\cal T}_{-N}(g)f(z)=f\bigl((a+zc)^{-1}(b+zd)\bigr)$$

We also define the {\it distinguished vector}
$\Xi_{-N}$. It is the function
$f(z)=\det(1+zz^*)^{-N}$.
The orbit of the vector $\Xi_{-N}$ consists of all elements
of the supercomplete system $\L_{-N}(z,a)$ and hence the vector
$\Xi_{-N}$ is cyclic.

The same arguments as in \Subsection 4.3 show that
the representation $\T_{-N}$ can be decomposed into
 a tensor product
of  $\tau_{-N}$ and the complex conjugate representation
$\ov \tau_{-N}$. Also the space $\H^\circ[{\cal L}_{-N}]$
can be canonically identified with the space of operators
 $\H^\circ[L_{-N}] \to \H^\circ[L_{-N}]$.

 \smallskip

 {\bf\punct. Why formulas in Sections
  4 and 9 are similar?}
 Representations $\tau_\alpha$
(see 3.1) make sense for arbitrary complex
 $\alpha$. For arbitrary real $\alpha$ the kernel
 $K_\alpha(z,u)=\det(1-uz^*)^{-\alpha}$
  defines some scalar product
 in a space of holomorphic functions,
 but this scalar product is not positive definite.
 The work in such spaces is possible but it is
 difficult from analytical point of view.

 For a negative integer $\alpha$ our representations are finite-dimensional
 and hence they can be  extended holomorphically to
 $\GL(p+q, \C)$.
  It is more natural to consider them as representations
 of  the compact form $\U(p+q)$ of the group $\GL(p+q, \C)$.

 In fact, the formulas of Sections 4 and 9
for actions of groups and scalar products really coincide
 and the reason of small difference
 in signs is the jump to another
 real form of the group $\GL(p+q, \C)$.

 \smallskip

 {\bf\punct. Plancherel formula.}
 We say that an irreducible representation
 of $\U(p+q)$ is $\U(p)\times \U(q)$-{\it spherical} if
 it contains a nonzero $\U(p)\times \U(q)$-fixed vector.
 By the same Gelfand Theorem 6.4, this vector is unique.

 The $\U(p)\times\U(q)$-spherical functions
 of $\U(p+q)$ are given
 by the same Berezin--Karpelevich formula%
 \footnote{In this case, $-1\le x_j\le 0$ are the eigenvalues
 of the matrix $\bigl(-z^*(1+zz^*)z\bigr)$,
 where $z\in \Mat_{p,q}$.}.

 \smallskip

 {\sc Lemma \fact.} {\it The kernel representation
 ${\cal T}_{-N}$ is a multiplicity free sum of
the spherical
 representations.}

 \smallskip

 {\sc Proof.} This is almost a special case of Theorem 6.3b.
 Nevertheless we shall give an independent proof.

 First the vector $\Xi_{-N}$ is cyclic and
 $\U(p)\times\U(q)$-invariant. Hence
 the orthogonal projection
 of $\Xi_{-N}$ to any subrepresentation is
$\U(p)\times\U(q)$-invariant and cyclic
 in the subrepresentation. Thus all irreducible
 subrepresentations of $\T_{-N}$ are spherical.

  Secondly, assume that $\T_{-N}$ contains two copies
  of the same irreducible
   representation $\rho^j$ of $\U(p+q)$.
 Let $V$, $V'$ be the corresponding
 orthogonal subspaces
  and let $\xi\in V$, $\xi'\in V'$
   be the spherical vectors. Let $J:V\to V'$
   be the unique intertwining  linear operator
   $V\to V'$ such that $J\xi=\xi'$.
   The projection of $\Xi_{-N}$ to $V\oplus V'$
   must have the form
   $\lambda \xi\oplus \mu \xi'$, where $\lambda,\mu\in\C$.
   But this vector is not cyclic in $V\oplus V'$
   (since it is contained in the graph of the operator
   $\mu\lambda^{-1}J:V\to V'$).
    \kvadrat

\smallskip

 The definition of the Plancherel measure given in
 \Subsection  7.1 is valid in our situation.
 But  direct integrals
 here reduce to finite sums and we shall give the definition again.

 Let $\rho^{j}$ be the spherical representations that occur
 in the spectrum of $\T_{-N}$, let $\Phi^j$ be the spherical
 function of $\rho^{j}$,
and
let $\xi^j$ be the projection of
the distinguished vector $\Xi_{-N}$ to the subspace $V^j$.
The Plancherel measure in our case is the collection of
numbers
$$\nu^j=\langle \xi^j,\xi^j\rangle$$
The matrix element of the distinguished vector
can be expanded into the sum of spherical functions:
$$
\langle {\cal T}_{-N}(g)\Xi_{-N},\Xi_{-N}\rangle=
\langle \sum  \rho^j(g) \xi^j,\sum\xi^j\rangle=
\sum \langle \rho^j(g) \xi^j,\xi^j\rangle=
\sum\langle \xi^j,\xi^j\rangle\Phi^j(g)$$
 As a result, we obtain

 \smallskip

 {\sc Theorem \fact.} {\it The Plancherel coefficients
 $\nu^j$  are the coefficients in formula {\rm (\ref{9.2})}.}

 \smallskip

 {\bf\punct. Limit as $N\to\infty$.} Consider the space
 $\H^\star[{\cal L}_{-N}]$ (see \Subsection 1.6).
 Consider
  measures
 on $\Mat_{p,q}$ of the form
 $\phi(z)\det(1+zz^*)^{-p-q}\{dz\}$, where
 $\phi$ are  compactly supported  smooth functions.
 Then we obtain a scalar product
 in the space of smooth compactly
 supported functions given by
 \begin{equation}
 \langle \phi,\psi\rangle_{-N}=
 C(-N)^{-1}
 \!\!\! \!\!\!\!\!   \!\!\!  \!
 \iint\limits_{\Mat_{p,q}\times\Mat_{p,q} }
  \!\!\! \!\!\!\!\!   \!\!\!
 {\cal L}_{-N}(z,u)
 \phi(z)\overline{\psi(u)}
 \frac{ \{du\}\,\{dz\} } { \det(1+zz^*)^{p+q}
 \det(1+uu^*)^{p+q}}
 \label{9.5}
 \end{equation}
 It is natural to choose the normalization constant by
 $$C(-N)=\int_{\Mat_{p,q}}(1+zz^*)^{-N-p-q}\{dz\}$$

 {\sc Proposition \fact.} {\it The limit
 of scalar products {\rm(\ref{9.5})}  as $N\to+\infty$
 is the $L^2$-scalar product with respect to
 the $\U(p+q)$-invariant measure {\rm(\ref{9.4})}.}

 {\sc Proof} is similar to the proof of  Theorem 4.3.
 \kvadrat

 \smallskip

 {\bf \punct. Pickrell formula.}

 \nopagebreak

 {\sc Theorem \fact.} {\it Let us omit the condition
 $N+p-1\ge  m_1$ in the summation in {\rm(\ref{9.2})}.
 Then  formula {\rm(\ref{9.2})} remains
 valid for all complex numbers $N$ such that
 $\Re N>-1$.}

 \smallskip

 For $p=q$ this gives the Pickrell formula \cite{Pic}
 (in this case, the
  Jacobi polynomials are the Legendre polynomials).

 \smallskip

 {\sc Remark.}  Let $N$ be a positive integer.
 Then omitting of the condition  $N+p-1\ge  m_1$
 does not change  formula (\ref{9.2}). All new summands,
 which appear in the formula,
 are zeros, since we have the factor $\Gamma(N+p-m_1)$
 in the denominator.     \hfill $\square$

 \smallskip

 {\sc Theorem \fact.} (Carlson, see \cite{AAR}, 2.8.1)
 {\it   Let $f(z)$ be holomorphic function for $\Re z>0$, let
 $f(z)=O(e^{(\pi-\epsilon)|z|})$,
 and $f(n)=0$ for all integers.
 Then $f(z)$ is identically zero.}

 \smallskip

 {\sc Proof of Theorem 9.6.}
 A simple calculation with Lemma 6.10 shows
 that the functions
  $\Phi_{p-h-m_1-1,\dots, p-h-m_p-1}(x_1,\dots,x_p)$
  constitute an orthogonal system in the
 space $L^2$ of symmetric functions on the cube $-1\le x_j \le 0$
 with respect to the measure
 $$d\sigma(x):=
 \prod_k (-x_k)^{q-p}\prod(x_k-x_l)^2 dx_1\dots dx_p$$
 Hence the coefficients  of the expansion
 of $\prod(1+x_j)^{N}$  in
 the spherical functions  $\Phi_{\dots}$
 are $L^2$-scalar products
 $$
 \widetilde c_m(N)=\const(m)\int_{[-1,0]^p} \prod(1+x_k)^{N}
 \Phi_{p-h-m_1-1,\dots, p-h-m_p-1}(x)\,d\sigma(x)
 $$
 Obviously,  these coefficients are bounded for
 $\Re N\le 0$.
Since
$$\Gamma(N+a)/\Gamma(N+b)\sim N^{a-b},\qquad N\to+\infty$$
(see \cite{HTF}, v.1 (1.19.4)), it follows that
 the coefficients $c_m(N)$ of (\ref{9.2})
  also are bounded.
 Now we apply the Carlson theorem to $\widetilde c_m(N)-c_m(N)$.
 \kvadrat

 \smallskip

 {\small
 {\bf \punct. Comments.}
 1) The formula (\ref{9.2}) itself  is an extension of
 the Pickrell
formula \cite{Pic}. Nevertheless, our method
of analytic continuation (see \cite{Ner7}) is
the same for all symmetric spaces and the
decompositions of the type (\ref{9.2}) follow automatically
from the $\B$-integrals evaluated in \cite{Ner5}.

\smallskip

 2) The limit of  formula (\ref{9.2}) as $N\to+\infty$ gives
 the Plancherel formula for
$L^2$ on the compact symmetric space
 $L^2\bigl(\U(p+q)/\U(p)\times\U(q)\bigr)$.

 \smallskip

 3) In particular, this gives the Helgason theorem
 on the description of spherical representations of
 compact groups, see \cite{Hel0},
 \cite{Hel2}, Theorem 5.4.1.

 \smallskip

  4) As in Subsection 4.6, we obtain a canonical basis
  in the subspace$\V_{-N}^\K$  of $\K$-invariant functions
in $\V_{-N}$
(since $\V_{-N}^\K$ is equivalent to the space of operators
in $\H^\circ[L_{-N}$). Its image under the spherical transform
  consists of determinants of dual Hahn polynomials
  (see \cite{KS} on the dual Hahn polynomials).

   For other series this gives multivariate dual Hahn polynomials
   or (for series $\GL_n$) the
    multivariate Krawtchouk polynomials.

\smallskip

  5) There exists a natural {\it inverse} (!)
   limit of the symmetric spaces
   $$
   \lim\limits_{\longleftarrow} \U(p+k+p)/\U(p)\times\U(k+p)
   \qquad\mbox{as} \,\,\, p\to+\infty
   $$
 It was constructed in the important work
 of  Pickrell \cite{Pic} (for $k=0$).
  Pickrell's type construction  exists for all classical
 compact symmetric spaces (\cite{non}, \cite{Ner9}, \cite{BO}).
 The are many similarities  between  harmonic analysis
 on the inverse limits of symmetric spaces and the analysis
 of Berezin kernels (see \cite{Pic}, \cite{Ner7}, \cite{Ner8},
 \cite{BO}).
        }

 \bigskip

{\large\bf\sect  Radial part of the spaces
$\V_\alpha$ and Gross--Richards kernels}

\nopagebreak

\medskip

{\bf\punct. An orbital  integral.}
Consider the space $\V_\alpha^\K$ consisting of
$\K$-invariant functions $f\in\V_\alpha$, see \Subsection  4.6.

\smallskip

{\sc Proposition \fact.}
{\it The space $\V_\alpha^\K$ has the form
$\H^\circ[R_\alpha]$, where the reproducing kernel
$R_\alpha(z,u)$, $z,u\in\B_{p,q}$, is given by
 \begin{equation}
 R_\alpha(z,u)=
   \det(1-zz^*)^\alpha\det(1-uu^*)^\alpha
 \!\!\!  \!\!\!    \!\!\!  \!\!\!
 \int\limits_{ {\bold k}_1\in \U(p),{\bold k}_2\in\U(q)}
\!\!\! \!\!\!    \!\!\!  \!\!\!
 |\det(1-{\bold k}_1z {\bold k}_2^{-1}u^*)
    |^{-2\alpha}d{\bold k}_1\,d{\bold k}_2
\label{10.1}
\end{equation}
  where $d{\bold k}_1$, $d{\bold k}_2$
   denote the Haar measures
  on the unitary groups $\U(p)$, $\U(q)$
 such that  the measure of the whole group is 1.}

 \smallskip

 {\sc Remark.} The kernel $R_\alpha(z,u)$ is invariant
 with respect to the transformations
 $z\mapsto {\bold k}_1z{\bold k}_2^{-1}$,
  $u\mapsto {\bold l}_1u{\bold l}_2^{-1}$,
  where $ {\bold k}_1, {\bold l}_1\in\U(p)$,
           ${\bold k}_2, {\bold l}_2\in\U(q)$.
 Hence all elements of the space $\H^\circ[R_\alpha]$
 are $\U(p)\times\U(q)$-invariant functions.
 \hfill $\square$

 {\sc Proof.}
Let
 $a\in \B_{p,q}$.
We want to find a function
 $\kappa_a\in\V_\alpha^\K$ such that
for any $f\in \V_\alpha^\K$,
\begin{equation}
f(a)=\langle f,\kappa_a\rangle_{\V_\alpha^\K}
\label{10.2}
\end{equation}

Let
$$
\theta_a(z)={\cal L }_\alpha(z,a)=
\frac{\det(1-aa^*)^\alpha\det(1-zz^*)^\alpha}
{|\det(1-za^*)|^{2\alpha}}
$$
Then, by the reproducing property (\ref{1.3})
\begin{equation}
 f(a)=\langle f,\theta_a\rangle_{\V_\alpha}
\label{10.3}
\end{equation}
This  implies boundedness
 of the linear functional $f\mapsto f(a)$.
  Thus we can apply \Subsection 1.3.
 Thus the reproducing kernel exists.

 The function $f$ is $\K$-invariant and hence
 we   can replace $\theta_a$ in  equation (\ref{10.3})
 by its average over the group $\K$.
 This gives the integral expression  (10.1)
 for the function $\kappa_a$
 and for the reproducing kernel $R_\alpha(z,u)$ of our space.
 \kvadrat

 \smallskip

 {\bf \punct. Evaluation of the reproducing kernel.}
 Let us write functions $f\in\V_\alpha^\K$ as
functions depending of the variables $x_1,\dots,x_p$,
 see (\ref{2.9}).

 \smallskip

 {\sc Theorem \fact.}
 \begin{multline}
 R_\alpha(x,y)=\frac{\Gamma^{2p}(\alpha-p)\prod_{j=0}^{q-1}j!}
             {(q-p)!^p\prod_{j=1}^p\Gamma^2(\alpha-j+1)}
  \prod_{k=1}^p(1+x_k)^{-\alpha+p-1}(1+y_k)^{-\alpha+p-1}
             \times \\ \times
           \frac{  \det\limits_{k,l}\Bigl\{
             \F\Bigl[\begin{matrix} \alpha-p,\alpha-p\\
               q-p+1\end{matrix}
               ;\frac{x_k y_l}{(1+x_k)(1+y_l)}\Bigr]\Bigr\}}
               {\prod_{1\le k<l\le p}(x_k-x_l)
               \prod_{1\le k<l\le p}(y_k-y_l)}
\label{10.4}
 \end{multline}

 {\sc Remark.} The kernel (\ref{10.4}) is a special case
 of the Gross-Richards kernels \cite{GR}.
 \hfill $\square$

 \smallskip

 {\sc Theorem \fact.} (S. Bergman, 1947)
 {\it  Let $K$ be a positive definite kernel on a set $X$.
 Let  $\zeta_1(x),\zeta_2(x),\dots$ be an orthonormal basis
 in the space $\H^\circ[K]$.  Then
 $$K(x,y)=\sum_j \ov{\zeta_j(x)}\zeta_j(y)$$
 and the series in the right-hand
 side converges on $X\times X$. }

 \smallskip

 {\sc Proof of Theorem 10.3.}
 Consider the expansion of the function
 $\theta_a(x):=K(a,x)$ with respect to the basis $\zeta_j$
 \begin{equation}
 \theta_a(x)=\sum c_j(a) \zeta_j(x)
 \label{10.5}
 \end{equation}
 This series converges in $\H^\circ[K]$
 and hence it converges
 pointwise.       Let us evaluate
  $\langle \theta_a(x),\zeta_k(x)\rangle$
 in two ways. By (\ref{10.5}), it is equal  $c_k(a)$.
 By the reproducing property (\ref{1.3}),
  it equals to  $\ov{\zeta_k(a)}$.
 \kvadrat

 \smallskip

{\sc Proof of Theorem 10.2.} We must evaluate
\begin{equation}
\sum\limits_{\mu_1\ge \mu_2\ge \dots\ge \mu_p}
\frac{ \Delta_{\mu_1,\dots,\mu_p}(x)
\Delta_{\mu_1,\dots,\mu_p}(y)}
 {     \langle  \Delta_{\mu_1,\dots,\mu_p},
      \Delta_{\mu_1,\dots,\mu_p}\rangle}
 \label{10.6}
 \end{equation}
where $\Delta_\mu$ is the canonical orthogonal basis
in $\V_\alpha^\K$
 defined by  (\ref{4.12}).
  This  reduces to the evaluation of the series
 $$\sum\limits_{\mu_1\ge  \dots\ge \mu_p\ge0}
 \prod_{j=1}^p\frac{\Gamma^2(\alpha+\mu_j-j+1)}
             {(\mu_j+p-j)!(\mu_j+q-j)!}
              s_{\mu_1, \dots, \mu_p}(X_1,\dots,X_p)
               s_{\mu_1, \dots, \mu_p}(Y_1,\dots,Y_p)
              $$
   where
$$
   X_k=x_k/(1+x_k),\qquad Y_k=y_k/(1+y_k)
$$
    and $s_\mu$ are
   the Schur functions.

   \smallskip

 {\sc Theorem \fact.} (Hua, \cite{Hua})
 $$\det\limits_{k,l}\bigl\{ \sum_{n\ge 0}
  a_n X^n_kY^n_k\bigr\}
           =
  \sum\limits_{n_1>n_2>\dots>n_p\ge 0} a_{n_1}\dots a_{n_p}
          \det\limits_{k,j}\{X_k^{n_j}\}
             \det\limits_{k,j}\{Y_k^{n_j}\}
  $$

  This can be verified directly.
  Let us write this identity
  in the form
  $$\frac{\det\limits_{k,l}\bigl\{ \sum_{n\ge 0}
  a_n X^n_kY^n_k\bigr\}}
{\prod\limits_{1\le k<l\le p}(X_k-X_l)
\prod\limits_{1\le k<l\le p}(Y_k-Y_l)}=
  \sum\limits_{\mu_1\ge \dots \ge \mu_p\ge 0}
  a_{\mu_1+p-1}a_{\mu_2+p-2}\dots a_{\mu_p}\,\, s_\mu(X) s_\mu(Y)
  $$

  Now Theorem 10.2 becomes obvious.\kvadrat

  \smallskip

  {\bf \punct. Spherical transform in the spaces
$\V_\alpha^\K$.}

\nopagebreak

\smallskip

{\sc Proposition \fact.} {\it
Let $\alpha>h=\frac12(q+p-1)$.
Then the spherical  transform is a unitary operator
from  $\V_\alpha^\K$ to the space $L^2$ with respect
to the Plancherel measure} (\ref{7.1}).

\smallskip

This is a rephrasing of Theorem 7.1.
This is also a simple corollary of the exotic
Plancherel formula (see 5.3) for
the index hypergeometric transform.

\smallskip

{\small
{\bf\punct. Comments.}
1) Analogues of the orbital integral (\ref{10.1})
(the average of the Berezin kernel over the compact
subgroup)
exist for all matrix balls.

\smallskip

2) The calculation of \Subsection
10.2  is valid
also for the series $\Sp(2n,\R)$, $\SOS(2n)$, $\GL(n,\C)$.
The last case was considered by Gross and Richards
\cite{GR}.

\smallskip

  3) For other matrix ball series, the Schur functions in (\ref{10.6})
  are replaced by Jack polynomials. I cannot
  identify these kernels with some known special
  functions.  It seems, that they define "new" scalar products
  in the space of symmetric functions.}

\bigskip

{\large\bf \sect \O rsted problem. Identification
of kernel representations and $L^2(\G/\K)$}

\nopagebreak

\medskip

{\bf \punct. The problem of unitary equivalence.}
As we have seen, for $\alpha>(q+p-1)/2$
the Plancherel measure $\nu_\alpha$ for $\T_\alpha$ differs
 from the Plancherel measure  $\nu_\infty$
for $L^2(\G/\K)$ by
a functional nonvanishing factor
(\ref{6.18}). Hence
for $\alpha>(q+p-1)/2$ the representation $\T_\alpha$ is
equivalent to the  representation of $\G$ in
$L^2(\G/\K)$.

It is easy to construct an intertwining operator
from $\T_\alpha$ to $L^2(\G/\K)$. The simplest possibility
is to consider the identity map $f\mapsto f$.\footnote{%
Scalar products in  $\T_\alpha$ and $L^2(\G/\K)$
are different and problem of  boundedness
of the identity map is nontrivial.
This identical operator $Id$ was discussed in many papers
(see \cite{Ber2},
\cite{Gut}, \cite{Rep}, \cite{OO}).
In \cite{OO} there was announced boundedness of $Id$
 for a large $\alpha$.
 Clearly, boundedness follows
from the explicit Plancherel formula (the expression (\ref{6.18})
is bounded for a fixed $\alpha$;
 this gives also an explicit formula
for the norm of $Id$.
On the other hand, the Plancherel formula for  vector-valued kernel
 representations is not known, and hence
a priori proofs of boundedness of $Id$
preserve sense.}

There arises a problem of construction of an
 explicit unitary operator
 \begin{equation}
 J_\alpha:  L^2(\G/\K) \to  \V_\alpha
 \label{11.1}
 \end{equation}

 {\bf \punct. $\Lambda$-function.}
 Our construction is based on one special function.
 Let $a,b,c>0$.
 Following \cite{Ner8}, we define
the $\Lambda$-function  by
 \begin{multline*}
 \Lambda^a_{b,c}(x)=\frac {1}{\Gamma(b+c)}
 \int\limits_0^\infty \Gamma(a+is)
 \frac {\Gamma(b+is)\Gamma(b-is)\Gamma(c+is)\Gamma(c-is)}
 {\Gamma(2is)\Gamma(-2is)}
    \times\\ \times
    \F(b+is,b-is;b+c;-x)\,ds
 \end{multline*}
 It seems that $\Lambda$-function cannot be expressed
in terms of the standard special functions
 by algebraic operations
 (except some special values of the parameters $b,c$).
 I discussed
  this function in detail in \cite{Ner8}.

  \smallskip

 {\bf\punct. Construction
 of the unitary intertwining operator.}  We preserve the notation
$\Xi_\alpha$ for the distinguished vector in   $\V_\alpha$
(see 4.4) and the notation ${\cal L}_\alpha(z,u)$
for the reproducing kernel of the space     $\V_\alpha$.

 Suppose we know the function
 ${\frak F}=J_\alpha^{-1}\Xi_\alpha\in L^2(\G/\K)$.
 For
 $g=
 \bigl(\begin{smallmatrix} a&b\\c&d\end{smallmatrix}\bigr)$,
  the image of the function
 $$\Xi_\alpha(z^{[g]})={\cal L}_\alpha(z,a^{-1}b)$$
 under $J_\alpha^{-1}$ is  ${\frak F}(z^{[g]})$.
 Hence if we know $J_\alpha^{-1}\Xi_\alpha$,
 then we know the $J_\alpha^{-1}$-image of
 all supercomplete system $\theta_u(z)={\cal L}_\alpha(z,u)$.
 We define the unitary $\G$-intertwining operator
 $J_\alpha:L^2(\G/\K)\to \V_\alpha$ by the formula
 \begin{equation}
 J_\alpha f(g)=\langle f, {\frak F}(z^{[g]}) \rangle_{L^2}=
 \int\limits_{\B_{p,q}} f(z) \ov{{\frak F}(z^{[g]})}
  \det(1-zz^*)^{-p-q} \{dz\}
 \label{11.2}
 \end{equation}
The distinguished vector $\Xi_\alpha$
 is $\K$-invariant, hence
the function ${\frak F}(z)$ also is $\K$-invariant. Thus
the function $J_\alpha f$
 is a left $\K$-invariant    on the group $G$
and hence $J_\alpha f$ is a function on $\G/\K$.

Now we want to find a $\K$-invariant function
${\frak F}$ from the equality
$$\langle {\frak F}(z), {\frak F}(z^{[g]})\rangle_{L^2(\G/\K)}=
 \langle \Xi(z), \Xi(z^{[g]})\rangle_{\V_\alpha}$$
 or in explicit form
 \begin{equation}
 \int\limits_{\B_{p,q}} {\frak F}(z)
 \ov{{\frak F}\bigl((a+zc)^{-1}(b+zd)\bigr)}
 \det(1-zz^*)^{-p-q}dz=\det(1-uu^*)^{\alpha}=
\prod(1+x_k)^{-\alpha}
\label{11.3}
 \end{equation}
 where $u=0^{[g]}=a^{-1}b$ and $x_j$ are the same as above
 (\ref{2.9}).

 \smallskip

 {\sc Theorem \fact.}
 {\it The function ${\frak F}(x)$ given by
 $${\frak F}(x)=
 \Lambda_{\G/\K}^\alpha(x):=
  \frac{\const}{\prod_{j=0}^{p-1}\Gamma(\alpha-j)}
 \det\limits_{k,j} \Bigl\{
   x_k^{1-2r} \frac{d^{j-1}}{dx^{j-1}_k} x^{2r+j-2}
   \Lambda_{r, r+j-1}^{\alpha-h}(x_k)\Bigr\}
   $$
   is a solution of the problem {\rm(\ref{11.3})}
    and hence the $\G$-intertwining
   operator $J_\alpha$ given by {\rm(\ref{11.2})} is unitary.}

   \smallskip

{\sc Remark.} This solution is not unique. I think
that it is the 'best' of possible solution (but this
 can be regarded as my own opinion,
 this subject is discussed in  \cite{Ner8}).
 \hfill $\square$

 \smallskip

 Let $B$ be a $p\times p$ matrix, let $x_j$ be its
 eigenvalues.
 We define   the {\it $\Lambda$-function of
 the symmetric space} $\G/\K=\U(p,q)/\U(p)\times\U(q)$ by
 $$\Lambda_{\G/\K}^\alpha(B):=\Lambda_{\G/\K}^\alpha(x)$$

 Formula (\ref{11.2}) written in an explicit form
 gives the following statement.

 \smallskip

{\sc Corollary \fact.} {\it The unitary $\G$-intertwining
 operator  $J_\alpha:L^2(\G/\K)\to \V_\alpha$ is given by
 $$J_\alpha f(u)=\int\limits_{\B_{p,q}}
  f(z) \Lambda_{\G/\K}^\alpha\bigl(
  Q(z,u)
  \bigr)
   \det(1-zz^*)^{-p-q}\{dz\}
 $$
 where}
  $$ Q(z,u):= (1-zz^*)^{-1}(z-u)(1-u^*u)^{-1}(z^*-u^*)$$

 {\sc Proof of Theorem 11.1.}
 Let $\xi_s$ be the spherical vector of  a spherical
 representation $\rho_s$ acting in the space $W_s$ (as above).

 Consider the  decomposition $\int W_s \, {\frak R}(s)ds$
 of $L^2(\G/\K)$ into the direct integral.
 The function $\frak F$ is $\K$-invariant. Hence  the
  image $\frak G$ of $\frak F$ in the direct integral is a function
of the form
 $s\mapsto\gamma(s)\xi_s$,
 where $\gamma(s)$ is the spherical transform of
 $\frak F$.
 Under the action of an element $g\in\G$ the function
 $\frak G$ transforms to
  $s\mapsto\gamma(s)\rho_s(g)\xi_s$.
  Hence the left-hand side of (\ref{11.3}) coincides
  with
  \begin{equation}
  \int |\gamma(s)|^2  \langle \xi_s,\rho_s(g)\xi_s\rangle_{W_s}
  {\frak R}(s)\,\{ds\}=
  \int  |\gamma(s)|^2 \Phi_s(g){\frak R}(s)\,\{ds\}
  \label{11.4}
  \end{equation}
  But the right-hand side of (\ref{11.3}) is the matrix
  element of the distinguished vector $\Xi_\alpha$.
 Hence $|\gamma(s)|^2$
  is the spherical transform of
  $\prod(1+x_k)^{-\alpha}$,
  and  hence $|\gamma(s)|^2$ is given by (\ref{6.18}).

  Now we assume
  $$\gamma(s)=\frac {1}{\prod_{j=0}^{p-1} \Gamma(\alpha-j)}
  \prod_{k=1}^p\Gamma(\alpha-(p+q-1)/2+s_k)
  $$
  (recall that $s$ is imaginary).
  It remains to evaluate the inverse spherical transform
  of $\gamma(s)$. This means
  that we must evaluate the integral
  \begin{multline*}
  \int\limits_{s_1\ge s_2\ge\dots\ge s_p\ge 0}
   \prod_{k=1}^p
   \Gamma(\alpha-h+is_k) \,\,\,
\det_{k,m}\bigl\{ \F(r+is_m,r-is_m;2r;-x_k)\Bigr\}
\times\\ \times
       \prod_{k=1}^p\frac
{\Gamma^2 (r+is_k) \Gamma^2 (r-is_k)}
                 {\Gamma(2is_k)\Gamma(-2is_k)}
                  \prod_{1\le l <m\le p} (s_l^2-s_m^2)
 ds_1\dots ds_p
 \end{multline*}
 The function
 $\Gamma(\alpha-h+is)$ is not even, hence we cannot change domain
of integration to $\R^p$. Nevertheless
the integrand is symmetric with respect to $s_k$
and we can replace the integral
by the integral
 $$\frac 1{n!} \int\limits_{s_1\ge 0,\dots,s_p\ge 0}$$
(with  the same integrand).

Let us convert the last factor
  of the integrand to the form
  $$\prod(s^2_l-s^2_m)=\det\limits_{k,j}
  \bigl\{ (b+is_k)_{j-1}(b-is_k)_{j-1}\bigr\}$$
  where $(b+is)_{j-1}$ is the Pochhammer symbol.
  By Lemma 6.10, the integral reduces to
  \begin{multline*}
 \det\limits_{k,j}\Bigl\{
  \int\limits_{0}^\infty
   \Gamma(\alpha-h+is) \,\,\,
 \F(r+is,r-is;2r;-x_k)
\times\\ \times
       \frac
{\Gamma (r+j-1+is) \Gamma (r+j-1-is)\Gamma (r+is) \Gamma (r-is)}
                 {\Gamma(2is)\Gamma(-2is)}
 ds\Bigr\}
 \end{multline*}
 Now we apply the formula (see \cite{HTF}, v.1, 2.8(22))
 $$(\gamma)_k y^{\gamma-1} \F(\alpha,\beta;\gamma;y)=
 \frac {d^k}{dy^k}
 y^{\gamma+k-1}\F(\alpha,\beta;\gamma+k;y)$$
 and obtain
 \begin{multline*}
 \const \cdot
\det\limits_{k,j}\Bigl\{ x^{1-2r}\frac {d^{j-1}}{dx^{j-1}}
x^{2r+j-2}
  \int\limits_{0}^\infty
   \Gamma(\alpha-h+is) \,\,\,
 \F(r+is,r-is;2r+j-1;-x_k)
\times\\ \times
       \frac
{\Gamma (r+j-1+is) \Gamma (r+j-1-is)\Gamma (r+is) \Gamma (r-is)}
                 {\Gamma(2is)\Gamma(-2is)}
 ds\Bigr\}
 \end{multline*}
 This completes the proof.\kvadrat

\smallskip

{\bf\punct. Hidden overgroup.}
The operator  $J_\alpha$ identifies the spaces $\V_\alpha$
and $L^2(\G/\K)$. We have seen in 4.2 that the
action of $\G$ in $\V_\alpha$ extends to the action
of the group $\widetilde \G=\G\times \G$ containing
$\G$ as the diagonal subgroup.  If we believe that our
operator $J$ is canonical, then we obtain the following
strange statement. {\it There exists
a natural one-parametric family {\rm(}depending of
$\alpha>h${\rm)}
of  actions
of the group $\G\times\G$ in the space $L^2(\G/\K)$.}

\smallskip

{\small
{\bf\punct. Comments.}
1. The problem of explicit unitary identification
goes back to the work \cite{OO}, in this paper there was
called attention to the coincidence of spectra
of $L^2(G/K)$ and restriction of a highest weight
representation of $\widetilde G$
to $G$.

 2. We define a
$\Lambda$-function of a classical
Riemannian noncompact
 symmetric space $G/K$ by
$$\Lambda^a_{G/K}(x)=
\int
\prod_k \Gamma(a+is_k)\Phi_s(x){\frak R}(s)ds$$
where $\Phi_s$ are spherical functions
of $G$ and $\frak R$ is the Gindikin--Karpelevich density.
For six series of the classical groups the
function $\Lambda^\alpha_{\G/\K}$ can be evaluated
explicitly.

\smallskip

3. {\it The case $G=\GL(q,\C)$, $\GL(q,\R), \GL(q,\HH)$}
is trivial. For definiteness, consider
 the kernel representation of $\GL(q,\C)$.
It is
 the restriction of the representation $\tau_\alpha$
(see (3.2))
of the group $\U(q,q)$ to the subgroup $\GL(n,\C)$
(see lists in \cite{Ner3}, \cite{Ner8}).
The symmetric space $\GL(q,\C)/\U(q)$
is the cone   $\Pos_q$ defined in \Subsection 3.6.
The Laplace transform is a unitary $\GL(q,\C)$-intertwining
operator from $L^2(\Pos_q)$ to the space $H_\alpha(\W_q)$
of the kernel representation.

\smallskip

4. The calculation given in 11.3 is valid also for
the  series
$G=\OO(n,\C)$ and $\Sp(2n,\C)$.

       }

\bigskip

{\large \bf Addendum. Pseudoriemannian symmetric spaces,
Berezin forms, and some problems of non $L^2$ harmonic analysis}

\medskip

\nopagebreak
\newcounter{add}
 \renewcommand{\theequation}{\Alph{add}.\arabic{equation}}
\setcounter{add}{1}
\setcounter{equation}{0}

\def\KK{{\Bbb K}}

{\small
 Molchanov (see \cite{Mol0}, \cite{Mol1}, \cite{MD1}, \cite{MD2})
 introduced
some types of representations   of certain groups $G$
("canonical representations") related to  expressions similar
to Berezin kernels (4.7).
I shall try to give a general scheme for pseudo-Riemannian
symmetric spaces, including
kernel representations together with
Molchanov's, van Dijk's, and others  examples   of
"canonical representations".
I shall  also try to  discuss similarities
and differences
   of these constructions
and the kernel representations (A.5--A.8).

In A.9--A.10 another (essentially different) problem of
harmonic analysis on pseudo-Riemannian symmetric
spaces is considered.

We need  uniform models of
all the classical pseudo-Riemannian
symmetric spaces obtained in \cite{pse}.
 First we consider some examples
and after this we describe the general construction.

\smallskip

{\bf A.1. Several examples.} In all examples below a point of a symmetric space
$G/H$ is an ordered pair of linear subspaces $W$, $Y$ in a fixed linear space $V$
such that $V=W\oplus Y$.

\smallskip

{\sc Example 1.}  The space
$\GL(n,\R)/\GL(k,\R)\times\GL(n-k,\R)$
 consists of pairs $(W,Y)$ of subspaces  in $\R^n$ such that
$\dim W=k$, $\dim Y=n-k$, $\R^n=W\oplus Y$.

\smallskip

{\sc Example 2.} Consider the linear space $\C^{2n}$ equipped with
the operator $J$ of the complex conjugation
$$J(z_1,\dots,z_{2n})=(\bar z_1,\dots,\bar z_{2n})$$
The group $G$ of all linear operators commuting with
$J$ is $\GL(2n,\R)$. A point of the symmetric space
$\GL(2n,\R)/\GL(n,\C)$ is a pair of subspaces $(W,Y)$ in $\C^{2n}$
such that $JW=Y$ (and hence $JY=W$) and $\C^{2n}=W\oplus Y$.

\smallskip

{\sc Example 3.} Consider the linear space $\R^{2n}=\R^n\oplus\R^n$
equipped with
the operator $J$ with the matrix
$\bigl(\begin{smallmatrix}1&0\\0&-1\end{smallmatrix}\bigr)$.
The group of all operators in $\R^{2n}$ commuting with $J$ is
$G=\GL(n,\R)\times\GL(n,\R)$.
Consider the set $\frak S$  of all pairs of subspaces $(W,Y)$ in $\R^{2n}$
such that $JW=Y$ (and hence $JY=W)$   and $W\oplus Y=\R^{2n}$.
Obviously,
$${\frak S}=G/H=\GL(n,\R)\times\GL(n,\R)/\GL(n,\R)$$
where the subgroup $H=\GL(n,\R)$ is the diagonal subgroup in $G$.

\smallskip

{\sc Example 4.} Consider the linear space $\C^{p}\oplus \C^q$
equipped with the Hermitian  form $Q$ with the matrix
$\bigl(\begin{smallmatrix}1&0\\0&-1\end{smallmatrix}\bigr)$
as above (Section 2). Consider the set $\frak S$  of all linear
subspaces $W\subset \C^{p}\oplus \C^q$  such that $\dim W=r$
and the form $Q$ is nondegenerate on $W$. Obviously
$$
{\frak S}=\bigcup\limits_{s,t: s+t=r, s\le p, t\le q}
\U(p,q)/\bigl(\U(s,t)\times\U(p-s,q-t)\bigr)
$$
(this construction includes the construction 2.4 above).
Let us also introduce  subspaces $Y$ which are the orthocomplements
$W^\bot$
to $W$.

\smallskip

{\sc Example 5.}
Consider the space
$\R^{2n}=\R^n\oplus\R^n$ equipped with the skew-symmetric
bilinear form with the  matrix
$\bigl(\begin{smallmatrix}0&1\\-1&0\end{smallmatrix}\bigr)$.
Consider the space $\frak S$ of pairs $(W,Y)$ of maximal isotropic subspaces
in $\R^{2n}$ such that $\R^{2n}=W\oplus Y$.
Obviously,
$${\frak S}=\Sp(2n,\R)/\GL(n,\R)$$

\smallskip

{\sc Example 6.}
Consider the space
$\R^{2n}=\R^n\oplus\R^n$ equipped with the skew-symmetric
bilinear form $B$ with the  matrix
$\bigl(\begin{smallmatrix}0&1\\-1&0\end{smallmatrix}\bigr)$
and the symmetric bilinear form
$Q=\bigl(\begin{smallmatrix}1&0\\0&-1\end{smallmatrix}\bigr)$.
Consider the space $\frak S$ of pairs $(W,Y)$ of maximal $B$-isotropic subspaces
in $\R^{2n}$ such that $\R^{2n}=W\oplus Y$
and $Y$ is the orthocomplement of $W$ with respect to $Q$.
It can easily be checked that
$${\frak S}=\bigcup\limits_{r=0}^n
\GL(n,\R)/\OO(r,n-r)$$

{\bf A.2. Uniform construction of classical pseudo-Riemannian symmetric spaces.}
In all examples above a point of a pseudo-Riemannian symmetric space  $G/H$
is a pair $(W,Y)$
 of transversal subspaces  satisfying  some simple conditions.
There are 3 types of conditions:

a) $W,Y$ are maximal isotropic subspaces

b) $Y$ is the orthocomplement of $W$

c) $Y=JW$, $W=JY$ for some fixed operator $J$.

\smallskip

For a uniform description
of all 54 series of  pseudo-Riemannian classical symmetric spaces
(see Berger classification, \cite{Berger})
we must fix   some notations and definitions from linear algebra
(this general construction is not necessary for understanding
Subsection A3--A10 below).

The term {\it linear space} below means a right finite-dimensional
module $V$ over $\R$, $\C$ or the quaternions $\HH$.

\smallskip

A {\it semiinvolution} $J$  is a linear or antilinear%
\footnote{A map $A:V\to V$ is called an antilinear operator
if $A(v+w)=Av+Aw$,  $A\mu v=\overline\mu Av$ for all $v,w\in V$,
$\mu\in \KK$.}
operator in a linear space $V$ satisfying the condition
$J^2=\pm 1$.\footnote{This definition is adapted to our
field $\KK=\R,\C\,\HH$, for a general definition of semiinvolution see
\cite{Died}.}
A semiinvolution $J$ in $V$ is {\it split} if there exists a
 subspace $W$ such that $V=W\oplus JW$.

\smallskip

The term {\it form} below means a nondegenerate form
on a linear space $V$ over $\KK=\R,\C,\HH$
of one of the following types

a) a symmetric or skew symmetric form over $\R$

b) a symmetric, skew symmetric or Hermitian form over $\C$

c) an Hermitian or anti-Hermitian
\footnote{A sesquilinear form $B(v,w)$ is called anti-Hermitian
if $B(w,v)=-\overline{B(v,w)}$.}
 form over $\HH$

\smallskip

\smallskip

A linear semiinvolution  $J$ is {\it consistent} with a form
$B(\cdot,\cdot)$ if $B(Jv,Jw)=\pm B(v,w)$.
 An antilinear semiinvolution $J$ is consistent with
$B$ if $B(Jv,Jw)=\pm \overline{B(v,w)}$.

\smallskip

A subspace  $W\subset V$ is called {\it isotropic}
with respect to a form $B(\cdot,\cdot)$ if $B(w,w')=0$ for all
$w,w'\in W$.

A form $B$ on $V$ is {\it split} if there exists
an isotropic subspace $W\subset V$ such that $\dim W=\frac 12\dim V$.

A {\it Grassmannian} is a set of all subspaces of a given dimension
 in a  linear space or  a set of all isotropic
 subspaces of a given dimension.

A real {\it classical group}  is a group of all
linear operators in a linear space over $\KK$ (i.e.,
$\GL(n,\R)$, $\GL(n,\C)$, $\GL(n,\HH)$) or the group
preserving some form in a linear space over $\KK$, i.e.,

-- the groups $\OO(p,q)$, $\Sp(2n,\R)$ over $\R$

-- the groups $\OO(n,\C)$, $\Sp(2n,\C)$, $\U(p,q)$ over $\C$

-- the groups $\Sp(p,q)$, $\SOS(2n)$ over $\HH$.

\smallskip

A {\it classical pseudo-Riemannian symmetric space} is a
homogeneous space of the form $G/H$, where $G$ is a classical
group, and a subgroup $H$ is the set of fixed points
of some automorphism  $\sigma$ of $G$ satisfying $\sigma^2=1$.

\smallskip

All classical symmetric spaces (up to covering and center)
can be obtain in the following way.

Consider the following types of data on a linear space $V$
over the field $\KK$.

STRUCTURE $1^\star$. A {\it basic form} $B$ is a split form on $V$.

STRUCTURE $2^\star$. A  {\it control semiinvolution} $J$
is a split semiinvolution on $V$.

STRUCTURE $3^\star$. A {\it control  form} is  arbitrary form on $V$.

Consider
 ordered pairs $(W,Y)$ of subspaces of $V$ such that
\begin{equation}
V=W\oplus Y
\end{equation}

We say that a pair $(W, Y)$   is {\it consistent with the basic} form
$B$ if the both subspaces $W,Y$ are $B$-isotropic.

We say that  a pair $(W,Y)$   is {\it consistent with the
semiinvolution} $J$ if $J W=Y$.

We say that  a pair $(W,Y)$   is {\it consistent with
the control form} $C$ if $Y$ is the orthocomplement to $W$
with respect to the form $C$.

We consider a linear space $V$ without any additional structure,
or a linear space $V$  equipped with one of the structures $1^\star-3^\star$,
or a linear space $V$  equipped with all structures  $1^\star-3^\star$.
In the last case we assume $J$ is consistent with $B$
and
$$C(v,w)=B(Jv,w)$$

\smallskip

{\sc Lemma A.1.} {\it If a pair  $(W,Y)$  is consistent
with two of the structures   $1^\star-3^\star$, then it is
consistent with the third structure.}

\smallskip

Fix a linear space $V$ with such structure.
Denote by $G(V)$ the groups of all linear operators
in $V$ preserving all structures on $V$.
Consider the  set $\frak S$ of all pairs $(W,Y)$ such that

1) $W\oplus Y=V$

2) $(W,Y)$ is consistent with the structure of $V$

3) $\dim Y$ is fixed

\smallskip

{\sc Observation A.2.} \cite{pse} a) {\it If the set $\frak S$
is nonempty, then it
  is a pseudo-Riemannian symmetric space
$G(V)/H$ or atheunion of a finite collection of
pseudo-Riemannian symmetric spaces
$G(V)/H_j$}

\smallskip

 b) {\it Each pseudo-Riemannian symmetric space
can be obtained in this way.}

\smallskip

Tables are contained in \cite{pse}. The examples given above
give  a
representative sample.

\smallskip

{\bf A.3. Matrix atlas.} Consider a symmetric space $G/H$
obtained in this way. Fix a point $(\widetilde W,\widetilde Y)\in G/H$.
Then for $(W,Y)\in G/H$ in  general position, the subspace
$W$ is a graph of an operator
$A:\widetilde W \to \widetilde Y$
and $Y$ is a graph of an operator $B:\widetilde Y \to \widetilde W$.
Thus for any point of $G/H$ we obtain a pair of operators
$$A:\widetilde W \to \widetilde Y;\qquad
B:\widetilde Y \to \widetilde W$$
Thus for any point  $(\widetilde W,\widetilde Y)$ we obtain a
coordinate system $(A,B)$
on $G/H$.

\smallskip

{\sc Example.} For the spaces $\U(p,q)/\U(p)\times \U(q)$
 we obtain
the Cartan matrix ball realization.

\smallskip

{\sc Remark.} The matrices $A,B$ are not arbitrary. For instance in Example 1
of A.1 a pair $(A,B)$  satisfy the unique condition $\det(1-AB)\ne 0$.
In Example 5 we have the same condition and also $A=A^t$, $B=B^t$.
Generally the pair $(A,B)$ ranges in an open subset in some linear
subspace in space of pairs of matrices.

\smallskip

{\bf A.4. Hua Loo Keng double ratio.}
 Fix two pairs of subspaces $(W_1,Y_1)$, $(W_2,Y_2)\in G/H$
in general position  (we also assume $\dim W_j\le\dim Y_j$).
 Then $W_2$ is the graph of some operator
$R:W_1\to Y_1$, and $Y_2$ is the graph of some operator
$S:Y_1\to W_1$. We define the Hua double ratio
operator
$${\cal D}(W_1,Y_1;W_2,Y_2):=SR:\quad W_1\to W_1$$
In  matrix coordinates, this operator is given by the formula
$${\cal D}(A_1,B_1;A_2,B_2)=(1-B_2A_1)^{-1}(B_1-B_2)(1-A_2B_1)^{-1}(A_1-A_2)$$

{\sc Lemma A.3.} a) {\it 1 is not an eigenvalue of
${\cal D}(W_1,Y_1;W_2,Y_2)$.}

\smallskip

  b) $1-{\cal D}(A_1,B_1;A_2,B_2)=(1-B_2A_1)^{-1}(1-B_2A_2)(1-B_1A_2)^{-1}(1-B_1A_1)$

\smallskip

c) ${\cal D}/(1-{\cal D})=
  (1-B_1A_1)^{-1}(B_2-B_1)(1-A_2B_2)^{-1}(A_2-A_1)$

\smallskip

{\sc Proof.} Statement a) is equivalent to $W_2\cap Y_2=0$,
and b), c) can checked by a simple calculation.

\smallskip

{\bf A.5. Berezin form.} Consider an arbitrary character
$\chi$ of the multiplicative group of $\KK$
(for instance, $\chi(z)=|z|^\alpha$). Consider
the kernel
$$
{\cal L}_\chi(W_1,Y_1;W_2,Y_2)=
\chi\bigl(\det\bigl[1- D(W_1,Y_1;W_2,Y_2)\bigr]\bigr)
$$
We define the Berezin form
on the space of smooth
compactly supported functions on $G/H$
 by
\begin{equation}
\langle f,g\rangle:=
\iint_{G/H\times G/H} {\cal L}_\chi(W_1,Y_1;W_2,Y_2)
f_1(W_1,Y_1)\overline{f(W_2,Y_2)}\,d\mu(W_1,Y_1)\,d\mu(W_2,Y_2)
\end{equation}
where $(W_1,Y_1)$, $(W_2,Y_2)$ are points of the symmetric
space and $\mu$ is the $G$-invariant measure on $G/H$.

\smallskip

{\bf A.6. Comparison of Riemannian and pseudo-Riemannian cases.}
Formally the construction A.5 in the Riemannian case
gives kernel representations. But a serious divergence
between the Riemannian and pseudo-Riemannian cases
appears immediately.

I emphasis that the kernel ${\cal L}_\chi$ is smooth
on the diagonal of $G/H\times G/H$ but
(for nonRiemannian case) it has singularities
outside the diagonal.
 It seems (but not proved carefully) that our Hermitian form
always (except for the Riemannian case) is indefinite. This
leads to serious technical difficulties.
In particular, {\it an  indefinite Hermitian form does not
define a topology in a functional space}.
So even the formulation of the specral problem is not obvious,
for a discussion see  A.7, A.8 below.

Many other phenomena existing
for the kernel representations do not survive
in the pseudo-Riemannian case. For instance, there is no realization
in holomorphic functions, there is no overgroup $\widetilde G$
described in 4.7 and Section 11%
\footnote{$\widetilde G$ has no relation to the overgroup $G^\circ$
described below}
etc. It seems that theory of the kernel representations
 cannot be a special case
of pseudo-Riemannian theory.

\smallskip

{\bf A.7. Indefinite harmonic analysis: approach related to
 Krein structures.}
  Consider  a linear  space
$X$ equipped with an indefinite Hermitian form
$Q$. Fix subspaces $X_+$, $X_-$ in $X$ such that the form $Q$
is positive definite on $X_+$, negative definite on $X_-$ and
$X=X_+\oplus X_-$. Let $\bar X_+$ be the completion of the pre-Hilbert
space $X_+$, let $\bar X_-$ be the completion of the prehibert space
$X_-$ with respect to the form $(-Q)$. Thus we obtain
a topological vector space $\bar X:=\bar X_+\oplus \bar X_-$ equipped with
the form $Q$ and with the fixed decomposition $\bar X_+\oplus \bar X_-$
(these datas are called the {\it Krein structure}, \cite{AI}).

A Krein structure is not canonically determined by the space $X$ and the
 form Q.\footnote{Let us decribe the simplest example.
Consider the space $X$ consisting of finite linear combinations
of vectors $e_1$, $e_2$,\dots; $f_1$,$f_2$,\dots. Assume
all these vectors are pairwise orthogonal and $\langle e_j,e_j\rangle=1$
and $\langle f_j,f_j\rangle=-1$. Consider the subspace
$X_+$ generated by $e_j$ and the subspace  $X_-$ generated by $f_j$.
Consider also the subspace
$X_+'$ generated by the vectors
$\sqrt{j+1}e_j+\sqrt{j}f_j$ and the subspace
$X_-'$ generated by the vectors
$\sqrt{j}e_j+\sqrt{j+1}f_j$.
Then  $\bar X_+\oplus \bar X_-$  and $\bar X_+'\oplus  \bar X_-'$
are different linear topological spaces containing $X$, i.e.,
the identical operator $X\to X$ cannot be extended to  a
bounded bijection $\bar X_+\oplus \bar X_-\to\bar X_+'\oplus  \bar X_-'$.}

It is not clear is  this approach useful in representation theory
of semisimple groups
or not. It seems that even a problem of existence of
nontrivial examples for groups of rank $>1$ is open.%
It seems to me that the following question is way to
understand this.

\smallskip

 Let $G$ be a semisimple group and $K$ be its
compact subgroup. Let an irreducible representation
(a Harish-Chandra module) $\rho$ of $G$ in the space $X$
 admit a $G$-invariant
 Hermitian form $Q$ and let the restriction of $\rho$ to
$K$ be mutiplicity free. Then $X$ admits
a canonical Krein structure
(since the restriction of $Q$ to any  irreducible
 $K$-subrepresentation
is positive definite or negative definite).
The simplest nontrivial examples of such a picture are

 1. the representations $\tau_\alpha(g)$
of the group  $\G=\U(p,q)$
for any noninteger real $\alpha$.

 2. Molchanov's degenerate
 representations of $\OO(p,q)$ (see \cite{cone}).

\smallskip

There arises the following problem.

\smallskip

{\sc Question A.4.} a) {\it Is it possible
in these two cases to write
the projections to the subspaces $\bar X_\pm$ explicitly?}

\smallskip

b) {\it Are the operators of the representation
continuous in the topology of the Krein space?}

\smallskip

{\bf A.8. Indefinite harmonic analysis. Molchanov's approach.}
This approach is not well-formilized, but there is a nice collection
of explicit nontrivial examples.%
\footnote{First example of this kind  were observed on "physical level"
in \cite{Mol01}, for mathematically rigorous way of formulation
of such problems see \cite{MD1}), see also
\cite{Mol0}, \cite{Mol1}, \cite{MD2}}
They have the following form.

Consider a representation $\zeta$ of a group $G$ in a topological
vector space $X$. Let $Q$ be a $G$-invariant Hermitian form
on $X$ (the basic  example is described above in A.5).

Also consider  some family $\rho_t$ of irreducible representations of
$G$ in the spaces $Y_t$, and assume that each representation
of this family admits a $G$-invariant Hermitian forms $R_t(\cdot,\cdot)$.
Consider a space $\cal Y$
 of functions $f$ that takes any $t$ to a vector $y_t\in Y_t$.
The space of such functions is some kind of a direct integral of
representations, but it is not a direct integral in the formal common
sense.

Consider a $G$-intertwining operator from $X$ to $\cal Y$.
The operator $J$ takes any vector $x\in X$ to some function
$Jx(t)$. An {\it indefinite Plancherel formula}
is the identity
$$Q(x_1,x_2)=\int R_t\bigl(Jx_1(t),Jx_2(t)\bigr)\,d\mu(t)$$
where $d\mu(t)$ is some ("Plancherel") measure.

In particular, this measure is obtained
for Berezin forms on certain rank 1 pseudo-Riemannian
symmetric spaces
(\cite{Mol1} \cite{MD1}, \cite{MD2}).

Fix arbitrary (for simplicity noninteger) $\alpha\in\R$.
Consider the space of smooth functions on $\B_{p,q}=\G/\K$
equipped with the (generally indefinite)
Hermitian form (4.11).

\smallskip

{\sc Conjecture A.5.} {\it
Our Plancherel formula
{\rm(}see Section {\rm 7)} is valid for arbitrary  $\alpha\in \R$
in Molchanov's sense}.

\smallskip

{\bf A.9. Conformal group of the symmetric space.}
Each classical pseudo-Riemannian symmetric space
$G/H$ admits a canonical open  embedding
to some space  $G^\circ/P^\circ$, where
the {\it conformal group}%
\footnote{The term "conformal" was introduced by Goncharov and Gindikin}
 $G^\circ\supset G$ is
a classical group  and $P^\circ\supset H$ is
some maximal parabolic subgroup in $G^\circ$%
\footnote{This fact
can be extracted from Makarevich's tables \cite{Mak},
but it never was claimed before \cite{pse};
for exceptional spaces analogy of this statement,
 in general, is false}.  In fact the space
$G^\circ/P^\circ$ is a Grassmannian or a product
of two Grassmannians.


\smallskip

{\sc Example.} In the situation described in 2.4
($G/H=\G/\K=\U(p,q)/\U(p)\times \U(q)$) the space
$G^\circ/P^\circ$ is the Grassmannian of $p$-dimensional
subspaces and $G^\circ=\GL(p+q,\C)$.

\smallskip

{\sc Example.} The conformal group in Example 1 in A.1 is
$G^\circ=\GL(n,\R)\times\GL(n,\R)$, in Example 2
we have $\GL(2n,\C)$,
in Example 3 we have $G^\circ=\GL(2n,\R)$,
 in Example 4 we have $G^\circ=\GL(p+q,\C)$,
in Example 5 we have $G^\circ=\Sp(2n,\R)\times\Sp(2n,\R)$, in Example 6
we have
$G^\circ=\Sp(2n,\R)$.

\smallskip

{\sc Remark.} The general algorithm for obtaining the
conformal group
$G^\circ$ is following.

1. Let
a control semiinvolution and a control form  be absent.
Then $G^\circ=G\times G$.  The space $G^\circ/P^\circ$
is the product of $B$-isotropic Grassmannians
if a basic form $B$ is present, and a product of
the usual Grassmannians
otherwise.

2. Otherwise, we "forget"   the control semiinvolution and the control form
on $V$ and consider the group of automorphisms
of the basic form $B$ if the basic form is present,
and the complete linear group if $B$ is absent.

\smallskip

{\sc Remark.} For 44 series of symmetric spaces,
 the subset $G/H$ is dense
in $ G^\circ/P^\circ$. For 10 series
$G/H$ is not dense in  $ G^\circ/P^\circ$ but
the group $G$ has a finite number of open orbits in $ G^\circ/P^\circ$,
all these orbits are symmetric spaces of the form $G/H_j$,
and the union of $G/H_j$ is dense in $ G^\circ/P^\circ$.
These 10 series are distinguished by the condition:
a control form is present and it is an Hermitian form
\footnote{An Hermitian form over $\R$ is a symmetric bilinear form.}

\smallskip

For Riemannian symmetric spaces  $G/H$ (i.e., $G/K$
in the notation of Section 11)
 {\it the conformal group $G^\circ$ and the hidden overgroup
$\widetilde G$ discussed in Section 11 are different objects.}
Hidden overgroup acts by  integral operators in $L^2$
and it does not act (even locally) on the symmetric space $G/K$
itself.

\smallskip

{\bf A.10. Deformation of $L^2$ on pseudo-Riemannian symmetric space.}
A natural representation $\rho_0$ of $G^\circ$ in $L^2(G^\circ/P^\circ)$
is a representation of a principal degenerate series.

\smallskip

{\sc Observation A.6.}(\cite{pse}) a) {\it
If $G$ has a dense orbit in  $G^\circ/P^\circ$, then
the restriction of $\rho_0$ to $G$ is $L^2(G/H)$.
Otherwise   this restriction is a direct sum
of several spaces $L^2(G/H_j)$.}

\smallskip

In many cases the representation $\rho_0$ can be included
in a degenerated complementary series $\rho_s$. Thus there arises
a problem of Plancherel formula for the restriction of  $\rho_s$
to $G$. This restriction is some kind of deformation of $L^2(G/H)$.

\smallskip

This problem survives even in the case in which the complementary series
is absent. Consider the
 $G$-invariant kernel on pseudo-Riemannian
symmetric space $G/H$  given by
$${\cal L}_\chi(W_1,Y_1;W_2,Y_2):=
\chi\Bigl(\det \frac{{\cal D} (W_1,Y_1;W_2,Y_2)}
  {1-{\cal D} (W_1,Y_1;W_2,Y_2)}\Bigr)$$
Then the exprssion (A.2) is an Hermitian form on the space
of compactly supported smooth functions on $G/H$.

\smallskip

{\sc Example.} For $G/H= \U(p,q)/\U(p)\times \U(q)$
we obtain the $\U(p,q)$-invariant kernel
$$\Bigl|\det(1-z^*z)^{-1}(z^*-u^*)(1-uu^*)^{-1}(z-u)\Bigr|^\mu=
\Bigl(\frac{\det(z-u)(z^*-u^*)}{\det(1-zz^*)\det(1-uu^*)} \Bigr)^\mu$$
This kernel has the main singularity on the diagonal $z=u$.
Evidentely the representation  of $\U(p,q)$  associated with
this kernel is not a kernel representation.

                  }

\def\itt{\it}

\sc Insitute of theoretical and experimental physics, Moscow, Russia

\sc Independent University of Moscow,
Bolshoi Vlas'evskii, 11,
Moscow -- 121002, Russia


\tt neretin@main.mccme.rssi.ru


\begin{thebibliography}{c}






\bibitem{AAR}
Andrews, G.R.,  Askey, R., Roy, R.,
{\it Special functions.}, Cambridge Univ. Press, 1999

\bibitem{AI}
Azizov, T.Ya., Iohvidov, I.S.,
{\it Foundations of the theory of linear operators
in spaces with indefinite metrics,}
Nauka, Moscow, 1986; English translation Wiley, New York, 1989.

\bibitem{Ber1}      
Berezin, F.A., {\it Quantization in complex symmetric spaces.}
Izv. Akad. Nauk SSSR, Ser. Math., 39, 2, 1362--1402 (1975);
English translation: Math USSR Izv. 9 (1976), No 2, 341--379(1976)

\bibitem{Ber2}     
Berezin, F.A. {\it On relations between covariant and contravariant
symbols of operators
for complex classical domains.}
 Dokl. Akad Nauk SSSR, 241, No 1 (1978), 15--17;
English translation: Sov. Math. Dokl. 19 (1978), 786--789

\bibitem{BK}
 Berezin, F.A., Karpelevich, F.I.,
{\it Zonal spherical functions and Laplace operators on
 some symmetric spaces}. Dokl. Akad. Nauk SSSR, 118 (1958), 9--12

\bibitem{Berger}
Berger, M., {\it Les espaces simmetrique noncompact.}
Ann. Sci. Ecole Norm. Super. 74 (1957), 85--177


\bibitem{BO} Borodin, A., Olshanski, G.,
{\it Infinite random matrices and ergodic measures,}
Preprint {\tt arXiv.org/abs/math/0010015}

\bibitem{Die}
van Diejen, J.F.,
{\it Properties of some families of hypergeometric polynomials
in several variables.}
Trans. Amer. Math. Soc., 351 (1999), 233--270

\bibitem{DS} van Diejen, J.F., Stockman, J.V.,
{\it Multivariate $q$-Racah polynomials.}
Duke M.J. 91 (1998),   89--136.

\bibitem{Died}
Dieudonne, J., {\it La geometrie des groupes classiques,}
Springer, 1971.

\bibitem{vD}          
van Dijk, H., Hille, S.C.,
 {\it Canonical representations related to hyperbolic
spaces},
J.Funct.Anal., 147, 109--139   (1997).


\bibitem{MD1}         
van Dijk, G., Molchanov, V.F.,
{\it The Berezin form for rank 1
 para-Hermitian symmetric spaces.}
J.Math.Pure.Appl.,IX ser, 77, N1, 747--799 (1998)



\bibitem{MD2}         
van Dijk, G., Molchanov, V.F.,
{\it Tensor products of maximal degenerate series representations
of the group ${\rm SL}(n,\R)$.}
J. Math. Pure Appl., 78 (1999), 99--119

\bibitem{Dix} Dixmier, J. {\it Les $C^*$-algebres et leurs representations,}
Paris, Gauthier--Villars, 1969.

\bibitem{FK}
Faraut, J., Koranyi, A.,
 {\it Analysis in symmetric cones.}
Oxford Univ.Press, (1994)


\bibitem{GG} Gelfand, I.M., Gindikin, S.G.,
{\it  Complex manifolds whose skeletons are real semisimple
Lie groups, and analytic discrete series of representations,}
Funkt. Anal i Prilozh., 11 (1977), 4, 19--27(Russian);
English translation in Funct. Anal. Appl., 11 (1977),
258--265.

\bibitem{GN} Gelfand, I.M., Naimark M.I.,
 {\it Unitary representations of classical groups.}
 Trudy Mat.Inst. Steklov.,v.36 (1950) (Russian);
 German translation: Gelfand I.N., Neumark M.A.,
 {\it Unitare Darstellungen der klassischen gruppen.},
 Akademie-Verlag, Berlin, 1957.


\bibitem{Gin1}       
Gindikin, S.G., {\it Analysis on homogeneous spaces}.
Uspehi mat. nauk,19, No 4, 3--92(1964);

\bibitem{Gin2}        
Gindikin, S.G.,
{\it Invariant distributions in homogeneous domains.}
Funkt. Anal. Prilozh. 9 (1975), No 1, 56--58;
 English translation:
Funct. Anal. Appl. 9 (1975), No.1, 50--52.

\bibitem{GK1}          
Gindikin, S.G., Karpelevich, F.I., {\it Plancherel measure for
Riemannian symmetric spaces of non-positive curvature.}
 Dokl. Akad Nauk SSSR,
 145, 252--255(1962); English translation Sov.Mat.Dokl,3, 962--965
 (1962).

\bibitem{GK2}            
Gindikin, S.G., Karpelevich, F.I.,
{\it On an integral connected with symmetric
 Riemannian space of non-positive curvature.}
 Izv. Akad. Nauk SSSR, Ser Mat., 30, 1147--1156(1966)
(Russian); English translation in
 Transl.Amer.Math.Soc., 85, 249--258(1969)


 \bibitem{GR}
 Gross, K.I., Richards, D.S.P.,
 {\it Total positivity, spherical series,
 and hypergeometric functions of
 matrix argument.}
 J. Approx. Theory 59 (1989), 224--246.

\bibitem{Gut}
Gutkin E.
 Coefficients of Clebsch-Gordan
for holomorphic discrete series.
Lett.Math.Phys. 3(1979), 185--192


\bibitem{HO}
 Heckman, G.I., Opdam, E.M.,  {\it
 Root systems and hypergeometric
functions}.I. Compositio Math, 64(1987), 329--352;
II. Compositio Math., 64 (1987), 353--373;
III. Compositio Math., 67 (1988), 21--49;
IV.  Compositio Math, 67(1988), 191--209.

\bibitem{Hel0}      
Helgason, S.,
{\it Group representations and symmetric spaces.}
Adv. Math., 5 (1970), 1--154.

\bibitem{Hel1}   
 Helgason, S.,
 {\it Differential geometry and symmetric spaces.}
 Acad. Press., New York, 1962.



\bibitem{Hel2}           
Helgason, S., {\it Groups and geometric analysis.}
 Acad. Press (1984).

\bibitem{Her}
Herz, C.S.,
{\it Bessel functions of matrix argument.}
Ann. of Math.(2) 61 (1955), 474--523.

\bibitem{Hil} Hille, S.C., {\it Canonical representations.}
Ph.D Thesis, Leiden University, June 1999.

\bibitem{HTF}            
{\it Higher transcendental functions}, v.1-3,
McGraw-Hill book company, 1953-1955

\bibitem{Hoo}
Hoogenboom, B., {\it Spherical functions
 and invariant differential
operators on complex Grassmann manifolds.} Ark.  Math.,
20(1982), 69--58.

\bibitem{Hua}
Hua Loo Keng, {\it Harmonic analysis
of functions of several complex
variables in classical domains%
}. Beijing, 1958(Chinese);
Russian translation: Moscow, Inostrannaya literatura (1959);
English translation: Amer. Math. Soc., Providence (1963)

\bibitem{JV}
Jakobsen, H.P., Vergne, M.,
{\it Restrictions and expansions of holomorphic
representations.}
J. Funct. Anal., 34 (1979), 29--53.


\bibitem{Kir} Kirillov, A.A., {\it Elements of
representation theory.} Nauka, Moscow, 1972; English translation:
                            Springer,1975;
                            French translation: Mir, Moscow, 1974




\bibitem{KS}
Koekoek, R., Swarttouw, R.F.,
{\it Askey scheme of hypergeometric orthogonal
polynomials and its $q$-analogues.} Delft University
of technology, 1994, available at
{\tt http:/aw.twi.tudelft.nl/$\widetilde{\phantom{a}}$%
koekoek/research.html}

\bibitem{Koo1}   
Koornwinder, T.H.,
{\it  Jacobi functions and analysis on noncompact symmetric
spaces.} in {\it Special functions: group theoretical aspects
and applications,}
eds. Askey, R., Koornwinder T.H., Schempp W., 1--85,
D. Reidel Publ. Co., Dordrecht--Boston, 1984.


\bibitem{Koo2}     
Koornwinder, T.H.,
{\it Askey--Wilson polynomials for root systems
of the type $BC$} in {\it Hypergeometric functions on domains of
positivity, Jack polynomials and applications}, 189--204,
Contemp. Math., 138, AMS, 1992

\bibitem{Kre} Krein, M.G., {\it
 Hermitian positive definite kernels on homogeneous spaces,}
I, II. Ukrain. Math. J., 1 (1949), 4, 64--98; 2(1950), 1, 10-59
(Russian); English translation in Amer.Math.Soc.Transl.




\bibitem{Mac} Macdonald, I.G.,
{\it Symmetric functions and Hall polynomials,}
Oxford University Press, 1995.

\bibitem{Mak}
Makarevich,  B.O.
{\it Open orbits of reductive groups in symmetric $R$-spaces.}
Mat. Sbornik, 91, 3 (1973), 390--401, English translation in
Math. USSR Sbornik, 20, 3 (1973), 406--418

\bibitem{Mol}    
 Molev, A.I., {\it Unitarizability of some Enright--Varadarajan
$U(p,q)$-modules.} in
{\it Topics in representation theory}, ed. A.A.Kirillov
(1991), 199--220, Adv. Sov. Math., v.2.,
 Amer. Math. Soc. Translations, 1991.

\bibitem{cone}
Molchanov, V.F.,
{\it Representations of the pseudoorthogonal group
related to cone,} Mat. Sbornik, 81 (1970), 3,358--375;
English translation: Math. USSR Sbornik, 10 (1970), 333--347

\bibitem{Mol0}
Molchanov, V.F.,
{\it Quantization on the imaginary Lobachevsky plane.}
Funkts. Anal. Prilozh., 14, 2 (1980), 73--74; English translation
Funct. Anal. Appl., 14 (1980), 142--144

\bibitem{Mol01}
Molchanov, V.F.,
{\it Unpublished notes on indefinite Plancherel formula,} 1980.

\bibitem{Mol1}
Molchanov, V.F.,
{\it Quantization on para-Hermitian symmetric
spaces}, Amer. Math. Soc. Transl., 1996, 175, 81 - 95.

\bibitem{Mol2}
Molchanov, V.F.,
{\it Separation of spectra for hyperboloids.}
Funkts. Analiz i Prilozh., 31,3 (1997), 35--43;
English translation in Funct. Anal. Appl. (1997)


\bibitem{ML}
Myller-Lebedeff, W.,
{\it Die Theorie der Integralgleichungen in
Awendung auf einge Reihenentwicklungen,}
 Math.Ann 64 (1907), 388--416

\bibitem{Ner}   
 Neretin, Yu.A.,
  {\it  On discrete occurrence of complementary
series representations in tensor products
 of unitary representations.}
Funct. Anal. Appl.., 20, 79-80 (1986)(Russian);
 English translation
in Funct.Anal.Appl.,20,68--70

\bibitem{Ner1}       
Neretin, Yu.A.,
 {\it Extension of representations of
 classical groups
to representations of categories%
.} Algebra i analiz, t.3, No 1, 176--202
(1991); English translation in St.Petersburg Math.J.,
Vol. 3(1992).

\bibitem{Ner2}             
 Neretin, Yu.A.,
   {\it Categories
  of symmetries and infinite-dimensional groups}
Oxford University Press (1996);
 Russian edition: Moscow, URSS(1998).



\bibitem{Ner3}         
Neretin, Yu.A.,
 {\it Boundary values
 of holomorphic functions and some
 spectral problems for unitary representations.}
In collection
 {\it Positivity in Lie groups: open problems}, 221--248,
  Hilgert J.,
J.D.Lawson, K.-H.Neeb, Vinberg E.B. (eds.),
Walter de Gruyter, Berlin (1998).

\bibitem{Ner4}      
Neretin, Yu.A.,
 {\it Restriction of  functions
 holomorphic in a domain
to a curve lying on the boundary and
discrete ${\rm SL}_2(\R)$-spectra.}
Izv. Ross. Akad. Nauk, Ser. mat., 62, 3(1998), 67--86;
English translation in Izv. Math., 62 (1998), 493--513.

\bibitem{pse}
Neretin, Yu.A.
{\it Pseudoriemannian symmetric spaces: one type realizations
and open embeddings to Grassmannians.},
Zapiski Nauchn. Semin. POMI, v. 256, 145--167;
English translation, to appear J.Math.Sci, New York;
preprint version is available via {\tt arXiv.org/abs/math/9905014}

\bibitem{Ner5}      
Neretin, Yu.A.,
 {\it Matrix analogs of $\B$-function
 and Plancherel formula for Berezin kernel representations,}
Mat. Sbornik, 191 (2000), 67--100 (Russian);
Preprint version  is available
via {\tt http://arXiv.org/abs/math/9905045}

\bibitem{Ner6}       
Neretin Yu.A.,
 {\it On separation of spectra
 in analysis of Berezin kernels}
 Funkts. Anal. Prilozh.,34 (2000), 3 (Russian);
Preprint version is available via
{\tt http://arXiv.org/abs/math/9906075}

\bibitem{Ner7}     
Neretin, Yu.A.,
{\it Plancherel formula for Berezin deformation
 of $L^2$ on Riemannian
symmetric space}, Preprint
 {\tt  http://arXiv.org/abs/math/9911020}
  (to appear J. Funct. Anal.)

\bibitem{Ner8} 
Neretin, Yu.A.,
 {\it Index hypergeometric transform
 and imitation of
analysis of Berezin kernels on hyperbolic spaces,}
Sbornik Math., to appear.

\bibitem{Ner9} 
Neretin, Yu.A.,
{\it Hua type integrals over unitary groups and over
inverse limits of unitary groups.}
Preprint  {\tt http://arXiv.org/abs/math/0010014}



\bibitem{NO}
Neretin, Yu.A., Olshanskii, G.I.,
{\it  Boundary values of
 holomorphic functions ,
 singular unitary representations of groups
$O(p,q)$  and their limits as  $q\to\infty$.}
Zapiski nauchn. semin. POMI RAN 223, 9--91(1995);
 English translation:
J.Math.Sci., New York, 87, 6 (1997), 3983--4035.





\bibitem{OO} \'Olafsson, G., \O rsted, B.,
{\it Bargmann transform for symmetric spaces,} in
{\it Lie theory and its applications in physics.},
 eds. Doebner H.,
Dobrev V.K., Hilgert J., World Scientific (1996), 3--14.









\bibitem{hardy} 
Olshanskii, G.I.,
{\it Complex Lie semigroups, Hardy spaces,
and Gelfand--Gindikin program.}
in {\it Group theory and cohomological algebra.}
Edition of Yaroslavl State University, 85--98
(1982)(Russian);
English translation in
Diff. Geometry and Appl., 1 (1991), 297--308.


\bibitem{Ols1}  
Olshanskii, G.I.,
 {\it Irreducible unitary representations
of the groups $\U(p,q)$ which admits the pass to the limit as
$q\to\infty$.} Zap. Nauchn. Semin. LOMI., 172, 114--120 (1989)
(Russian); English translation in J.Sov.Math., 59,
1102--1107(1992).

\bibitem{non}
Olshanskii, G.I.,
 {\it  Unpublished notes on inverse limits of symmetric spaces.}


\bibitem{OZ} \O rsted B., Zhang G.,
 {\it Tensor products of analytic continuations
of discrete series}. Can.\ J.\ Math.,49, N 6, 1224--1241(1997)


\bibitem{Pic}
Pickrell, D.,
{\it Measures on infinite-dimensional Grassmann
manifold.} J. Funct. Anal., 70 (1987), 323--356.


\bibitem{Puk}
Pukanszky, L., {\it On the Kronecker products of irreducible
 unitary representations  of the
 $2\times2$ real unimodular  group.}
  Trans. Amer. Math. Soc., 100
(1961), 116--152

\bibitem{RS}
Reed, M., Simon,B.,
{\it Methods of modern mathematical physics.}
V.1-2; Springer, 1972, 1975.

\bibitem{Rep}
Repka J., {\it Tensor products of holomorphic discrete series.}
Canad. J. Math., 31(1979), 836--844



\bibitem{Sie} Siegel, C.L.,
{\it Uber die analytische Theorie der quadratischen Formen,}
Ann. Math. 36 (1935), 527-606



\bibitem{Sch}
{ Schoenberg, I.J.,} {\it Metric spaces and positive
 definite functions.}
Trans. Amer. Math. Soc., 44 (1938), 522--536


\bibitem{Ter}
Terras  A., {\it Harmonic analysis on symmetric
spaces and applications, v. 2.}, Springer, 1988

\bibitem{UU}
Unterberger, A., Upmeier, H., {\it The Berezin transform and
 invariant differential operators}.
Comm.Math.Phys.,164, 563--597(1994)

 \bibitem{RV}Vergne, M., Rossi, H,
{\it Analytic continuations of holomorphic discrete series
of semisimple Lie groups.} Acta Math., 136, N1-2, 1-59 (1976)


\bibitem{VGG}
Vershik, A.M., Gelfand, I.M., Graev, M.I.,
{\it Representations of $SL_2(R)$, where $R$ is function ring.}
Uspehi Mat. Nauk 28 (1973), No 5,83--128(Russian);
English translation: Russian Math. Surveys 28 (1973), No 5, 87--132

\bibitem{Vla}
Vladimirov, V.S.,
{\it   Generalized functions in mathematical physics.}
Nauka, Moscow , 1976 (Russian)
English translation: Mir, Moscow, 1979

\bibitem{Wal} Wallach, N.R.,
{\it Analytic  continuation of discrete series.}
Trans. Amer. Math. Soc.,
251, 19-37 (1979)

\bibitem{Wey}
Weyl, H.,
{\it Uber gewonliche lineare Differentialgleichungen
mis singularen Stellen und ihre Eigenfunktionen {\rm 2} Note}.
Nachr. Konig. Gess. Wissen. Gottingen, Math.-Phys.,
1910, 442-467; Reprinted in
Weyl, H., {\it Gessamelte Abhandlungen}, Bd. 1, 224--247,
Springer, 1968.

\bibitem{Wil}
Wilson, J.A,
{\it Some hypergeometric polynomials.}
SIAM J. Math. Anal., 11, 690--701.


\bibitem{Zha} Zhang, Genkaj,
{\it Berezin transform on line bundles over bounded symmetric domains}
Preprint, 1999.

 \bibitem{Zhe}
 \v Zelobenko, D.P.,
 { \it Compact Lie groups and their representations.}
 Nauka, Moscow, 1970 (Russian); English translation
 Amer. Math. Soc.,  1973.



\end{thebibliography}
\end{document}